\documentclass[a4paper,titlepage]{report}

\usepackage{latexsym}
\usepackage{amssymb}
\usepackage{amsfonts}
\usepackage{amsmath}
\usepackage{amsthm}
\usepackage{mathtools}
\usepackage{comment}
\usepackage{booktabs}
\usepackage[T1]{fontenc} 
\usepackage[utf8]{inputenc} 
\usepackage[english]{babel}
\usepackage[shortlabels]{enumitem}
\usepackage[autostyle,italian=guillemets]{csquotes}
\usepackage[style=numeric]{biblatex}
\usepackage{todonotes}
\usepackage{geometry}
\usepackage{frontespizio}
\usepackage{titlesec}
\usepackage{url}
\usepackage{indentfirst}

\usepackage{tikz}

\theoremstyle{definition}
\newtheorem{definition}{Definition}[section]
\theoremstyle{plain}
\newtheorem{theorem}[definition]{Theorem}
\newtheorem{lemma}[definition]{Lemma}
\newtheorem{proposition}[definition]{Proposition}
\newtheorem{corollary}[definition]{Corollary}
\theoremstyle{remark}
\newtheorem{remark}[definition]{Remark}
\newtheorem{example}[definition]{Example}

\newcommand{\numberset}{\mathbb}
\newcommand{\N}{\numberset{N}}
\newcommand{\Z}{\numberset{Z}}
\newcommand{\ua}{u_{\alpha}}
\newcommand{\Ua}{U_{\alpha}}
\newcommand{\pp}{PP_n}
\newcommand{\pf}{PF_n}
\newcommand{\pr}{PR_n}
\newcommand{\nap}{PF_{n,1}}
\newcommand{\knap}{PF_{n,k}}
\newcommand{\A}{\abs{\alpha}}
\newcommand{\uaJ}{u_{\alpha_{|_J}}}
\newcommand{\AJ}{\abs{\alpha_{|_J}}}
\newcommand{\tnk}{\vartheta_{n,k}}
\newcommand{\Teq}{\Theta_{n,k}^=}
\newcommand{\Tnk}{\Theta_{n,k}^{\le}}
\newcommand{\unk}{\upsilon_{n,k}}
\newcommand{\R}{\mathcal{R}}

\DeclarePairedDelimiter{\abs}{\lvert}{\rvert}

\usepackage[hyperindex,breaklinks]{hyperref}
\usepackage{xurl}

\hypersetup{
	breaklinks=true,   % splits links across lines
	colorlinks=true,   % displays links as colored text instead of blocks
	pdfusetitle=true,  % \title and \author values into pdf metadata
}

\begin{document}

%\begin{frontespizio}
	\Preambolo{\renewcommand{\fronttitlefont}{\fontsize{16}{20}\scshape}}
	\Universita{Firenze}
	%\Logo{Università_Firenze_Logo}
	\Facolta{Scienze Matematiche, Fisiche e Naturali}
	\Corso[Laurea Magistrale]{Matematica}
	\Titoletto{Tesi di Laurea Magistrale}
	\Titolo{On the combinatorics of $k$-Naples parking functions and parking strategies}
	\Candidato{Francesco Verciani}
	\Relatore{Luca Ferrari}
	\Annoaccademico{2021-2022}
%\end{frontespizio}

\hypersetup{linkcolor=black}
\tableofcontents
\setcounter{page}{1}

\hypersetup{linkcolor=red}

\chapter*{Introduction}

\addcontentsline{toc}{chapter}{Introduction}

In 1966, Alan Konheim and Benjamin Weiss published the paper \emph{An Occupancy Discipline and Applications} \cite{doi:10.1137/0114101}, where they presented a mathematical problem concerning the speed of storing and retrieval of records in a computer filing system. Further ahead in the article, the same problem is presented in terms of $n$ cars parking in $n$ labelled spots on a one-way street. Each car has an associated preference - a favourite spot, if you will - and as a first attempt, it tries to park in that spot. If it is not available, the car then proceeds to drive forwards and park in the first available spot, if there are any. This is the usual way to describe the problem nowadays, and lists of preferences such that all cars are able to park are called \emph{parking functions}.

In the following years, parking functions have been generalized in various ways. One such way to generalize them leads to the notion of \emph{$k$-Naples parking functions} \cite{https://doi.org/10.48550/arxiv.1908.07658}, where cars are allowed backward movement. The ability to drive backwards complicates matters: in particular, while all rearrangements of a parking function are still parking functions, the same does not hold for $k$-Naples parking functions. Thus, the way the cars are listed needs to be taken into account.

This thesis is structured as follows. We describe the parking problem in detail in Chapter \ref{chapter1}, where we also introduce the relevant notations to our purposes. In Chapter \ref{chaptercharacterization}, we introduce a subclass of parking preferences, which we call \emph{complete} parking preferences. We then propose a characterization of $k$-Naples parking functions, in terms of subsequences that are complete $k$-Naples parking functions. Along the way, we are also able to determine which $k$-Naples parking functions are such that any rearrangement is also a $k$-Naples parking function. We extend our results in cases where we consider less cars than the total number available spots, or where the number of spots and the number of cars are both infinite. We provide some enumerative results regarding complete $k$-Naples parking functions and permutation invariant $k$-Naples parking functions in Chapter \ref{chapterenumeration}.

In \cite{https://doi.org/10.48550/arxiv.2001.04817}, the author suggests a generalization of the $1$-Naples parking problem, where only even-numbered cars (or only odd-numbered cars) are able to drive backwards. A natural further generalisation is then to consider different parking rules for each car. We study such a situation in Chapter \ref{chapterstrategies}, and define a given list of rules to be a \emph{parking strategy} if it allows all cars to park. Given a preference list, we find parking strategies that minimize the number of steps taken by the cars, as well as parking strategies where the minimum number of cars is allowed to drive backwards. Finally, we completely characterize the case where the chosen parking rules are restricted to the standard rule and the $1$-Naples parking rule.

\chapter{Parking functions}
\label{chapter1}
\section{Classical parking functions}

We start by giving a formal definition of the parking problem.

\begin{definition}%{\rm{\bf (Parking functions)}}
	Let $n$ be a positive integer, and let $[n]\coloneqq\{1,2,\dots,n\}$. Consider $n$ cars, denoted by $c_1,c_2,\dots,c_n$, that attempt to park, one at a time, in order, on a one-way street with $n$ parking spots, labelled $1,2,\dots,n$. Let $\alpha\coloneqq(a_1,a_2,\dots,a_n)$, where $a_i\in[n]$ is the \textit{parking preference} of car $c_i$. Each car $c_i$ behaves according to the following rule: it drives up to parking spot $a_i$, and parks there if it is not occupied; otherwise, it proceeds forwards and parks in the first available spot (i.e. the minimum $j>a_i$ that is unoccupied), and is unable to park if there are none. We will refer to this rule as the  \emph{standard parking rule}.\\
	A \emph{parking function} is defined as a parking preference that allows all cars to park. We denote by $\pp\coloneqq[n]^n$ the set of all parking preferences of length $n$, and by $\pf$ the set of parking functions of length $n$. We also set
	\begin{equation}
		\abs{\alpha}_i\coloneqq\abs*{\{j\in[n]\mid a_j=i\}},
	\end{equation}
	which is the number of cars that have parking preference $i$.
\end{definition}

\begin{comment}
\begin{example}
	Let $\alpha=(3,1,3,4,1)$\\
	\begin{tikzpicture}[line cap=round,line join=round,>=triangle 45,x=1cm,y=1cm]
		\clip(0.8105853041362491,-1.356082851905133) rectangle (5.230625401459856,1.325336444768069);
		\draw [line width=2pt] (1,0)-- (5,0);
		\draw (1.0005789174032005,-0.4419626371301777) node[anchor=north west] {1};
		\draw (2.000733975921681,-0.4419626371301777) node[anchor=north west] {2};
		\draw (3.0008890344401618,-0.4419626371301777) node[anchor=north west] {3};
		\draw (4.001044092958642,-0.4419626371301777) node[anchor=north west] {4};
		\draw (5.001199151477122,-0.4419626371301777) node[anchor=north west] {5};
		\draw (3.0008890344401618,0.5581924213883027) node[anchor=north west] {$c_1$};
		\draw (1.0005789174032005,0.5581924213883027) node[anchor=north west] {$c_2$};
		\draw (4.001044092958642,0.5581924213883027) node[anchor=north west] {$c_3$};
		\draw (5.001199151477122,0.5581924213883027) node[anchor=north west] {$c_4$};
		\draw (2.000733975921681,0.5581924213883027) node[anchor=north west] {$c_5$};
		\begin{scriptsize}
			\draw [fill=xdxdff] (2,0) circle (2.5pt);
			\draw [fill=xdxdff] (3,0) circle (2.5pt);
			\draw [fill=xdxdff] (4,0) circle (2.5pt);
			\draw [fill=xdxdff] (5,0) circle (2.5pt);
			\draw [fill=xdxdff] (1,0) circle (2.5pt);
		\end{scriptsize}
	\end{tikzpicture}
dyfdfjfdgjfdgf
\end{example}
\end{comment}

An elegant characterization of parking functions is given by the following well known result.

\begin{theorem}%{\rm{\bf(Characterization of parking functions)}}\\
	\label{standard pf}
	Let $\alpha=(a_1,a_2,\dots,a_n)\in\pp$ be a parking preference of length $n$, and $b_1\le b_2\le\dots\le b_n$ be its rearrangement in increasing order. Then $\alpha$ is a parking function if and only if $b_j\le j$ for all $j\in[n]$. 
\end{theorem}

\begin{proof}
	We can rewrite the desired condition as
	\begin{equation}
		\label{charpf}
		\sum_{i=1}^j \A_i\ge j \qquad \text{for all $j\in[n]$}.
	\end{equation}
	Suppose $\alpha\in\pf$, and let $j\in[n]$. It is clear that, if car $c_i$ has parking preference $a_i$, then it certainly cannot park in any spot strictly smaller than $a_i$. That means that parking spots $1$ through $j$ can only be filled by cars with parking preference at most $j$; since $\alpha$ is a parking function, all of those spots are filled, which means there are at least $j$ cars with preference at most $j$, and we have \eqref{charpf}. \\
	Vice versa, if  \eqref{charpf} is true, suppose that $\alpha$ is not a parking function, then there exists $j\in[n]$ such that parking spot $j$ is empty after all cars have attempted to park. Consider $i\in[n]$ such that car $c_i$ has preference $a_i\le j$; car $c_i$ has to have parked somewhere in $[1,j-1]$, because otherwise it would have arrived at parking spot $j$ and, finding it empty, would have parked there. So all cars with parking preference in $[1,j]$ park in $[1,j-1]$, meaning there can be at most $(j-1)$ cars with preference at most $j$, i.e. $\sum_{i=1}^j \A_i\le (j-1)<j$, and we have a contradiction.
\end{proof}

\begin{remark}
	Theorem \ref{charpf} provides an intuitive way to characterize parking functions: $\alpha\in\pp$ is a parking function if and only if there is at least $1$ car with preference $1$, at least $2$ cars with preference in $[1,2]$, and so on. Of particular note, changing the order of the preferences in $\alpha$ does not influence whether $\alpha$ is a parking function, as the only thing that matters is the number of cars with certain preferences. Thus, we say that parking functions are \emph{permutation invariant}. That is, let $\alpha=(a_1,\dots,a_n)\in\pp$ and $\sigma\in S_n$, where $S_n$ is the set of permutations of length $n$, and denote
	\[
	\sigma(\alpha)=(a_{\sigma(1)}\dots,a_{\sigma(n)}).
	\]
	Then $\alpha$ is a parking function if and only if $\sigma(\alpha)$ is a parking function for all $\sigma\in S_n$.\\
\end{remark}

\begin{theorem}
	\label{numpf}
The number of parking functions of length $n$ is $(n+1)^{n-1}$.
\end{theorem}

This result is proved by Konheim and Weiss in \cite{doi:10.1137/0114101}, but there exists a number of alternative proofs. We will outline two of them that we find especially noteworthy.
The first one involves modifying the parking problem, and leads to Theorem \ref{standard pf count Pollak}. The second one relies on a recursive formula found in Theorem \ref{recursivepfcount}.\\

For the first proof, illustrated by J. Riordan in \cite{RIORDAN1969408} and due to H. O. Pollak, we start by slightly modifying the classical parking problem. Suppose that the $n$ spots are arranged in a circle, instead of a one way street: in other words, suppose that, in the original problem, if a car does not find any available spot from its preference to spot $n$, then it starts checking from spot $1$ onwards. Clearly, in this version of the problem each car is eventually able to park, as long as not all spots are already occupied.

\begin{lemma}
	\label{circularcount}
	Consider the circular parking problem described as above, with $n$ parking spots and $m<n$ cars. Then the number of parking preferences ${\alpha=\{a_1,\dots,a_m\}}$ such that no car occupies spot $n$ is $(n-m)n^{m-1}$.
\end{lemma}

\begin{proof}
	Consider $\alpha=(a_1,\dots,a_m)\in[1,n]^m$. In the circular problem, there will be some $n-m$ spots left unoccupied. For each $j\in[n]$, consider the map $\overline{\tau}_j$ defined as
	\[
	\overline{\tau}_j(\alpha)=(\overline{a_1-j},\dots,\overline{a_m-j}),
	\]
	where $\overline{i}\in[n]$ is the representative of $i$ modulo $n$. Due to the cyclical nature of the modified problem, clearly if $\alpha$ leaves spots $i_1,\dots,i_{n-m}$ unoccupied, then $\overline{\tau}_j(\alpha)$ leaves spots $\overline{i_1-j},\dots,\overline{i_m-j}$ unoccupied. Thus, there are exactly $n-m$ distinct indices $j=j_1,\dots,j_{n-m}$ such that $\overline{\tau}_j(\alpha)$ leaves spot $n$ unoccupied.\\
	The maps $\overline{\tau}_j$, for $j\in[n]$, induce an equivalence in $[1,n]^l$ such that $\alpha_1$ and $\alpha_2$ are equivalent if and only if $\alpha_1=\overline{\tau}_j(\alpha_2)$ for some $j\in[n]$. Furthermore, each equivalence class contains exactly $n$ elements and, by our previous considerations, $n-m$ of those leave spot $n$ unoccupied. Thus, a fraction of $\frac{n-m}{n}$ of the total $n^m$ ways to choose $\alpha$ is such that spot $n$ remains free, so the number of desired elements is $\frac{n-m}{n}n^m=(n-m)n^{m-1}$.
\end{proof}

\begin{theorem}
	\label{standard pf count Pollak}
	Consider the standard parking problem with $n$ parking spots and $1\le m\le n$ cars. Then the number of parking preferences $\alpha=(a_1,\dots,a_m)\in[1,n]^m$ such that all cars are able to park is $(n-m+1)(n+1)^{m-1}$. In particular, $\abs{\pf}=(n+1)^{n-1}$.
\end{theorem}

\begin{proof}
	Consider the circular version of the parking problem with $n+1$ spots. Note that a parking preference $\alpha\in[1,n]^m$ is such that all cars are able to park in the standard problem if and only if it is such that it leaves spot $n+1$ unoccupied in the modified problem. This  follows from the fact that, if spot $n+1$ is never reached, then the two versions of the problem are equivalent. Thus, by Lemma \ref{circularcount}, the number of parking preferences such that $m$ cars park in $n$ spots is $(n-m+1)(n+1)^{m-1}$.
\end{proof}

Using the binomial expansion, we get $\abs{\pf}=(n+1)^{n-1}=\sum_{i=0}^{n-1} \binom{n-1}{i}n^{n-i+1}=\sum_{m=1}^{n} \binom{n-1}{m-1}n^{n-m}$. It is interesting to note that the terms of this sum actually correspond to a specific partition of $\pf$.

\begin{lemma}
	\label{1preference}
	Let $n\ge m\ge 1$. Then the number of parking functions of length $n$ such that $m$ cars have preference $1$ is $\binom{n-1}{m-1}n^{n-m}$.
	%\begin{equation}
	%	\label{count0}
	%	\vartheta_{n,0}(m,1)=\binom{n-1}{m-1}n^{n-m}.
	%\end{equation}
\end{lemma}

\begin{proof}
	Clearly, there is only $1=\binom{n-1}{n-1}n^{n-n}$ parking function such that all $n$ cars have preference $1$. Given $l\in[1,n-1]$, we now evaluate how many parking functions have $n-l$ cars with preference $1$.\\
	Let $\beta=(b_1,\dots,b_l)\in[1,n]^l$ be a parking preference such that $l$ cars leave spot $1$ free in the circular problem with $n$ spots. Note that such preferences are counted as in Lemma \ref{circularcount}, since it does not actually matter which spot is the one to leave unoccupied, so there is a total of $(n-l)n^{l-1}$ such preferences. Consider the standard parking problem with preference $\alpha=(b_1,\dots,b_l,1,\dots,1)\in\pp$. For all $i\le l$, car $c_i$ is able to park: this can be shown by observing that, in the modified problem, no cars reach spot $n$ and start over from spot $1$, since spot $1$ has remained unoccupied. Moreover, it is clear that all cars $c_j$ ($l+1\le j\le n$) park somewhere, so $\alpha$ is a parking function. Finally, since $\beta$ left spot $1$ unoccupied in the modified problem, trivially $b_i\ne1$ for all $i\le l$, so $\A_1=n-l$. \\
	As we have already observed, as a consequence of Theorem \ref{charpf}, since $\alpha$ is a parking function, it is permutation invariant. Thus, if $\alpha'\in\pp$ is a parking preference that is a rearrangement of $\alpha$ where the elements $b_1,\dots,b_l$ remain in the same relative order, then $\alpha'$ is also a parking function.
	Conversely, it is easy to check that if $\alpha\in\pf$ is such that $\A_1=n-l$, then the other elements correspond to a sequence $\beta=(b_1,\dots,b_l)\in[1,n]^l$ as above.\\
	Thus, a parking function $\alpha\in\pf$ with $n-l$ preferences equal to $1$ is uniquely determined by choosing $\beta$ and $l$ indices $j_1<\dots,j_l$, so that $a_{j_i}=b_i$ for all $i=1,\dots,l$, and $a_j=1$ otherwise. Therefore, the number of parking functions such that $\A_1=n-l$ is
	\[
	(n-l)n^{l-1}\binom{n}{l}=\binom{n-1}{n-l-1}n^l.
	\] 
	We conclude by setting $m=n-l$.
\end{proof}

We now proceed with the second proof of Theorem \ref{numpf}, while illustrating 
a recursive aspect of parking functions highlighted in \cite{https://doi.org/10.48550/arxiv.1908.07658}.

\begin{theorem}
	\label{recursivepfcount}
	The number of parking functions of length $n$ follows the recursive relation
	\begin{equation}
		\abs{\pf}=\sum_{j=0}^{n-1} \binom{n-1}{j}(j+1)\abs{PF_{j}}\abs{PF_{n-j-1}}.
	\end{equation}
\end{theorem}

\begin{proof} Consider $\alpha\in PF_{n}$, $n\ge2$, and let $j+1\in[n]$. Suppose that the first $n-1$ cars occupy all spots except for spot $j+1$, then clearly $a_{n}\in[1,j+1]$ since $\alpha$ is a parking function.
Moreover, the $j$ cars occupying spots $[1,j]$ correspond to a parking function of length $j$, and the same is true for the cars occupying spots $[j+2,n]$. Accounting for all ways to choose $j$ cars among the first $n-1$, and for all values of $j$, we get
\[
\abs{\pf}=\sum_{j=0}^{n-1} \binom{n-1}{j}(j+1)\abs{PF_{j}}\abs{PF_{n-j-1}}.
\]
\end{proof}

Now, there is clearly only $1$ parking function of length $1$, and we can consider the empty parking preference as being the only parking function of length $0$. Suppose that, for all $i\le n-1$, $PF_{i}=(i+1)^{i-1}$. We recall Abel's generalization of the binomial formula (see, for example, \cite[18]{riordan1968combinatorial}):
\begin{equation}
	\label{Abel}
(z+w+m)^m=\sum_{j=0}^m\binom{m}{j} w (w+m-j)^{m-j-1} (z+j)^j.
\end{equation}
Setting $m=n-1$, $w=1$ and $z=1$ gives
\[
\abs{\pf}=\sum_{j=0}^{n-1} \binom{n-1}{j}(j+1)(j+1)^{j-1}(n-j)^{n-j-2}=(n+1)^{n-1}.
\]
\\

By Cayley's formula, $(n+1)^{n-1}$ is also the number of labelled rooted forests on $n$ vertices, or the number of labelled trees on $n+1$ vertices, and it is known that this also corresponds to the number of acyclic maps $\beta\colon[n]\to[n]$. Multiple bijections between these classes of objects have been found \cite{schutzenberger1968enumeration,kreweras1980famille,foata1974mappings,https://doi.org/10.48550/arxiv.0810.0427}.\\

\section{Naples parking functions}

The standard parking problem has been generalized in various ways. Probability is the most common way to shake things up, and indeed one such generalization was already presented when parking functions were first introduced in \cite{doi:10.1137/0114101}, involving choosing the parking preferences in some random way. Other generalizations include: having the parking preferences avoid certain patterns \cite{https://doi.org/10.48550/arxiv.2209.04068}, changing the size of the cars \cite{Richard_Ehrenborg_2016},  and even arranging the parking spots on rooted labelled trees \cite{https://doi.org/10.48550/arxiv.1504.04972}. The article \cite{https://doi.org/10.48550/arxiv.2001.04817} gives a brief overview of some further ways to generalize the parking problem.\\

From this point on, we are going to study \emph{$k$-Naples parking functions}, a generalization allowing cars to move backwards, introduced by Baumgardner \cite{baumgardner2019naples} and later further generalized in \cite{https://doi.org/10.48550/arxiv.1908.07658}.
\begin{definition}%{\rm{\bf ($k$-Naples parking functions)}}
	Modify the standard parking rule as follows: if car $c_j$ drives up to parking spot $a_j$ and finds it occupied, it drives backwards and checks spot $a_j -1$ (if it exists), filling it if it is free, before proceeding forwards like the standard parking rule. This is called the \emph{1-Naples parking rule}, or simply \emph{Naples parking rule} if there is no risk of confusion. Given $k\in[n]$, we can further modify the rule by allowing cars to check up to $k$ spots prior to their preference, again filling the first one available, before proceeding forwards if necessary. This is called the \emph{$k$-Naples parking rule}.\\
	A \emph{$k$-Naples parking function} of length $n$ is a parking preference $\alpha=(a_1,a_2,\dots,a_n)$ such that all cars $c_1,c_2,\dots,c_n$ are able to park following the $k$-Naples parking rule; we denote by $\knap$ the set of $k$-Naples parking functions.
\end{definition}

Note that, while Theorem \ref{standard pf} guarantees that parking functions are permutation invariant, this is generally not the case for k-Naples parking functions: for example, $(3,3,2)\in\nap$, but $(2,3,3)\notin\nap$.\\

It is clear that if $k\ge n-1$, then each car is able to check all possible spots, so $\abs{PF_{n,n}}=\abs{PF_{n,n-1}}=n^n$. Following a similar argument to the recursion previously shown for classical parking functions, the authors in \cite{https://doi.org/10.48550/arxiv.1908.07658} provide the following formula for enumerating $k$-Naples parking functions:

\begin{theorem}
Let $n\ge k\ge1$. Then:
\begin{equation}
	\label{knap recursive}
	\abs{PF_{n+1,k}}=\sum_{i=0}^n\binom{n}{i}\min((i+1)+k,n+1)\abs{PF_{i,k}}(n-i+1)^{n-i-1}.
\end{equation}
\end{theorem}

\begin{proof}
	Let $\alpha=\{a_1,\dots,a_{n+1}\}\in PF_{n+1,k}$, and $i\in[0,n]$. Suppose that the first $n$ cars occupy all spots except for spot $(i+1)$. The $i$ cars occupying spots $[1,i]$ correspond to a $k$-Naples parking function of length $i$. \\
	On the other hand, cars occupying spots $[i+2,n+1]$ correspond to a $k$-Naples parking function following the additional condition that no car tries to  drive backward from spot $1$: these are sometimes called \emph{contained} parking functions. Observe that a sequence $\beta=(b_1,\dots,b_m)\in[m]^m$ is a contained $k$-Naples parking function if and only if, when considering the circular problem with $m+1$ spots, spot $(m+1)$ is left unoccupied: this is because in that case, spot $m+1$ is occupied if and only if it is reached either by a car driving forwards from spot $m$ (and thus $\beta$ would not be a parking function), or by a car driving backwards from spot $1$ (which contradicts the additional condition we imposed). Thus, following the same argument used in Lemma \ref{circularcount}, contained $k$-Naples parking functions are equinumerous to standard parking functions, and consequently there are $\abs{PF_{n-i}}=(n-i+1)^{n-i-1}$ possible parking preferences for cars occupying spots $[i+2,n+1]$.\\
	Finally, car $c_{n+1}$ can reach spot $(i+1)$ either by simply driving forwards, if $a_{n+1}\in[1,i+1]$, or by driving backwards for up to $k$ spots, if $a_{n+1}\in[i+2,i+1+k]\cap[1,n]$. Accounting for all ways to choose $i$ cars among the first $n$, and for all values of $i$, we get \eqref{knap recursive}. 
\end{proof}

An explicit formula for the number of $k$-Naples parking functions has also been found \cite{https://doi.org/10.48550/arxiv.2009.01124}.\\

In \cite{carvalho:online}, a thorough examination of ascending and descending $k$-Naples parking functions is conducted, finding bijections with combinatorial objects enumerated by Catalan numbers. A characterization of permutation invariant $k$-Naples parking functions is also established, and we will provide a similar result in the next chapter (Corollary \ref{perminv}).

\section{Notation and preliminaries}
In this section, we introduce a new notation regarding the number of cars that have specific preferences, as well as its relevant arithmetical properties.

\begin{definition}
Let $\alpha=(a_1,a_2,\dots,a_n)\in\pp$, we define for all $j\in[n]$:
\begin{equation}
\label{ua}
\ua(j)\coloneqq \sum_{i=j}^n \A_i -(n-j+1).
\end{equation}
In other words, $\ua(j)$ counts the excess of cars having parking preference in $[j,n]$ (i.e. $\sum_{i=j}^n \A_i$) with respect to the number of available parking spots in $[j,n]$ (which is $n-j+1$). We also define
\begin{equation}
\Ua\coloneqq\{j\in[n] \mid \ua(j)\ge1\}.
\end{equation}
Clearly $\ua(1)=0$, and consequently $1\notin\Ua$.\\
\end{definition}

The next proposition records some basic properties of the quantity $\ua(j)$ that will be very useful in the sequel.

\begin{proposition}
Let $\alpha\in\pp$, then for all $j\in[n]$, and $j_1,j_2\in[n]$ such that $j_1<j_2$:
\begin{subequations}
\begin{align}
\label{propua} &\text{(i)} \quad  \ua(j)=j-1- \sum_{i=1}^{j-1}\A_i;\\
\label{propua2}&\text{(ii)}  \quad\ua(j_2)-\ua(j_1)=(j_2-j_1)-\sum_{i=j_1}^{j_2-1}\A_i;\\
\label{propua3}&\text{(iii)}\quad \ua(j)=\ua(j+1)+\A_j-1\qquad \text{if $j<n$};\\
\label{propua4}&\text{(iv)}\quad \ua(j)<j.
\end{align}
\end{subequations}

\end{proposition}

\begin{proof}
Observe that $n= \sum_{i=1}^{j-1}\A_i+ \sum_{i=j}^{n}\A_i$; substituting in \eqref{ua} we get \emph{(i)}. Furthermore, exploiting \eqref{ua} we get
\[
\ua(j_2)-\ua(j_1)=\sum_{i=j_2}^n \A_i -(n-j_2+1) - \sum_{i=j_1}^n \A_i +(n-j_1+1)=(j_2-j_1)-\sum_{i=j_1}^{j_2-1}\A_i.
\]
Finally, \emph{(iii)} is just a special case of \emph{(ii)}, with $j_1=j$ and $j_2=j+1$, and \emph{(iv)} clearly follows from \emph{(i)}.
\end{proof}

\begin{definition}
A \emph{maximal interval} of a finite subset $A$ of $\N$ is an interval $I\subseteq A$ which is not strictly contained in any other interval $J\subseteq A$. In particular, if $[p,q]\subseteq A$ is a maximal interval of $A$, then $(p-1)\notin A$ and $(q+1)\notin A$.
\end{definition}

\begin{proposition}
	\label{propUa}
Let $\alpha\in\pp$, and $[p,q]\subseteq\Ua$ a maximal interval of $\Ua$. Then:
\begin{subequations}
\begin{align*}
&\text{(i)}\quad\! \ua(p)=1;    %\qquad
&&\text{(ii)}\quad\!  \ua(p-1)=0;\\
&\text{(iii)} \quad\! \A_{p-1}=0; %\qquad
&&\text{(iv)}\quad\!\A_q\ge2.
\end{align*}
\end{subequations}
\end{proposition}

\begin{proof}
We recall that $1\notin\Ua$, so $(p-1)\ge1$. Since $[p,q]\subseteq\Ua$ is a maximal interval, we have $(p-1)\notin\Ua$, so $\ua(p-1)\le0$, whereas $\ua(p)\ge1$. By \eqref{propua3} we get
\[
1\le\ua(p)=\ua(p-1)-\A_{p-1}+1\le1;
\]
\[
0\ge\ua(p-1)=\ua(p)+\A_{p-1}-1\ge1-1=0.
\]
Thus, $\ua(p)=1$ and $\ua(p-1)=0$, and we obtain
\[
0= \ua(p-1) = \ua(p)+\A_{p-1}-1= \A_{p-1}.
\]
Now, note that $\ua(q)\ge1$. Suppose that $q=n$, then by definition $1\le \ua(n)=\A_n-1$, so $\ua(q)=\ua(n)\ge2$. On the other hand if $q<n$, then $\ua(q+1)\le0$, and consequently
\[
1\le\ua(q)= \ua(q+1)+\A_q-1\le \A_q-1,
\]
and we get \emph{(iv)}.\\
\end{proof}

\begin{remark}
Using formula \eqref{propua}, the inequality $\sum_{i=1}^j\A_i\ge j$  is equivalent to $\ua(j+1)\le0$. Note that $\ua(1)\le0$ and $\sum_{i=1}^n\A_i\ge n$ are always satisfied. Thus,  condition \eqref{charpf} is equivalent to $\ua(j)\le0$ for all $j\in[n]$, and Theorem  \ref{standard pf} can be reformulated as:
\begin{equation}
\alpha\in\pf\quad \text{if and only if} \quad \Ua=\varnothing.
\end{equation}
\end{remark}

Our aim is to find an analogous characterization for \emph{k}-Naples parking functions, so that determining whether a parking preference is a $k$-Naples parking function can be done without actually checking if and where the cars are going to park. Note that, since $k$-Naples parking functions are generally \emph{not} permutation invariant, a characterization that solely relies on $\ua$ and $\Ua$ - which are computed independently of the order in a parking preference - cannot exist, so further considerations need to be made.\\
We start by proving the following necessary condition:

\begin{theorem}
	\label{necessk}
Let $k\ge1$ and $\alpha\in\knap$, then $\ua(j)\le k$ for all $j\in[n]$.
\end{theorem}

\begin{proof}
Suppose that there exists $j\in[n]$ such that $\ua(j)>k$. Note that for all $i$, $\ua(i)\le i-1$, and consequently $j\ge k+2$, in particular there exist at least $k+1$ parking spots prior to parking spot $j$. Observe that car $c_i$, following the \emph{k}-Naples parking rule, certainly cannot park in any parking spot that is strictly smaller than $a_i-k$. Moreover, since $\alpha\in\knap$, all cars with preference in $[j,n]$ will park somewhere in $[j-k,n]$. However, $[j-k,n]$ contains $(n-j+k+1)$ spots, whereas the number of cars with preference in $[j,n]$ is
\begin{equation}
\sum_{i=j}^n \A_i > (n-j+1)+k,
\end{equation}
where we used $\ua(j)>k$ and \eqref{ua}, thus we have a contradiction.
\end{proof}

The previous theorem expresses the fact that if $\alpha\in\knap$, then $k\ge \max_j \ua(j)$. This is obviously not a sufficient condition: for example, $\alpha=(2,3,3)$ is such that $\ua(j)=1$ for all $j\in[5]$, but it is not a $1$-Naples parking function. \\
Given $\alpha\in\pp$ and $k\ge \max_j \ua(j)$, the following theorem shows that there does  exist a rearrangement of $\alpha$ such that it becomes a $k$-Naples parking function: specifically, the most favourable arrangement for $\alpha$ is the descending order, meaning that such a parking preference is a $k$-Naples parking function if and only if $k\ge \max_j \ua(j)$.

\begin{theorem}
	\label{easydescend}
	Let $\alpha\in\pp$ such that $\max_j \ua(j)=M\ge0$. If $\alpha$ is in descending order, then $\alpha\in\knap$ for all $k\ge M$.
\end{theorem}

\begin{proof}
	Suppose that $\alpha\notin\knap$, then there exists $i\in[n]$ such that car $c_i$ that cannot park. Then all spots in $[a_i-k,n]$ have already been occupied. Since $\alpha$ is in descending order, spots $[a_i-k,n]$ are occupied by cars with preference in $[a_i,n]$. Thus, the number of cars with preference in $[a_i,n]$ is at least the number of spots in $[a_i-k,n]$, that is, $n-a_i+k+1$, plus one (car $c_i$ that does not park). As a result, we get
	\[
	\ua(a_i)=\sum_{j=a_i}^n\A_j-(n-a_i+1)\ge (n-a_i+k+2)-(n-a_i+1)=k+1,
	\]
	which contradicts the condition $k\ge \max_j \ua(j)$.
\end{proof}

\chapter{A characterization of $k$-Naples parking functions}
\label{chaptercharacterization}
\section{Complete parking preferences}

We now focus on a special subclass of parking preferences.
\begin{definition}
Let $n\ge2$, and $\alpha\in\pp$. We say that $\alpha$ is \emph{complete} if $\Ua=[2,n]$.
\end{definition}

\begin{remark}
As previously observed $1\notin\Ua$, so $\alpha$ is complete if $\Ua$ has maximum size. 
\end{remark}

\begin{definition}
	\label{outcomemap}
Let $\alpha=(a_1,\dots,a_n)\in\pp$, and assume all cars follow the the standard parking rule. We define the function $\psi_0$, which maps $\alpha$ to $\psi_0(\alpha)=(p_1,\dots,p_n)$, where for each $j\in[n]$:
\begin{itemize}
\item if car $c_j$ can park, then $p_j$ is the parking spot that car $c_j$ occupies;
\item if car $c_j$ cannot park, then $p_j=\infty$.
\end{itemize}
Given $\alpha\in\pp$, with abuse of notation we will write $\psi_0(c_j)=p_j$.
For $k\ge1$, we define functions $\psi_k$ analogously, assuming all cars follow the \emph{k}-Naples parking rule.
\end{definition}

\begin{example}
	Let $\alpha=(5,3,3,5,4)\in PP_5$. Computing $\ua(j)$ for $j\in[2,5]$ gives:
	\[
	\ua(5)=1;\quad \ua(4)=1; \quad \ua(3)=2;\quad \ua(2)=1.
	\]
	Thus $\alpha$ is complete; furthermore, it is a $3$-Naples parking function, since $\psi_3(\alpha)=(5,3,2,4,1)$. Observe that it is not $2$-Naples, since, assuming all cars follow the $2$-Naples parking rule, car $c_5$ is not able to park, i.e. $\psi_2(c_5)=\infty$. Note that its rearrangement $(5,5,4,3,3)$ is a $2$-Naples parking function, and that there is no rearrangement that yields a $1$-Naples parking function, since $\ua(3)=2$, and that would contradict Theorem \ref{necessk}.
\end{example}

\begin{remark}
The maps $\psi_k$ just defined are an extension of the \emph{outcome map} shown in \cite{colaric:online}, which was defined only for standard parking functions. It can be useful to define the outcome even if not all cars are able to park: specifically, setting $\psi_k(c_i)=\infty$ when a car does not park, we get that $\psi(c_i)>a_i$. This will be especially useful for the study of complete parking preferences.\\
\end{remark}

\begin{lemma}
	\label{spot1}
	Let $\alpha=(a_1,\dots,a_n)\in\pp$ be complete. If $\alpha\in\knap$, then $\psi_k(c_n)=1$.
\end{lemma}

\begin{proof}
	Suppose that $\alpha\in\knap$. Consider the first $n-1$ cars, and let $h$ be the only spot left unoccupied by those cars. If $h>1$, then spots $[1,h-1]$ are all occupied by cars with preference in $[1,h-1]$, since any car with preference at least $h$ parking in $[1,h-1]$ would have also checked spot $h$, which is unoccupied. Thus, the number of cars with preference in $[1,h-1]$ is at least $h-1$, and so $\ua(h)=h-1-\sum_{i=1}^{h-1}\A_i\le0$, which is absurd since $\alpha$ is complete. Therefore, $h=1$, so the last car will fill spot $1$.
\end{proof}

\begin{lemma}
	\label{menouno}
	Let $\alpha\in\pp$, $k\ge1$, and assume that all cars follow the $k$-Naples parking rule. Let $i,j\in[n]$, and suppose that car $c_i$ is such that $\psi_k(c_i)=j$, and has preference $a_i<j$. Then there are no cars with preference in $[j,n]$ that park in any spot in $[1,j]$.
	Moreover, if $\alpha\in\knap$, then $\ua(j)\le-1$.
\end{lemma}

\begin{comment}
\begin{proof}
	Suppose $a_i\le k$, then when car $c_i$ attempts to park, firstly it checks all parking spots in $[1,a_i]$, finding them occupied, then it checks spots $[a_i+1,j-1]$, finding them occupied as well, and finally fills spot $j$. Observe that parking spots $[1,j-1]$ cannot have been filled by cars with preference at least $j$, because such cars would have checked spot $j$ (and filled it) before checking and filling any of the spots in $[1,j-1]$. By hypothesis, spot $j$ is also filled by a car with preference lower than $j$, thus spots $[1,j]$ are only filled by cars with preference in $[1,j-1]$, so  $\sum_{i=1}^{j-1}\A_i\ge j$, and consequently using \eqref{propua} 
	\[
	\ua(j)=j-1-\sum_{i=1}^{j-1}\A_i\le-1.
	\]
	Suppose now that $a_i> k$, then when car $c_i$ attempts to park it finds that all spots in $[a_i-k,a_i]$ are occupied, and so are spots $[a_i+1,j-1]$. In particular, parking spots $[j-k,j-1]$ are already occupied when car $c_i$ parks in spot $j$, thus like in the previous case, they can only have been occupied by cars with preference in $[1,j-1]$. Note that cars with preference at least $j$ cannot park in any parking spot prior to $j-k$, so any spot in $[1,j-k-1]$ can also only have be occupied by a car with preference in $[1,j-1]$. Hence, like in the previous case, no spot in $[1,j]$ can be filled by a car with preference in $[j,n]$.\\
	Finally, if $\alpha\in\knap$, then specifically all spots in $[1,j]$ are occupied, thus by \eqref{propua} we again obtain $\ua(j)\le-1$.
\end{proof}
\end{comment}

\begin{proof}
	Consider $h\in[n]$ such that $a_h\in[j,n]$. If $h<i$, then spot $j$ has not yet been occupied by car $c_i$, thus car $c_h$ cannot park in any spot in $[1,j]$ without occupying spot $j$ as well. \\
	On the other hand, suppose that $h>i$. Note that car $c_i$ has parked by driving forwards, meaning that  it has checked the $k$ spots prior to $a_i$, finding them already occupied, and then it has also checked all spots in $[a_i,j-1]$, finding them full as well. In particular, after car $c_i$ has parked, all spots in $[j-k,j]$ (or $[1,j]$, if $j\le k$) are full. Thus since $h>i$, car $c_h$ cannot possibly occupy any spot prior to $j$, since it can only drive backwards up to $k$ spaces.\\
	Finally, suppose that $\alpha\in\pf$, then all spots  in $[1,j]$ are full, and since they are all occupied by cars with preference in $[1,j-1]$, it means that there are at least $j$ cars with preference in $[1,j-1]$, and by \eqref{propua} we get
	\[
	\ua(j)=j-1 -\sum_{i=1}^{j-1}\A_j\le -1.
	\]
\end{proof}

\begin{remark}
	The inequality $\ua(j)\le-1$ is just a technical condition, though it often allows us to find possible contradictions when assuming that certain cars park by driving forwards. On the other hand, the observation that no cars with preference in $[j,n]$ park in $[1,j]$ is quite significant: specifically, it expresses the fact that, whenever a car $c_i$ parks in a spot $j$ by driving forwards, then some sort of separation between spots $[1,j]$ and spots $[j+1,n]$ occurs, just like a one-way valve, that prevents any backward movement by any successive car.\\
	Furthermore, note that the conditions found in Lemma \ref{menouno} strongly rely on the fact that all cars follow the same parking rule: if cars are allowed to have different parking rules (as we will see in a later chapter), then those restrictions are lifted.\\
	Finally, we want to emphasize the fact that Lemma \ref{menouno} holds for general parking preferences, not just when they are complete, as it will be frequently used in  the sequel.
\end{remark}

The next lemma shows that a complete parking preference is a $k$-Naples parking function if and only if all cars park by driving backwards (or, equivalently, all spots are occupied by cars driving backwards).

\begin{proposition}
\label{backprop}
Let $\alpha\in\pp$ be complete, and assume that all cars follow the $k$-Naples parking rule. Then the following are equivalent:
\begin{enumerate}[(i)]
\item $\alpha\in\knap$;\\
\item for all $j\in[n]$, parking spot $j$ is filled by a car whose preference is at least $j$;\\
\item for all $j\in[n]$, $\psi_k(c_j)\le a_j$.
\end{enumerate}
\end{proposition}

\begin{proof}
Suppose $\alpha\in\knap$. If, for some $j\in[n]$, parking spot $j$ is filled by a car $c_i$ such that $a_i<j$ (and so $j\ge2$), then, by Lemma \ref{menouno}, $\ua(j)\le-1$. However $\alpha$ is complete, so $\ua(j)\ge1$ for all $j\in[2,n]$, and we have a contradiction. Hence $(i)\implies(ii)$.\\
Now, clearly if for some $j\in[n]$, $\psi_k(c_j)>a_j$, then either $\psi_k(c_j)=\infty$, meaning that some cars do not park and thus some spots remain unoccupied, or $\psi_k(c_j)<\infty$, meaning that spot $\psi_k(c_j)$ is filled by a car with preference smaller than $\psi_k(c_j)$. In both cases, that contradicts \emph{(ii)}, thus $(ii)\implies(iii)$.\\
Finally, if for all $j\in[n]$, $\psi_k(c_j)\le a_j<\infty$, in particular all cars park, i.e. $\alpha\in\knap$ and we have that $(iii)\implies(i)$.\\
\end{proof}

The next couple of results put restrictions on the number of cars taking specific backward steps. Specifically, if $\alpha\in\pp$ is complete, then the number of cars taking the backward step from spot $j$ to spot $j-1$ is at most $\ua(j)$, and this upper bound is attained for all $j$ only when $\alpha$ is a $k$-Naples parking function. Notably, this restrictions hold only for complete parking preferences: for generic $k$-Naples parking functions, it may happen that more than $\ua(j)$ cars can take the backward step from $j$ to $j-1$, and this is compensated by the fact that there can also be cars that take the \emph{forward} step from spot $j-1$ to $j$.

\begin{proposition}
	\label{mostprop}
Let $\alpha\in\pp$ be complete, and assume all cars follow the \emph{k}-Naples parking rule. Then for all $j\in[n]$ the number of cars with preference at least $j$ that park in spots strictly smaller than $j$ is at most $\ua(j)$. In other words,
\begin{equation}
\label{eqmost}
\abs*{\{i\in[n]\mid a_i\ge j,\, \psi_k(c_i)<j\}}\le \ua(j) \quad \textit{for all $j\in[n]$}.
\end{equation}
\end{proposition}

\begin{proof}
Suppose $j=n$, then by definition $\A_n=\ua(n)+1$, where $\ua(n)\ge1$ because $\alpha$ is complete. Suppose that parking spot $n$ is filled by a car $c_\lambda$ such that $a_\lambda<n$, then, 
%necessarily no car with preference $n$ has yet attempted to park (otherwise spot $n$ would have been filled by one of those cars). Furthermore, just as we remarked in Lemma \ref{menouno}, for car $c_\lambda$ to be able to park in spot $n$ then all spots in $[n-k,n-1]$ must have already been filled (by cars with preference in $[1,n-1]$). So any car with preference $n$ finds spots $[n-k,n]$ already occupied, and consequently is not able to park; in particular, 
by Lemma \ref{menouno}, no car with preference $n$ parks in $[1,n]$. In particular, no car with preference $n$ parks in any spot smaller than $n$, and we get
\[
\abs*{\{i\in[n]\mid a_i\ge n,\, \psi_k(c_i)<n\}}=0< \ua(n).
\]
Assume instead that spot $n$ is filled by a car with preference $n$; aside from that car, there are only $\ua(n)$ more cars with preference $n$, so trivially at most $\ua(n)$ cars with preference $n$ park in spots strictly smaller than $n$.\\
Now let $2\le j<n$ and suppose by induction that \eqref{eqmost} is satisfied for $j+1$. Assume spot $j$ is filled by a car $c_{\lambda}$ with preference strictly smaller than $j$: then, again by Lemma \ref{menouno}, no car with preference at least $j$ parks in any spot in $[1,j]$, and \eqref{eqmost} is trivially satisfied.\\
Assume instead that spot $j$ is filled by a car with preference at least $j$. Note that 
\begin{multline*}
 \{i\in[n]\mid a_i\ge j,\, \psi_k(c_i)<j\}\\
 =\{i\in[n]\mid a_i= j,\, \psi_k(c_i)<j\}\cup\{i\in[n]\mid a_i>j,\, \psi_k(c_i)<j\},
 \end{multline*}
and that the union is disjoint. Clearly $\abs*{\{i\in[n]\mid a_i= j,\, \psi_k(c_i)<j\}}\le\A_j$, moreover
 \begin{equation*}
 \abs*{\{i\in[n]\mid a_i> j,\, \psi_k(c_i)<j\}}
 \le\abs*{\{i\in[n]\mid a_i\ge j+1,\, \psi_k(c_i)<j+1\}}\le\ua(j+1),
 \end{equation*}
where the last inequality holds by induction. Now, since we assumed that spot $j$ is filled by a car with preference at least $j$, that car must either have preference $j$, or be one of the cars in $\{i\in[n]\mid a_i\ge j+1,\, \psi_k(c_i)<j+1\}$ (that are at most $\ua(j)$). As such, in the former case clearly $\abs*{\{i\in[n]\mid a_i= j,\, \psi_k(c_i)<j\}}\le\A_j-1$, while in the latter case $\abs*{\{i\in[n]\mid a_i> j,\, \psi_k(c_i)<j\}}\le\ua(j+1)-1$. In any case, we obtain
 \[
 \abs{\{i\in[n]\mid a_i\ge j,\, \psi_k(c_i)<j\}}\le \A_j+\ua(j+1)-1=\ua(j),
 \]
 where the last equality is given by \eqref{propua3}.
\end{proof}

\begin{proposition}
\label{exactprop}
Let $\alpha\in\knap$ be complete. Then for all $j\in[n]$ the number of cars with preference at least $j$ that park in spots strictly smaller than $j$ is \emph{exactly} $\ua(j)$. In other words,
\begin{equation}
\abs*{\{i\in[n]\mid a_i\ge j,\, \psi_k(c_i)<j\}}= \ua(j), \qquad \textit{for all $j\in[n]$}.
\end{equation}
\end{proposition}

\begin{proof}
Let $j\in[n]$, since $\alpha\in\knap$, Proposition \ref{backprop} implies that all spots in $[j,n]$ are filled by cars with preference at least $j$. Note that $[j,n]$ contains $n-j+1$ spots. Furthermore all cars with preference at least $j$ are able to park, and consequently the number of cars with preference at least $j$ that park in spots strictly smaller than $j$ is exactly
\[
\sum_{i=j}^n\A_i-(n-j+1)=\ua(j).
\]
\end{proof}

\section[Characterization of complete $k$-Naples parking functions]{Characterization of complete $k$-Naples parking \\functions}

In this section, we will characterize complete $k$-Naples parking functions, focusing on conditions on the order of cars in a complete parking preference. For example, a complete $1$-Naples parking function needs to be in descending order (for a proof, see Section \ref{section1nap}). 
If $k\ge2$, conditions are not so strict: the next theorem shows that, in essence, the parking preference needs to have "not too many" cars with small preferences coming before cars that have larger preferences, in order to allow the latter cars to park by driving backwards.

\begin{theorem}
\label{charUcomp}
Let $\alpha\in\pp$ be complete, and $k\ge1$. Then the following are equivalent:
\begin{enumerate}[(i)]
\item $\alpha\in\knap$;\\
\item for all $h\in[n]$ such that $\A_h\ne0$, there exists $ 0\le\lambda_h\le k-\ua(h)$ such that
\begin{equation}
	\label{Ucomp2}
\abs*{\{i<j_h\mid a_i\in[h-\ua(h)-\lambda_h,h-1]\}}\le\lambda_h,
\end{equation}
where $j_h$ is the last index such that $a_{j_h}=h$, i.e. 
\[
j_h\coloneqq\min\{j\in[n]\mid a_i\ne h\quad\forall\, i>j\};
\]
\item for all $h\in[n]$, there exists $0\le\lambda_h\le k-\ua(h)$ such that
\begin{equation}
	\label{Ucomp3}
\abs*{\{i<j_h^*\mid a_i\in[h-\ua(h)-\lambda_h,h-1]\}}=\lambda_h,
\end{equation}
where $j_h^*$ is the last index such that $a_{j_h^*}\ge h$, i.e. 
\[
j_h^*\coloneqq\min\{j\in[n]\mid a_i< h\quad\forall\, i>j\}.
\]
\end{enumerate}
\end{theorem}

\begin{proof}
Note that clearly $j_h\le j_h^*$, so the implication \emph{(iii)}$\implies$\emph{(ii)} is trivial.\\
{\bf Case \emph{(i)}$\implies$\emph{(iii)}}.
Let $h\in[n]$. Note that if $h=1$, then $\lambda_h=0$ trivially satisfies \eqref{Ucomp3}, thus we can suppose $h\ge2$. Let $j_h^*\in[n]$ be defined as above, and consider car $c_{j_h^*}$. Suppose that $\psi_k(c_{j_h^*})\ge h$: Proposition \ref{backprop} implies that spots $[\psi_k(c_{j_h^*}),n]$ are all filled by cars with preference in $[\psi_k(c_{j_h^*}),n]$. All cars that precede $c_{j_h^*}$ and have preference in $[\psi_k(c_{j_h^*}),n]$ cannot have filled any spots smaller than $\psi_k(c_{j_k^*})$, since otherwise one of them would have also occupied spot $\psi_k(c_{j_h^*})$. Recalling that $c_{j_h^*}$ is the last car with preference in $[h,n]$, in particular it is the last car with preference in $[\psi_k(c_{j_h^*}),n]$ (Proposition \ref{backprop} implies that $a_{j_h^*}\ge\psi_k(c_{j_h^*})$, and we are assuming $\psi_k(c_{j_h^*})\ge h$), and consequently all cars with preference in $[\psi_k(c_{j_h^*}),n]$ park in spots $[\psi_k(c_{j_h^*}),n]$. Thus $\ua(\psi_k(c_{j_h^*}))=0$, and we have a contradiction since $\alpha$ is complete, and $\psi_k(c_{j_h^*})\ge h \ge2$.\\
Therefore we can assume that $\psi_k(c_{j_h^*})< h$. Proposition \ref{exactprop} ensures that exactly $\ua(h)$ cars with preference at least $h$ park in spots strictly smaller than $h$: since $c_{j_h^*}$ is the last one of such cars, the others must already have parked, so certainly car $c_{j_h^*}$ finds (at least) all parking spots in $[h-(\ua(h)-1),h]$ already occupied, hence $\psi_k(c_{j_h^*})\le h-\ua(h)$. Furthermore, $\psi_k(c_{j_h^*})\ge a_{j_k*}-k\ge h-k$.\\
Let $\lambda_h=h-\psi_k(c_{j_h^*})-\ua(h)$. The previous inequalities imply that $\lambda_h$ satisfies $0\le\lambda_h\le k-\ua(h)$. Consider parking spots $[\psi_k(c_{j_h^*}),h-1]$: $\ua(h)$ spots (one of which is spot $\psi_k(c_{j_h^*})$) are occupied by cars with preference at least $h$. The remaining 
\[
\abs{[\psi_k(c_{j_h^*}),h-1]}-\ua(h)= h-\psi_k(c_{j_h^*})-\ua(h)=\lambda_h
\]
spots are therefore occupied by cars with preference in $[\psi_k(c_{j_h^*}),h-1]$. Let $i<j_h^*$ be such that $a_i\in[\psi_k(c_{j_h^*}),h-1]$: by Proposition \ref{backprop} car $c_i$ parked in $[1,h-1]$, but spot $\psi_k(c_{j_h^*})$ has to remain empty, so $\psi_k(c_i)\in[\psi_k(c_{j_h^*})+1,h-1]$. Thus, recalling that $\psi_k(c_{j_h^*})=h-\ua(h)-\lambda_h$ we get:
\begin{multline*}
\abs*{\{i<j_h^*\mid a_i\in[h-\ua(h)-\lambda_h,h-1]\}}\\
=\abs*{[h-\ua(h)-\lambda_h,h-1]}-\ua(h)\\
=h-1-(h-\ua(h)-\lambda_h)+1-\ua(h)=\lambda_h.
\end{multline*}
{\bf Case \emph{(ii)}$\implies$\emph{(i)}}.
Suppose $\alpha\notin\knap$, then Proposition \ref{backprop} implies that there exist cars such that $\psi_k(c_j)>a_j$. Let $h\in[n]$ and $j\in[n]$ be such that $h$ is the minimum number such that there is a car $c_j$ with preference $h$ that does not park backwards (in particular, $\A_h\ne0$). In other words, let $h,j\in[n]$ 
be such that $\psi_k(c_j)>h=a_j$, and $\psi_k(c_i)\le a_i$ for all $i\in[n]$ such that $a_i<h$. By hypothesis \emph{(ii)}, there exists $0\le\lambda_h\le k-\ua(h)$ such that 
\begin{equation*}
\abs*{\{i<j\mid a_i\in[h-\ua(h)-\lambda_h,h-1]\}}
\le\abs*{\{i<j_h\mid a_i\in[h-\ua(h)-\lambda_h,h-1]\}}\le\lambda_h.
\end{equation*}
Since car $c_j$ parks in a spot strictly larger than its preference, on its turn it finds all spots in $[h-\ua(h)-\lambda_h,h-1]$ already occupied, and since we supposed that all cars with preference strictly smaller than $h$ park in spots smaller than their preference, at most $\lambda_h$ spots in $[h-\ua(h)-\lambda_h,h-1]$ have been filled by such cars. Consequently, at least $\abs{[h-\ua(h)-\lambda_h,h-1]}-\lambda_h=\ua(h)$ spots have been filled by cars with preference at least $h$. On the other hand, from Proposition \ref{mostprop} we get that at most $\ua(h)$ cars with preference at least $h$ can have parked in $[h-\ua(h)-\lambda_h,h-1]$, so the only possibility is that \emph{exactly} $\ua(h)$ cars with preference at least $h$ have parked in $[h-\ua(h)-\lambda_h,h-1]$.\\ 
Finally, note that spot $h$ has also been filled by a car with preference at least $h$, because by construction of $h$ no cars with preference at most $h-1$ could have occupied it. In particular, it has been filled either by one of the $\A_h$ cars with preference $h$, or (again using Proposition \ref{mostprop}) by one of the cars with preference at least $h+1$ parking in spots at most $h$ (of which there are at most $\ua(h+1)$). Note that (by \eqref{propua3}) $\ua(h+1)+\A_h=\ua(h)+1$, meaning that $\ua(h)+1$ spots in $[h-\ua(h)-\lambda_h,h-1]\cup\{h\}$ are filled by $\A_h$ cars with preference $h$, and $\ua(h+1)$ cars with preference at least $h+1$. In particular, \emph{all} cars with preference $h$ park in spots at most $h$, and we have a contradiction with the definition on $h$.
\end{proof}

\begin{remark}
	\label{remarksatisfy}
	Note that condition \eqref{Ucomp2} is always satisfied for $h\le k+1$. Let $\lambda_h=(h-1)-\ua(h)\le k-\ua(h)$, where \eqref{propua4} ensures that $\lambda_h\ge0$. Then
	\begin{multline*}
	\abs{\{i<j_h\mid a_i\in[h-\ua(h)-\lambda_h,h-1]\}}\\
	=\abs{\{i<j_h\mid a_i\in[1,h-1]\}}\le \sum_{i=1}^{h-1}\A_i=h-1-\ua(h)=\lambda_h.
	\end{multline*}
\end{remark}

\begin{remark}
	\label{reformulate Ucomp}
	The significance of condition \eqref{Ucomp3} is the following. Assuming $\alpha\in\knap$, in order to be able to park "correctly", the last car with preference at least $h$ needs to park "backwards", i.e. in a spot smaller than $h$ (if it doesn't, we have a contradiction). To be able to do so, at least one of the $k$ spots prior to $h$ needs to be still free: assuming this spot is $h-\ua(h)-\lambda_h$, then no more than $\lambda_h$ cars with preference in $[h-\ua(h)-\lambda_h,h-1]$ can already have had their turn.\\
	Moreover, using \eqref{propua2}, condition \eqref{Ucomp3} can easily be reformulated as follows (note that $a_{j_h}\notin[h-\ua(h)-\lambda_h, h-1]$):
	\begin{multline*}
		\abs*{\{i>j_h^*\mid a_i\in[h-\ua(h)-\lambda_h,h-1]\}}\\
		=\abs*{\{i\mid a_i\in[h-\ua(h)-\lambda_h,h-1]\}}
		-\abs*{\{i<j_h^*\mid a_i\in[h-\ua(h)-\lambda_h,h-1]\}}\\
		=\sum_{i=h-\ua(h)-\lambda_h}^{h-1} \A_i-\lambda_h
		=\sum_{i=h-\ua(h)-\lambda_h}^{h-1} \A_i -(\lambda_h+\ua(h))+\ua(h)\\
		=-\ua(h)+\ua(h-\ua(h)-\lambda_h)+\ua(h)
		=\ua(h-\ua(h)-\lambda_h),
	\end{multline*}
	where $ 0\le\lambda_h\le k-\ua(h)$.
	Replacing $\eta_h=\lambda_h+\ua(h)$, we obtain
	\begin{equation}
		\abs*{\{i>j_h^*\mid a_i\in[h-\eta_h,h-1]\}}= \ua(h-\eta_h),
	\end{equation}
	where $\ua(h)\le\eta_h\le k$.\\
	In the same way, \eqref{Ucomp2} can be rewritten as 
	\[
	\abs*{\{i>j_h^*\mid a_i\in[h-\eta_h,h-1]\}}\ge \ua(h-\eta_h),
	\]
	where $\ua(h)\le\eta_h\le k$.\\
\end{remark}

This reformulation allows us to focus on cars coming \emph{after} car $c_{j_h}$, rather than before. With this in mind, we recall the concept of \emph{RTL-maximum}.

\begin{definition}
	Let $\alpha\in\pp$. We call $a_j$ a \emph{RTL-maximum} (right-to-left-maximum) if 
	\begin{equation}
		a_j>\max{\{a_i\mid i>j\}},
	\end{equation}
i.e. $a_j$ is greater than all elements to its right.
\end{definition}

The following theorem is an improvement over Theorem \ref{charUcomp}, where condition \eqref{Ucomp3} (or \eqref{Ucomp2}) is checked only for RTL-maxima.
\begin{theorem}
	\label{rtltheorem}
	Let $\alpha\in\pp$ be complete, and $k\ge1$. Then $\alpha\in\knap$ if and only if for all $j\in[n]$ such that $a_j$ is a RTL-maximum, there exists $0\le\lambda_j\le k-\ua(a_j)$ such that
	\begin{equation}
		\label{Ucomp4}
		\abs*{\{i<j\mid a_i\in[a_j-\ua(a_j)-\lambda_j,a_j-1]\}}	\le\lambda_j.
	\end{equation}
	Equivalently, $\alpha\in\knap$ if and only if for all $j\in[n]$ such that $a_j$ is a RTL-maximum, there exists $\ua(a_j)\le\eta_j\le k$ such that
	\begin{equation}
		\label{Ucomp5}
		\abs*{\{i>j\mid a_i\in[a_j-\eta_j,a_j-1]\}}\ge\ua(a_j-\eta_j).
	\end{equation}
	Finally, if $\alpha\in\knap$, then \eqref{Ucomp4} and \eqref{Ucomp5} are equalities for $\lambda_j=a_j-\psi_k(c_j)-\ua(a_j)$ and $\eta_j=a_j-\psi_k(c_j)$ respectively.
\end{theorem}

\begin{proof}
	We showed in Remark \ref{reformulate Ucomp} that \eqref{Ucomp4} and \eqref{Ucomp5} are equivalent. Furthermore, the fact that \eqref{Ucomp4} is an equality when $\lambda_j=a_j-\psi_k(c_j)-\ua(a_j)$ (and therefore \eqref{Ucomp5} is also an equality) follows from the proof of the implication (i)$\implies$(iii) in Theorem \ref{charUcomp}. 
	Furthermore, the fact that $\alpha\in\knap$ implies \eqref{Ucomp4} is a simple consequence of Theorem \ref{charUcomp}. Thus, all we need to prove is that, if \eqref{Ucomp4} holds (for some $\lambda_j$) for all RTL-maxima $a_j$, then $\alpha$ is a $k$-Naples parking function.\\
	Suppose $\alpha\notin\knap$; by Proposition \ref{backprop}, there exist cars that park in spots larger than their preference. Note that, if $h\in[n]$ is such that $\psi_k(c_p)>h$ for some cars $c_p$ such that $a_p=h$, then in particular $\psi_k(c_{j_h})>h$, where $j_h$ is defined as in Theorem \ref{charUcomp}. Let $h$ be the biggest number such that $\psi_k(c_{j_h})>h$.\\
	Assume that $a_{j_h}$ is not a RTL-maximum; given $j^*=j_h^*$ defined as in Theorem \ref{charUcomp}, clearly by definition of $j^*$ we have that $a_{j^*}$ is RTL-maximum, and furthermore $j^*>j_h$ and $a_{j^*}>h$. By construction of $h$, this implies $\psi_k(c_{j^*})\le a_{j^*}$.
	Suppose $\psi_k(c_{j^*})\in[h+1,a_{j^*}]$, then since $c_{j^*}$ is the last car with preference at least $\psi_k(c_{j^*})$, all cars with preference in $[\psi_k(c_{j^*}),n]$ park in $[\psi_k(c_{j^*}),n]$, and so $\ua(\psi_k(c_{j^*}))=0$, which is a contradiction, because $\alpha$ is complete.
	Furthermore, $\psi_k(c_{j^*})\le h$ is not possible, since if there were a free parking spot available for $c_{j^*}$ by driving backwards from $a_{j^*}>h$, then car $c_{j_h}$ could have reached it as well ($j^*>j_h$). So assuming that $a_{j_h}$ is not a RTL-maximum leads to a contradiction.\\
	Therefore, $a_{j_h}$ is a RTL-maximum; following the proof of \emph{(ii)$\implies$(i)} in Theorem \ref{charUcomp} we obtain that
	\[
	\abs*{\{i<j_h\mid a_i\in[h-\ua(h)-\lambda_h,h-1]\}}>\lambda_h
	\]
	for all $0\le\lambda_h\le k-\ua(h)$. We conclude by setting $j=j_h$ and $h=a_{j_h}$.\\
\end{proof}

Analysing complete preferences with respect to RTL-maxima also gives the following easy necessary condition.

\begin{lemma}
	\label{lemmartl}
	Let $\alpha\in\knap$ be complete, and denote $j_1<\dots<j_m$ all the indices such that $a_{j_i}$ ($1\le i\le m$) are RTL-maxima. Then $a_{j_{i}}-a_{j_{i+1}}\le k$ for all $1\le i\le m-1$, and $a_{j_m}\le k+1$.
\end{lemma}

\begin{proof}
	Let $1\le i < m$ such that $a_{j_{i}}-a_{j_{{i+1}}}> k$. Then all other cars with preference in $[a_{j_i}-k,a_{j_i}]$ have parked before $c_{j_i}$. If $\alpha\in\knap$, then $\psi_k(c_{j_i})\in[a_{j_i}-k,a_{j_i}]$ by Proposition \ref{backprop}. Then in an analogous way as the previous theorems, $\ua(\psi_k(c_{j_i}))=0$, and we have a contradiction since $\psi_k(c_{j_i})\ge a_{j_i}-k>a_{j_{i+1}}\ge 1$.\\
	If $a_{j_m}>k+1$, then again all other cars with preference in $[a_{j_m}-k,a_{j_m}]$ have parked before $c_{j_m}$ and in the same way we obtain $\ua(\psi_k(c_{j_m}))=0$, which is a contradiction because $\psi_k(c_{j_m})\ge a_{j_m}-k>1$.\\
\end{proof}

To show some applications of the previous theorems, as well as to examine some aspects that may not be obvious, we provide the following example.

\begin{example}
	Let $\alpha=(7,8,7,5,8,4,5,2)$. Observe that $\alpha$ is not a parking function, and that it is a complete parking preference. Its RTL-maxima are $a_5=8$, $a_7=5$ and $a_8=2$, where $\ua(8)=1$, $\ua(5)=2$ and $\ua(2)=1$. Note that the last RTL-maximum is such that $a_8=2\le k+1$ for all $k\ge1$ so, as we have noticed in Remark \ref{remarksatisfy}, condition \eqref{Ucomp4} can always be satisfied, therefore we only need to check what happens for the other two RTL-maxima. Note that, by Theorem \ref{necessk},  $\alpha$ is not a $1$-Naples parking function since $\ua(5)=2$.\\
	Set $k=2$, we check condition \eqref{Ucomp4} with $a_7=5$: since $\ua(5)=k=2$, the only choice for $\lambda_7$ is $\lambda_7=0$, and we get
	\begin{equation*}
	\abs*{\{i<7\mid a_i\in[a_7-\ua(a_7)-\lambda_7,a_7-1]\}}
	=\abs*{\{i<7\mid a_i\in[3,4]\}}=1>0=\lambda_7.
	\end{equation*}
	Thus $\alpha\notin PF_{8,2}$. Indeed, checking where cars park gives
	\[
	\psi_2((7,8,7,5,8,4,5,2))=(7,8,6,5,\infty,4,3,2).
	\]
	
	Note that, in this example, even though condition \eqref{Ucomp4} fails for $a_7=5$, that does not imply that car $c_7$ could not park: that is because $a_7=5$ was not the only case where \eqref{Ucomp4} fail. In fact, as we will see, \eqref{Ucomp4} also fails for $0\le\lambda_5\le2= k-\ua(8)$. Now, recall that a key factor in the last theorems is that, when supposing $\alpha$ to be a $k$-Naples parking function, we assume that a total of $\ua(j)$ cars with preference at least $j$ parks in spots smaller than $j$ (see Proposition \ref{exactprop}). This does not happen here for $j=7$, since some cars before car $c_7$ fail to park (i.e. car $c_5$), and that is why car $c_7$ succeeds in parking driving backwards.\\
	To better illustrate this, consider $\alpha_1=(8,8,7,7,5,4,5,2)$, that is obtained from $\alpha$ by rearranging all cars coming before car $c_7$ in the most favourable way, i.e. in descending order. In this case, the only instance in which condition \eqref{Ucomp4} fails is for $j=7$, and indeed we get
	\[
	\psi_2(\alpha_1)=(8,7,6,5,4,3,\infty,2).
	\]
	Thus, in general, condition \eqref{Ucomp4} failing for a specific RTL-maximum does not necessarily mean that the associated car does not park, but it does still imply that something does go wrong for some car.\\

	Now, we get back to the original parking preference $\alpha=\{7,8,7,5,8,4,5,2\}$, and set $k=3$. Setting $\lambda_7=1=k-\ua(a_7)$ gives
	\[
	\abs*{\{i<7\mid a_i\in[a_7-\ua(a_7)-\lambda_7,a_7-1]\}}
	=\abs*{\{i<7\mid a_i\in[2,4]\}}=1=\lambda_7,
	\]
	so \eqref{Ucomp4} is satisfied for $j=7$. Finally, we focus on the last remaining RTL-maximum ($a_5=8$): we need $0\le\lambda_5\le k-\ua(8)=2$. We obtain:
	\begin{itemize}
		\item if $\lambda_5=0$, then 
		\begin{equation*}
		\abs*{\{i<5\mid a_i\in[a_5-\ua(a_5)-\lambda_5,a_5-1]\}}
		=\abs*{\{i<5\mid a_i=7\}}=2>0=\lambda_5;
		\end{equation*}
		\item if $\lambda_5=1$, then 
		\begin{equation*}
			\abs*{\{i<5\mid a_i\in[a_5-\ua(a_5)-\lambda_5,a_5-1]\}}
			=\abs*{\{i<5\mid a_i\in[6,7]\}}=2>1=\lambda_5;
		\end{equation*}
		\item if $\lambda_5=2$, then 
		\begin{equation*}
			\abs*{\{i<5\mid a_i\in[a_5-\ua(a_5)-\lambda_5,a_5-1]\}}
			=\abs*{\{i<5\mid a_i\in[5,7]\}}=3>2=\lambda_5.
		\end{equation*}
	\end{itemize}
	Hence, \eqref{Ucomp4} always fails for $j=5$, and we notice that $\psi_3(c_5)=\infty$.
	
	Finally, let $k=4$. In this case, setting $\lambda_5=3=k-\ua(a_5)$ we get
	\[
	\abs*{\{i<5\mid a_i\in[a_5-\ua(a_5)-\lambda_5,a_5-1]\}}=\abs*{\{i<5\mid a_i\in[4,7]\}}=3=\lambda_5.
	\]
	Thus, $\lambda_5=3$ satisfies \eqref{Ucomp4} for $j=5$, and indeed $\alpha$ is a $4$-Naples parking function, and
	\[
	\psi_4(7,8,7,5,8,4,5,2)=(7,8,6,5,4,3,2,1).
	\]
\end{example}

\begin{remark}
	As a final observation, consider a parking preference $\alpha\in\pp$ such that $a_{j_1},\dots,a_{j_m}$ ($j_1<\dots<j_m$) are its RTL-maxima. For all $i=1,\dots,m$, let $k_i$ be the minimum value such that \eqref{Ucomp4} holds for $j=j_i$ and $\lambda_{j_i}=k_i-\ua(a_{j_i})$, and let $k=\min_i k_i$. Then, as it happens in the previous example, $k$ is the smallest value such that $\alpha\in\knap$.\\
	Moreover, having $0\le\lambda_j< k-\ua(j)$ such that $\lambda_j$ satisfies \eqref{Ucomp4} for the RTL-maximum $a_j$ does not imply that \eqref{Ucomp4} is also satisfied for larger values of $\lambda_j$: in other words, given $k$, 
	\emph{all} values of $0\le\lambda_j\le k-\ua(j)$ need to be checked, it is not sufficient to just check $\lambda=k-\ua(j)$. For example, let $\alpha=(7,8,5,5,8,7,4,2)$, where $a_5=8$ is a RTL-maximum, and let $k=3$. Setting $\lambda_5=1<2=k-\ua(8)$ gives
	\[
	\abs*{\{i<5\mid a_i\in[a_5-\ua(a_5)-\lambda_5,a_5-1]\}}=\abs*{\{i<5\mid a_i\in[6,7]\}}=1=\lambda_5.
	\]
	However, for $\lambda_5=2=k-\ua(8)$ we get
	\[
	\abs*{\{i<5\mid a_i\in[a_5-\ua(a_5)-\lambda_5,a_5-1]\}}=\abs*{\{i<5\mid a_i\in[5,7]\}}=3>2=\lambda_5.
	\]
	
	Finally, note that while condition \eqref{Ucomp4} does impose some restrictions on the order of $\alpha$, these are not so rigid: for example, one can freely rearrange all preferences of cars parking before the leftmost RTL-maximum, without affecting whether $\alpha$ is a $k$-Naples parking function.
\end{remark}

\section[Characterization of $k$-Naples parking functions through complete subsequences]{Characterization of $k$-Naples parking functions \\through complete subsequences}
\label{sectionknaples}

We now move back to general parking preferences. We start with the following generalization of Proposition \ref{exactprop}.

\begin{proposition}
	\label{leastprop}
	Let $\alpha\in\knap$, then for all $j\in[n]$ such that $j\in\Ua$, the number of cars with preference at least $j$ that park in spots strictly smaller than $j$ is \emph{at least} $\ua(j)$, i.e.
	\begin{equation}
		\abs*{\{i\in[n]\mid a_i\ge j,\, \psi_k(c_i)<j\}}\ge \ua(j) \qquad \textit{for all $j\in\Ua$}.
	\end{equation}
\end{proposition}

\begin{proof}
	For $j\in\Ua$, by definition we have $\sum_{i=j}^n\A_i=\ua(j)+(n-j+1)$. At most $n-j+1$ cars with preference at least $j$ can park in $[j,n]$, and since all cars are able to park, the remaining cars (which are at least $\sum_{i=j}^n\A_i-(n-j+1)=\ua(j)$) need to park in spots strictly smaller than $j$.
\end{proof}

\begin{remark}
	Note that Proposition \ref{exactprop} follows trivially from Propositions \ref{mostprop} and \ref{leastprop}.\\
\end{remark}

The next few results will provide a link between arbitrary $k$-Naples parking functions and complete $k$-Naples parking functions.

\begin{theorem}
	\label{theoprec}
	Let $\alpha\in\pp$ and $k\ge1$. Then $\alpha\in\knap$ if and only if for all $[p,q]\subseteq\Ua$ such that $[p,q]$ is a maximal interval, spot $p-1$ is occupied. In particular, if $\alpha\in\knap$ then there exists $i\in[n]$ such that $a_i\ge p$ and $\psi_k(c_i)=p-1$.
\end{theorem}

\begin{proof}
	Clearly if $\alpha\in\knap$ then all spots are occupied, in particular each spot $p-1$ immediately preceding a maximal interval of $\Ua$ is filled as well. Furthermore, by Proposition \ref{propUa}, $\ua(p-1)=0$ and $\A_{p-1}=0$, thus spot $p-1$ must be filled by a car with preference at least $p$ by Lemma \ref{menouno}.\\
	Let $\alpha\notin\knap$, and let $h\in[n]$ be the smallest parking spot that remains unoccupied. All spots in $[1,h-1]$ are therefore occupied, in particular they must have been filled by cars with preference in $[1,h-1]$ (otherwise $h$ would have been filled as well); conversely, any car with preference in $[1,h-1]$ must have parked in $[1,h-1]$ (again because otherwise $h$ would have been filled). Hence, using \eqref{propua}
	\[
	h-1=\sum_{i=1}^{h-1}\A_i=h-1-\ua(h),
	\] 
	from which we get $\ua(h)=0$, which also implies $h\notin\Ua$.\\
	Finally, clearly $\A_h=0$, so  using \eqref{propua3}
	\[
	\ua(h+1)=\ua(h)-\A_h+1=1,
	\]
	so $h+1\in\Ua$, in particular it must be the lower extreme of a maximal interval of $\Ua$.
\end{proof}

\begin{remark}
	Note that the above proof does not exploit the fact that cars follow the $k$-Naples parking rule for a specific $k$, thus the result holds independently from $k$.
	Therefore, given a parking preference $\alpha\in\pp$, Theorem \ref{theoprec} identifies some \emph{critical spots} that are required to be filled by cars driving backwards, and this can be accomplished if $k$ is sufficiently large and/or if the preferences are suitably ordered. Filling such spots is sufficient for all cars to be able to park.\\
	In fact, from this point of view, standard parking functions are characterized by the absence of any such critical spots, since no car can drive backwards. This is consistent with the result found in the previous chapter, i.e. $\Ua=\varnothing$.\\
	Finally, it will be useful to observe that if independent parking rules are set for each car, Theorem \ref{theoprec} would still yield the desired critical spots. However, in that case it  not necessary that all of them are reached by driving backwards. Some results in this direction will be studied in the last chapter.\\
\end{remark}

\begin{definition}
	Let $\alpha\in\pp$ and $w\in\Z$, we define the \emph{translation operator} $\tau_w$ such that $\tau_w(\alpha)=(a_1-w,\dots,a_n-w)$. Typically, we will consider cases where $0\le w<\min{\{a_1,\dots,a_n\}}$, so that $\tau_w(\alpha)\in PP_{n-w}$.\\
	We also extend this definition to sets in the following way: let $W\subset\Z$ be a finite subset of integers, and $w\in\Z$. We define
	\[
	\tau_w(W)=\{i\in\Z\mid i+w\in W\}.
	\]
\end{definition}

\begin{definition}
	Let $\alpha=(a_1,\dots,a_n)\in\pp$, and $J=\{j_1,\dots,j_{\abs{J}}\}\subseteq[n]$, where $j_1<\dots<j_{\abs{J}}$. We denote the \emph{restriction of $\alpha$ to $J$} with
	\begin{equation}
		\alpha_{|_J}=(a_{j_1},\dots,a_{j_{\abs{J}}}).
	\end{equation}
\end{definition}

\begin{lemma}
	\label{translatelemma}
	Let $\alpha\in\pp$, and $j\in[n]$ such that $\ua(j)=0$. Let $J=\{i\in[n]\mid a_i\ge j\}$, and $J^c=[1,n]\setminus J$. Then $\abs{J}=n-j+1$, so $\alpha_{|_{J^c}}\in PP_{j-1}$, and $\tau_{j-1}(\alpha_{|_J})\in PP_{n-j+1}$.\\
	Furthermore, $u_{\alpha_{|_{J^c}}}(i)=\ua(i)$ for all $i\in[1,j-1]$, thus $U_{\alpha_{|_{J^c}}}=\Ua\cap[1,j-1]$, and $u_{\tau_{j-1}(\alpha_{|_J})}(i)=\ua(i+j-1)$, thus $U_{\tau_{j-1}(\alpha_{|_J})}=\tau_{j-1}(\Ua)\cap[1,n-j+1]=\tau_{j-1}(\Ua\cap[j,n])$.
\end{lemma}

\begin{proof}
	Since $\ua(j)=0$, by definition of $\ua(j)$ we get
	\[
	\abs{J}=\sum_{i=j}^n\A_i=\ua(j)+(n-j+1)=n-j+1.
	\]
	Furthermore, for all $i\in[1,j-1]$,
	\[
	u_{\alpha_{|_{J^c}}}(i)=(i-1)-\sum_{k=1}^{i-1}\abs{\alpha_{|_{J^c}}}_k
	=(i-1)-\sum_{k=1}^{i-1}\A_k=\ua(i).
	\]
	Finally, since $\abs{\tau_{j-1}(\alpha_{|_J})}_k=\A_{k+j-1}$ for all $k\in[1,n-j+1]$, we get
	\begin{multline*}
	u_{\tau_{j-1}(\alpha_{|_J})}(i)=\sum_{k=i}^{n-j+1}\abs{\tau_{j-1}(\alpha_{|_J})}_k-((n-j+1)-i+1)\\
	=\sum_{k=i}^{n-j+1}\A_{k+j-1}-(n-(i+j-1)+1)\\
	=\sum_{k=i+j-1}^n\A_k -(n-(i+j-1)+1)=\ua(i+j-1)
	\end{multline*}
	for all $i\in[1,n-j+1]$.
\end{proof}

\begin{lemma}
	\label{restrlemma}
	Let $\alpha\in\knap$, and $j\in[n]$ such that $\ua(j)=0$. Then given $J=\{i\in[n]\mid a_i\ge j\}$, we have that $\tau_{j-1}(\alpha_{|_J})\in PF_{n-j+1,k}$.
\end{lemma}

\begin{proof}
	We can suppose without loss of generality that $j\ge2$, since if $j=1$ then $J=[n]$ and $\tau_{j-1}(\alpha_{|_J})=\alpha$. We denote $\hat{a}_i$, $\hat{c}_i$ and $\hat{\psi}_k$ the preferences, cars and the outcome map of the problem for $\tau_{j-1}(\alpha_{|_J})$ respectively.\\
	Suppose that $\tau_{j-1}(\alpha_{|_J})\notin PF_{n-j+1,k}$. Then there are some unoccupied spots: let $h$ be the largest such spot. If $u_{\tau_{j-1}(\alpha_{|_J})}(h)\le-1$, then by \eqref{propua}, 
	\[
	\sum_{i=1}^{h-1}\abs{\tau_{j-1}(\alpha_{|_J})}_i=(h-1)-u_{\tau_{j-1}(\alpha_{|_J})}(h) \ge h,
	\]
	and since there are only $h-1$ spots in $[1,h-1]$ necessarily spot $h$ would have been filled. So $u_{\tau_{j-1}(\alpha_{|_J})}(h)\ge0$, and thus also $\ua(h+j-1)\ge0$; moreover, clearly $\abs{\tau_{j-1}(\alpha_{|_J})}_h=0$, and so also $\A_{h+j-1}=0$. 
	Recall that $\alpha\in\knap$, and consider spot $h+j-1$ for the original problem: suppose that it is filled by a car with preference strictly smaller than $h+j-1$, then by Lemma \ref{menouno} we would get a contradiction since $\ua(h+j-1)\ge0$.\\
	Thus, spot $h+j-1$ is occupied by a car $c_m$ such that $a_m\in[h+j,n]$. In particular, there is no $i<m$ such that $a_i<h+j-1$ and $\psi_k(c_i)\ge h+j-1$, otherwise spot $h+j-1$ would have already been filled. Conversely, for all $i< m$ such that $a_i\ge h+j-1$, we have $\psi_k(c_i)\in[h+j-1,n]$ (again because otherwise spot $h+j-1$ would be occupied). Let 
	\[
	H=\{i\in[n]\mid i\le m,\, a_i\ge h+j-1\}=\{h_1,\dots,h_{\abs{H}}\},
	\]
	where $h_1<\dots<h_{\abs{H}}$, and consequently $h_{\abs{H}}=m$. Clearly, $\psi_k(c_{h_1})=a_{h_1}$; moreover, for $\tau_{j-1}(\alpha_{|_J})$, $\hat{a}_{h_1}=a_{h_1}-j+1$. Hence, after car $\hat{c}_{h_1}$'s turn, spot $\hat{a}_{h_1}=\psi_k(c_{h_1})-j+1$ has certainly been filled.\\
	Suppose by induction that, for all $i<\lambda$ ($\lambda\le\abs{H}$), after car $\hat{c}_{h_i}$'s turn, spot $\psi_k(c_{h_i})-j+1$ is occupied. Consider car $c_{h_{\lambda}}$: if $\psi_k(c_{h_{\lambda}})=a_{h_{\lambda}}$, then the induction holds in the same way as the case $\lambda=1$. Otherwise, car $c_{h_{\lambda}}$ first checks some  spots $[s,t]$, finding them occupied, before parking in $\psi_k(c_{h_{\lambda}})$. Here, either $t=a_{h_{\lambda}}$ and $s\ge a_{h_{\lambda}}-k$, or $s=a_{h_{\lambda}}-k$ and $t>a_{h_{\lambda}}$; moreover, $s>h+j-1$ if $\lambda<\abs{H}$, because spot $h+j-1$ only gets filled by car $c_{h_{\abs{H}}}=c_m$. Note that spots $[s,t]$ are all of the form $\psi_k(c_{h_i})$, for some $i<\lambda$. By induction, spots $[s-j+1,t-j+1]$ are filled in the problem for $\tau_{j-1}(\alpha_{|_J})$, before car $\hat{c}_{h_{\lambda}}$ tries to park. Hence, since car $\hat{c}_{h_{\lambda}}$ follows the same rule as car $c_{h_{\lambda}}$, car $\hat{c}_{h_{\lambda}}$ certainly reaches spot $\psi_k(c_{h_{\lambda}})-j+1$, and so after its turn that spot is certainly occupied.\\
	Thus, the induction hypothesis holds for all elements of $H$, in particular after car $\hat{c}_{h_{\abs{H}}}$'s turn, some car must have certainly occupied spot
	\[
	\psi_k(c_{h_{\abs{H}}})-j+1=\psi_k(c_m)-j+1=(h+j-1)-j+1=h,
	\]
	and we have a contradiction, since we supposed that $h$ remains unoccupied in the problem for $\tau_{j-1}(\alpha_{|_J})$. Hence,  $\tau_{j-1}(\alpha_{|_J})\in PF_{n-j+1,k}$.
\end{proof}

\begin{remark}
	Lemma \ref{restrlemma} does \emph{not} imply that if $\ua(j)=0$, then all cars with preference in $[j,n]$ park in $[j,n]$, only that the restriction to the set of cars with preference in $[j,n]$, translated by a factor of $j-1$, is still a $k$-Naples parking function. For example, consider $\alpha=(4,4,3,2,3)$, and $k=1$: $\alpha\in PF_{5,1}$. Here $\ua(4)=0$, $\{i\in[5]\mid a_i\ge4\}=\{1,2\}$, and clearly $\tau_3(\alpha_{|_{\{1,2\}}})=(1,1)$ is a 1-Naples parking function of length $2$. However, in the original problem, $\psi_1(c_2)=3\notin[4,5]$.\\
	Moreover, while it is true that the cars having preference in $[j,n]$ correspond to a $k$-Naples parking function, the same is not necessarily true for the cars having preference in $[1,j]$. In other words, given $\alpha$, $j$ and $J$ defined as in Lemma \ref{restrlemma}, let $J^c=[1,n]\setminus J$: then $\tau_{j-1}(\alpha_{|_J})\in PF_{\abs{J},k}$, but not necessarily $\alpha_{|_{J^c}}\in PF_{\abs{J^c},k}$. In the previous example, $\alpha=(4,4,3,2,3)\in PF_{5,1}$ is such that $\ua(4)=0$, and $\tau_3(\alpha_{|_{\{1,2\}}})\in PF_{2,1}$, but $\alpha_{\{3,4,5\}}=(3,2,3)\notin PF_{3,1}$.\\
\end{remark}

We can now proceed with characterizing $k$-Naples parking functions with respect to the existence of complete subsequences. Specifically, we prove that $\alpha\in\pp$ is a $k$-Naples parking function if and only if for each maximal interval $[p,q]\subseteq\Ua$, there exists a subsequence $\alpha_{|_J}\subseteq [p,p+\abs{J}-1]$ that corresponds to a complete $k$-Naples parking function.\\
For ease of reference, we split the double implication into two separate theorems.

\begin{theorem}
	\label{theosuff}
	Let $\alpha\in\pp$ and $k\ge1$. Suppose that for all maximal intervals $[p,q]\subseteq\Ua$ there exists $J\subseteq[n]$, $\abs{J}\ge2$, such that $a_j\in[p-1,p-2+\abs{J}]$ for all $j\in J$, $\tau_{p-2}(\alpha_{|_J})$ is complete and $\tau_{p-2}(\alpha_{|_J})\in PF_{\abs{J},k}$. Then $\alpha\in\knap$.
\end{theorem}

\begin{proof}
	Let $[p,q]\subseteq\Ua$ be a maximal interval, and $J\in[n]$ as above. Consider $\tau_{p-2}(\alpha_{|_J})$: $\tau_{p-2}(\alpha_j)\in[1,\abs{J}]$ for all $j\in J$, furthermore $\tau_{p-2}(\alpha_{|_J})\in PF_{\abs{J},k}$, so it fills all of $[1,\abs{J}]$. We denote $\tilde{\psi}_k$ the outcome map of the restriction $\tau_{p-2}(\alpha_{|_J})$, $\tilde{c}_j$ the cars for the restricted problem, and $\tilde{a}_j=a_j-(p-2)$ their corresponding preferences.\\
	Consider car $c_{j_1}$: clearly $\tilde{\psi}_k(\tilde{c}_{j_1})=a_{j_1}-(p-2)$. In the original problem, either $\psi_k(c_{j_1})=a_{j_1}$, or car $c_{j_1}$ finds spot $a_{j_1}$ already occupied, and parks somewhere else. In any case, after car $c_{j_1}$'s turn, spot $a_{j_1}= \tilde{\psi}_k(\tilde{c}_{j_1})+(p-2)$ is certainly occupied.\\
	By induction, suppose that for all $i\in[1,h-1]$, after car $c_{j_i}$'s turn, spot $\tilde{\psi}_k(\tilde{c}_{j_i})+(p-2)$ is occupied, and consider car $c_{j_h}$. If $\tilde{\psi}_k(\tilde{c}_{j_h})=a_{j_h}-(p-2)$ then we conclude in the same way as above. Otherwise, by Proposition \ref{backprop}, car $\tilde{c}_{j_h}$ parks $m$ spots backwards from $\tilde{a}_{j_h}$, with $m\le k$. Note that all occupied spots $[\tilde{a}_{j_h}-(m-1),\tilde{a}_{j_h}]$ must be $\tilde{\psi}_k(\tilde{c}_{j_i})$ for some $i\in[1,\abs{J}]$. Therefore by induction, on its turn, car $c_{j_h}$ finds all spots in $[\tilde{a}_{j_h}-(m-1)+(p-2),\tilde{a}_{j_h}+(p-2)]$ all occupied: thus, it certainly reaches spot $\tilde{a}_{j_h}-m+(p-2)=\tilde{\psi}_k(\tilde{c}_{j_h})$, and after its turn that spot is necessarily occupied (it doesn't matter whether it is filled by $c_{j_h}$ or by another car that came before).\\
	Therefore, for all spots in $[1,\abs{J}]$ for the restricted problem, the corresponding spots $[p-1,\abs{J}+(p-2)]$ are occupied in the original problem for $\alpha$. In particular, spot $(p-1)$ is filled, and that is the spot immediately preceding maximal interval $[p,q]\subseteq\Ua$. We conclude by Theorem \ref{theoprec}.
\end{proof}

\begin{theorem}
	\label{theoexist}
	Let $\alpha\in\pp$ and $k\ge1$. If $\alpha\in\knap$, then for all maximal intervals $[p,q]\subseteq\Ua$ there exists $J\subseteq[n]$, such that $a_j\in[p-1,p-2+\abs{J}]$ for all $j\in J$, $\tau_{p-2}(\alpha_{|_J})$ is complete and $\tau_{p-2}(\alpha_{|_J})\in PF_{\abs{J},k}$. Furthermore, $p-1+\abs{J}\notin \Ua$ and in particular $q\le p-2+\abs{J}$.
\end{theorem}

\begin{proof}
	Let $[p,q]\subseteq\Ua$ be a maximal interval, note that $\ua(p-1)=0$, and $\A_{p-1}=0$ by Proposition \ref{propUa}. We start by proving the Theorem when $p=2$, before proceeding with the general case.\\
	{\bf Case $\mathbf{p=2}$:}   Suppose $p=2$. Let $M$ be the smallest number such that all spots in $[1,M]$ are filled by cars with preference in $[1,M]$: clearly $M>1$, because $\A_{p-1}=\A_1=0$, furthermore $M$ exists, since at most $M=n$ ($\alpha\in\knap$). Let $J=\{j_1,\dots,j_M\}$, $j_1<\dots<j_M$, be such that $c_{j_1},\dots,c_{j_M}$ are the cars parking in $[1,M]$.
	Observe that $\ua(M+1)\notin\Ua$: this is trivial if $M=n$, otherwise using \eqref{propua}
	\[
	\ua(M+1)=(M+1)-1-\sum_{i=1}^M\A_i\le M-M=0,
	\]
	where $\sum_{i=1}^M\A_i\ge M$ since cars with preference in $[1,M]$ fill at least spots $[1,M]$.\\
	Consider $\alpha_{|_J}$: clearly $\alpha_{|_J}\in PF_{M,k}$. We will now prove that $\alpha_{|_J}$ is complete. Suppose that $\AJ_M=0$, then $[1,M]$ is filled by cars with preference in $[1,M-1]$; in particular, $[1,M-1]$ is filled by cars with preference in $[1,M-1]$ which is absurd by construction of $M$. Suppose that $\AJ_M=1$, and let $j\in[n]$ such that $a_j=M$: if spot $M$ was filled by a car with preference in $[1,M-1]$, then by Lemma \ref{menouno}, car $c_j$ would not be able to park in $[1,M]$, which leads to a contradiction. Hence, spot $M$ is filled by car $c_j$, and thus spots $[1,M-1]$ are filled by cars with preference in $[1,M-1]$, which is again absurd by construction of $M$.\\
	Thus, $\AJ_M\ge 2$, in particular $\uaJ(M)\ge1$. By induction, let $h\in[2,M-1]$ and suppose that $\uaJ(h+1)\ge1$: we want to prove that $\uaJ(h)\ge1$. Recall that $\uaJ(h)=\uaJ(h+1)+\AJ_h -1$. If $\uaJ(h+1)\ge2$, then $\uaJ(h)\ge1$, and if $\AJ_h\ne0$, then $\uaJ(h)\ge\uaJ(h+1)\ge1$. Assume therefore $\uaJ(h+1)=1$, and $\AJ_h=0$. If spot $h$ is filled by a car with preference in $[1,h-1]$, then $[1,h]$ is filled only by cars with preference in $[1,h]$, contradicting the definition for $M>h$.\\
	Hence, spot $h$ is filled by a car with preference in $[h+1,n]$. Assume that there exists $a^*\le h-1$ such that there is a car with preference $a^*$ which parks in a spot $h^*\in[h+1,n]$. Then, again by Lemma \ref{menouno}, we get $\uaJ(h^*)\le-1$, which contradicts our induction hypothesis. Therefore, all cars in $[1,h-1]$ park in $[1,h-1]$. Now, if $\uaJ(h)=0$, then $\sum_{i=1}^{h-1}\A_i=h-1$, so $[1,h-1]$ is completely filled by cars with preference in $[1,h-1]$: since $h\ge2$, this contradicts the construction of M. Thus, $\uaJ(i)\ge1$ for all $i\in[2,n]$, meaning that $\alpha_{|_J}$ is complete.\\
	{\bf General case:} Let $[p,q]$ be an arbitrary maximal interval. Let $\hat{J}\subseteq[n]$ be the set such that $a_i\ge(p-1)$ for all $i\in\hat{J}$, and consider $\tau_{p-2}(\alpha_{|_{\hat{J}}})$. By Lemma \ref{restrlemma}, $\tau_{p-2}(\alpha_{|_{\hat{J}}})$ is a parking function of length $n-(p-2)$, furthermore the maximal interval $[p,q]\subseteq\Ua$ has been translated to the maximal interval $[2,q-p+2]\in U_{\tau_{p-2}(\alpha_{|_{\hat{J}}})}$. Thus, using what we proved for case $p=2$, we can find a nonempty subset $J\subseteq \hat{J}\subseteq[n]$ such that $\tau_{p-2}(\alpha_{|_J})\in PF_{\abs{J},k}$, and is complete. Moreover, for all $j\in J$, $\hat{a}_j\in[1,\abs{J}]$, thus $a_j\in[p-1,\abs{J}+p-2]$, and since $(\abs{J}+1)\notin U_{\alpha_{|_{\hat{J}}}}$, $(\abs{J}+p-1)\notin\Ua$, which also implies $q\le(\abs{J}+p-2)$.
\end{proof}

\begin{example}
	Let $\alpha=(8,4,7,1,6,8,7,5,10,1)$, and $k=2$. Here $\Ua=[4,7]$, and $\psi_2(\alpha)=(8,4,7,1,6,9,5,3,10,2)$, so $\alpha\in PF_{10,2}$. Let $J=\{2,3,5,7,8\}$, where for all $j\in J$, $a_j\in[4,7]$: then $\tau_2(\alpha_{|_{\{2,3,5,7,8\}}})=(2,5,4,5,3)$ which is in $PF_{5,2}$, furthermore $U_{\tau_2(\alpha_{|_{\{2,3,5,7,8\}}})}=[2,5]$ so $\tau_2(\alpha_{|_{\{2,3,5,7,8\}}})$ is complete.\\
	While a set $J$ as described in Theorem	\ref{theoexist} certainly exists, it is not necessarily unique. For example, let $J'=\{1,2,3,6,7,8\}$: for all $j\in J'$, $a_j\in[4,8]$, and $\tau_2(\alpha_{|_{\{1,2,3,6,7,8\}}})=(6,2,5,6,5,3)$ is also complete, and a 2-Naples parking function (of length $6$). Moreover, set $J''=\{1,2,5,6,7,8\}$ also works.\\
\end{example}

\begin{corollary}
	\label{perminv}
	Let $\alpha=(a_1,\dots,a_n)\in\pp$, and $k\ge1$. Then the following are equivalent:
	\begin{enumerate}[(i)]
		\item for all maximal intervals $[p,q]\subseteq\Ua$, $\abs{[p,q]}\le k$;\\
		\item for all $\sigma\in S_n$, where $S_n$ is the set of permutations of length $n$,
		\begin{equation}
			\sigma(\alpha)\coloneqq(a_{\sigma(1)},\dots,a_{\sigma(n)})
		\end{equation}
		is a $k$-Naples parking function of length $n$. In particular, $\alpha\in\knap$.
	\end{enumerate}
\end{corollary}

\begin{proof}
	Let $\alpha=(a_1,\dots,a_n)\in\pp$ be such that for all maximal intervals $[p,q]\subseteq\Ua$, $\abs{[p,q]}\le k$. Let $[p,q]$ be one of such maximal intervals: $\ua(q+1)\le0$, and by \eqref{propua3}, $\ua(q)=\ua(q+1)+\A_q-1$, thus $\A_q\ge \ua(q)+1$. Choose $\ua(q)+1$ distinct cars with preference $q$, and let $H$ be the set of their indices. Now, let
	\[
	J=\{i\in[n]\mid a_i\in[p,q-1]\}\cup H.
	\]
	By Proposition \ref{propUa}, $\ua(p)=1$, therefore using \eqref{propua2}:
	\[
	\abs{J}=\sum_{i=p}^{q-1}\A_i + \ua(q)+1=(q-p)-\ua(q)+\ua(p)+\ua(q)+1=q-p+2=\abs{[p-1,q]}.
	\]
	Consider $\tau_{p-2}(\alpha_{|_J})\in PP_{\abs{J}}$:
	\begin{multline*}
	u_{\tau_{p-2}(\alpha_{|_J})}(q-p+2)=\abs{\tau_{p-2}(\alpha_{|_J})}_{q-p+2}-1\\
	=\AJ_q-1=\abs{H}-1=\ua(q)+1-1=\ua(q)\ge1,
	\end{multline*}
	and by induction if $u_{\tau_{p-2}(\alpha_{|_J})}(i)=\ua(i+p-2)$ for all $i\in[j+1,q-p+2]$
	\begin{multline*}
	u_{\tau_{p-2}(\alpha_{|_J})}(j)=u_{\tau_{p-2}(\alpha_{|_J})}(j+1)+\abs{\tau_{p-2}(\alpha_{|_J})}_{j+1}-1\\
	=\ua(j+1+(p-2))+\AJ_{j+1+(p-2)}-1\\
	=\ua(j+1+(p-2))+\A_{j+1+(p-2)}-1=\ua(j+(p-2)).
	\end{multline*}
	Hence, since $\ua(i)\ge1$ for all $i\in[p,q]$, $u_{\tau_{p-2}(\alpha_{|_J})}(j)\ge1$ for all $j\in[2,\abs{J}]$, meaning that $\tau_{p-2}(\alpha_{|_J})$ is complete. Moreover, since $\abs{J}=(q-p+1)+1=\abs{[p,q]}+1\le k+1$, trivially $\tau_{p-2}(\alpha_{|_J})\in PF_{\abs{J},k}$. Thus, for any maximal interval of $\Ua$ we have found $J$ as required by Theorem \ref{theosuff}, so $\alpha\in\knap$. Furthermore, we constructed $J$ without regard to the order of preferences/cars involved, so any rearrangement of the preferences/cars is still a $k$-Naples parking function, and \emph{(ii)} holds.\\
	Conversely, suppose that there is a maximal interval $[p,q]\subseteq\Ua$ such that $q-p+1>k$. We can assume that $\alpha$ is in increasing order, since to disprove \emph{(ii)} it is sufficient to prove that a certain rearrangement of $\alpha$ is not a $k$-Naples parking function.\\
	Suppose $\alpha\in\knap$, then by Theorem \ref{theoexist} there exists $J\subseteq[n]$ such that $\tau_{p-2}(\alpha_{|_J})\in PF_{\abs{J},k}$, and is complete. Since $\alpha$ is in increasing order, $\tau_{p-2}(\alpha_{|_J})$ is as well, in particular the preference of the last car is $\abs{J}$ ($\abs{\tau_{p-2}(\alpha_{|_J})}_{\abs{J}}>0$ because $\tau_{p-2}(\alpha_{|_J})$ is complete). This preference is clearly a RTL-maximum, but $\abs{J}\ge q-p+2>k+1$, which contradicts Lemma \ref{lemmartl}. Thus $\alpha\notin\knap$.
\end{proof}

\begin{remark}
	An independent proof of a slightly different version of Corollary \ref{perminv} is found in \cite{carvalho:online}, from the point of view of Dyck paths representing ascending and descending parking functions. The authors show that a $k$-Naples parking function is permutation invariant if and only if its ascending rearrangement is also a $k$-Naples parking function. Subsequently, ascending parking preferences are connected to Dyck paths, and it is established that an ascending parking preference is a $k$-Naples parking function if and only if the corresponding Dyck path does not remain under the $y=0$ line for more than $2k$ consecutive steps.\\
\end{remark}

As a final remark, observe that Theorems \ref{theoprec}, \ref{theosuff} and \ref{theoexist} are related in the following way. Let $\alpha\in\pp$, and consider a given maximal interval $[p,q]\subseteq\Ua$. By Theorem \ref{theoprec}, if $\alpha$ is a $k$-Naples parking function, then spot $p-1$ needs to be occupied by a car whose preference is at least $p$: in particular, when only taking into account cars with preference at least $p$, spot $p-1$ is reached. If this happens, Theorem \ref{theoexist} guarantees the existence of a subsequence of $\alpha$ that, when suitably translated by a factor of $p-2$, is a complete $k$-Naples parking function. Conversely, Theorem \ref{theosuff} shows that the presence of such a subsequence implies that spot $p-1$ is reached by some backward driving car.
Furthermore, the proof of Corollary \ref{perminv} shows that if $\abs{[p,q]}\le k$, such a subsequence certainly exists, so in this case spot $p-1$ will surely get filled. Thus, we can sum up our characterization of $k$-Naples parking functions in the following way.

\begin{theorem}
	\label{sumuptheorem}
	Let $\alpha\in\pp$, $k\ge1$, and assume that all cars follow the $k$-Naples parking rule. Let $[p,q]\subseteq\Ua$ be a maximal interval, then the following are equivalent:
	\begin{enumerate}[(i)]
		\item considering only the behaviour of cars with preference at least $p$, spot $p-1$ is reached by some car;
		\item there exists $J\subseteq[n]$, $q-p+2\le\abs{J}$, such that $\tau_{p-2}(\alpha_{|_J})$ is a complete $k$-Naples parking function. 
	\end{enumerate}
	The stated properties are guaranteed to be satisfied if $\abs{[p,q]}\le k$.\\
	Furthermore, $\alpha$ is a $k$-Naples parking function if and only if all maximal intervals of $\Ua$ of size greater than $k$ satisfy the above conditions. 
\end{theorem}

\section{Further generalizations}

Theorem \ref{sumuptheorem} provides a characterization of $k$-Naples parking functions, with $n$ cars parking in $n$ spots. We will extend such a result for some generalizations of the standard parking problem.

\begin{theorem}
	\label{lesscars}
Consider the parking problem with $n$ parking spots and $m\le n$ cars. Let $\alpha=(a_1,\dots,a_m)\in[1,n]^m$ be their preferences, and assume that all cars follow the $k$-Naples parking rule. Then all cars are able to park if and only if $\beta=(a_1,\dots,a_m,1,\dots,1)\in\pp$ is a $k$-Naples parking function.
\end{theorem}

\begin{proof}
	Consider the parking problem with $n$ cars and $n$ parking spots, with preference $\beta=(a_1,\dots,a_m,1,\dots,1)$. Suppose that the first $m$ cars are all able to park, then there is a total of $n-m$ spots left unoccupied. Since the remaining $n-m$ cars all have preference $1$, they are clearly able to park in those spots, so $\beta\in\knap$. Vice versa, if $\beta\in\knap$, then in particular the first $m$ cars are able to park.
\end{proof}

We recall the definition of $\ua(j)$, $j\in[n]$, given by \eqref{ua}:
\[
\ua(j)\coloneqq \sum_{i=j}^n \A_i -(n-j+1).
\]
Note that this definition is expressed in terms of cars having preference at least $j$. Thus, we can use this definition for the problem with $m\le n$ cars, and moreover, adding $n-m$ cars with preference $1$ does not change the value of $\ua(j)$, for all $j\in[2,n]$. Therefore, we get a straightforward generalization of Theorem \ref{sumuptheorem} for this case.

\begin{theorem}
	\label{version2sumuptheorem}
	Consider the parking problem with $n$ parking spots and $m\le n$ cars. Let $\alpha=(a_1,\dots,a_m)\in[1,n]^m$ be their preferences, and assume that all cars follow the $k$-Naples parking rule. Let $[p,q]\subseteq\Ua$ be a maximal interval, then the following are equivalent:
	\begin{enumerate}[(i)]
		\item considering only the behaviour of cars with preference at least $p$, spot $p-1$ is reached by some car;
		\item there exists $J\subseteq[m]$, $q-p+2\le\abs{J}$, such that $\tau_{p-2}(\alpha_{|_J})$ is a complete $k$-Naples parking function. 
	\end{enumerate}
	The stated properties are guaranteed to be satisfied if $\abs{[p,q]}\le k$.\\
	Furthermore, all cars are able to park if and only if all maximal intervals of $\Ua$ of size greater than $k$ satisfy the above conditions. 
\end{theorem}

\begin{proof}
	Consider $\beta=(a_1,\dots,a_m,1,\dots,1)\in\pp$, and its parking problem. Note that $\Ua=U_{\beta}$, and consider a maximal interval $[p,q]\subseteq\Ua$. When only considering cars with preference at least $p$ (which is also what is needed in point \emph{(ii)}), $\alpha$ and $\beta$ behave in the same way. As a result, the equivalence of \emph{(i)} and \emph{(ii)} follows from Theorem \ref{sumuptheorem}, as well as the fact that those conditions are satisfied if $\abs{[p,q]}\le k$.\\
	Finally, by Theorem \ref{lesscars} all cars are able to park for $\alpha$ if and only if $\beta\in\knap$, and by Theorem \ref{sumuptheorem} that holds if and only if the above conditions are satisfied for all maximal intervals of $\Ua$.\\
\end{proof}

We now generalize the concept of parking functions for an infinite number of cars.

\begin{definition}
	We say that an infinite sequence $\alpha=(a_n)_{n\ge1}$, such that $a_n\in\N\setminus\{0\}$ for all $n\ge1$ is a \emph{parking preference of infinite length}. We denote the corresponding set of parking preferences by $PP_{\infty}$.\\
	Consider the parking problem with a countably infinite number of cars, with preference $\alpha=(a_n)_{n\ge1}\in PP_{\infty}$, and a one-way street with an infinite number of spots, numbered starting from $1$. Clearly, all cars are always able to park. We denote by $\mathcal{P}_n$ the state of the problem after the first $n$ cars have parked. We define $\mathcal{P}$ as the final state of the problem, after all cars have parked, where a spot is occupied if and only if there exists $n\in\N$ such that it is occupied in $\mathcal{P}_n$. The outcome maps $\psi_k$ are defined in the usual way.\\
	Assume all cars follow the standard parking rule. We say that $\alpha=(a_n)_{n\ge1}\in PP_{\infty}$ is a \emph{parking function of infinite length} if all spots are occupied in $\mathcal{P}$, and we write $\alpha=(a_n)_{n\ge1}\in PF_{\infty}$. We define the set $PF_{\infty,k}$ in an analogous way for cars following the $k$-Naples parking rule.\\
	Let $\alpha=(a_n)_{n\ge1}\in PP_{\infty}$, and $j\ge1$. Since we cannot define $\ua(j)$ by means of the usual definition \eqref{ua}, we define it using \eqref{propua}, i.e.
	\[
	\ua(j)\coloneqq j-1-\sum_{i=1}^{j-1}\A_j.
	\]
	The set $\Ua$ is again defined by $\Ua\coloneqq\{j\in[n]\mid \ua(j)\ge1\}$. Note that all other properties of $\ua(j)$ and $\Ua$ still hold, except for \eqref{ua}. Furthermore, when we consider a maximal interval $[p,q]\subseteq\Ua$, we allow $q$ to take the value $q=\infty$, to express the case where $j\in\Ua$ for all $j\ge p$.
\end{definition}

Just like in the previous case, we will again find a characterization of $k$-Naples parking functions similar to Theorem \ref{sumuptheorem}, but we first need a few additional results.

\begin{theorem}
	\label{inftytheoprec}
	Let $\alpha\in PP_{\infty}$ and $k\ge1$. Then $\alpha\in PF_{\infty,k}$ if and only if for all $[p,q]\subseteq\Ua$ such that $[p,q]$ is a maximal interval, spot $p-1$ is occupied. In particular, if $\alpha\in PF_{\infty,k}$ then there exists $i\ge1$ such that $a_i\ge p$ and $\psi_k(c_i)=p-1$.
\end{theorem}

\begin{proof}
	The proof of this theorem is analogous to the proof of Theorem \ref{theoprec}, since we look for the minimum free spot, and find that it is a spot immediately preceding a maximal interval of $\Ua$.
\end{proof}

As previously noted, maximal intervals in $\Ua$ can be infinite. The next result shows that this cannot happen if $\alpha$ is a $k$-Naples parking function.

\begin{lemma}
	\label{noinfinitylemma}
	Let $\alpha\in PF_{\infty,k}$, and let $h\in\N\setminus\{0\}$. Then there exist $N,\overline{n}\ge h$, such that in state $\mathcal{P}_{N}$, all spots in $[1,\overline{n}]$ are occupied by cars with preference in $[1,\overline{n}]$. 
	In particular, all maximal intervals $[p,q]\subseteq\Ua$ are finite, i.e. $q<\infty$.
\end{lemma}

\begin{proof}
	Let $h\in\N\setminus\{0\}$. Since $\alpha\in PF_{\infty,k}$, all spots in $[1,h]$ are eventually occupied. Let $N\ge1$ be the smallest number such that all spots in $[1,h]$ are occupied in state $\mathcal{P}_N$. Note that clearly $N\ge h$, since at least $h$ cars have parked. Consider state $\mathcal{P}_N$, and let $\overline{n}$ be the largest number such that $[1,\overline{n}]$ are all occupied. It is clear that $\overline{n}$ exists, and that $h\le \overline{n}$. By definition of $\overline{n}$, spot $\overline{n}+1$ is unoccupied in state $\mathcal{P}_N$, so no car with preference at least $\overline{n}+1$ could have parked in $[1,\overline{n}]$. Thus, $[1,\overline{n}]$ is filled by cars with preference in $[1,\overline{n}]$.\\
	Finally, let $[p,q]\in\Ua$ be a maximal interval. For $h=p-1$, there exist $N,\overline{n}\ge p-1$ such that in state $\mathcal{P}_N$, spots $[1,\overline{n}]$ are occupied by cars with preference in $[1,\overline{n}]$. Then 
	\[
	\ua(\overline{n}+1)=\overline{n}-\sum_{i=1}^{\overline{n}}\A_i\le\overline{n}-\overline{n}=0.
	\]
	Thus, since $\overline{n}+1\ge p$, and $\ua(j)\ge1$ for all $j\in[p,q]$, we get $q\le\overline{n}$, so $[p,q]$ cannot be an infinite maximal interval.
\end{proof}

\begin{corollary}
	\label{subsetcoroll}
	Let $\alpha\in PF_{\infty,k}$, and $h\ge1$. Then there exists a finite subset $J\subset\N\setminus\{0\}$ such that $\abs{J}\ge h$, and $\alpha_{|_J}\in PF_{\abs{J},k}$.
\end{corollary}

\begin{proof}
	Let $h\ge1$. By Lemma \ref{noinfinitylemma}, there exist $N,\overline{n}\ge h$ such that in state $\mathcal{P}_N$, spots $[1,\overline{n}]$ are all occupied by cars with preference in $[1,\overline{n}]$. Let 
	\[
	J=\{i\in[1,N]\mid \psi_k(c_i)\in[1,\overline{n}]\},
	\]
	then $\abs{J}=\overline{n}$, and trivially $\alpha_{|_J}\in PF_{\overline{n},k}$.
\end{proof}

\begin{theorem}
	\label{inftysumuptheorem}
	Let $\alpha\in PP_{\infty}$, $k\ge1$, and assume that all cars follow the $k$-Naples parking rule. Let $[p,q]\subseteq\Ua$ be a maximal interval, then the following are equivalent:
	\begin{enumerate}[(i)]
		\item considering only the behaviour of cars with preference at least $p$, spot $p-1$ is reached by some car;
		\item there exists a finite subset $J\subset\N\setminus\{0\}$, $q-p+2\le\abs{J}$, such that $\tau_{p-2}(\alpha_{|_J})$ is a complete $k$-Naples parking function. 
	\end{enumerate}
	The stated properties are guaranteed to be satisfied if $\abs{[p,q]}\le k$, and are certainly not satisfied if $\abs{[p,q]}=\infty$.\\
	Furthermore, $\alpha$ is a $k$-Naples parking function if and only if all maximal intervals of $\Ua$ of size greater than $k$ satisfy the above conditions. 
\end{theorem} 

\begin{proof}
	Let $\alpha\in PP_{\infty}$, and $[p,q]\subseteq\Ua$ be a maximal interval. Consider only the behaviour of cars with preference at least $p$, and suppose that spot $p-1$ is occupied by a car $c_m$. Let $n+1$ be the smallest unoccupied spot larger than $p-1$. Note that since $\A_{p-1}=0$, $a_m\ge p$, and thus $n+1>p$. Let $H=\{i\in[m]\mid a_i\in[p,n]\}$, $h=\abs{H}$.  Define $\beta=(b_1,\dots,b_n)\in\pp$ by $\beta_{|_{[1,h]}}=\alpha_{|_H}$, and $b_i=1$ for all $i\ge h$, and consider the finite parking problem with $n$ spots for $\beta$. Since all cars with preference at least $p-1$ (i.e. the first $h$ cars) occupy exactly spots $[p-1,n]$, we get $\ua(p-1)=0$, and since $\abs{\beta}_{p-1}=0$, we have $\ua(p)=1$. Moreover, since all other cars have preference $1$, clearly $\beta\in\knap$. Thus, we get \emph{(ii)} by Theorem \ref{sumuptheorem}.\\
	Vice versa, suppose there exists a finite subset $J\subset\N\setminus\{0\}$, $q-p+2\le\abs{J}$, such that $\tau_{p-2}(\alpha_{|_J})$ is a complete $k$-Naples parking function. Note that $\tau_{p-2}(\alpha_{|_J})\in PP_{\abs{J}}$, so we get $a_i\in [p-1,p-2+\abs{J}]$ for all $i\in J$. Consider the parking problem with $p-2+\abs{J}$ spots and $\abs{J}$ cars, with preference $\alpha_{|_J}$. Since $\tau_{p-2}(\alpha_{|_J})$ is complete, we have $U_{\alpha_{|_J}}=[p,p-2+\abs{J}]$, and \emph{(i)} follows by Theorem \ref{version2sumuptheorem}.\\
	If $\abs{[p,q]}=\infty$, then \emph{(ii)} certainly fails, as it requires a finite subset $J$ such that $q-p+2\le\abs{J}$. If $\abs{[p,q]}\le k$, then since $\ua(q)\ge1$ and $\ua(q+1)\le0$, we get $\sum_{i=1}^{q-1}\A_i\le q-2$, and $\sum_{i=1}^q\A_i\ge q$. Let $\tilde{H}$ be the set of all indices such that $a_i\le q-1$, as well as the first $q-\sum_{i=1}^{q-1}\A_i$ indices such that $a_i=q$. Then $\alpha_{|_{\tilde{H}}}\in PP_q$, and $U_{\alpha_{|_{\tilde{H}}}}=\Ua\cap[1,q]$, and by Theorem \ref{sumuptheorem} the above conditions are satisfied.\\
	Finally, suppose $\alpha\notin PF_{\infty,k}$, then by Theorem \ref{inftytheoprec} there exists a maximal interval $[p,q]\in\Ua$ such that spot $p-1$ is never occupied, thus \emph{(i)} fails for that maximal interval. Vice versa, suppose $\alpha\in PF_{\infty,k}$, and consider a maximal interval $[p,q]\in\Ua$. By Theorem \ref{inftytheoprec}, spot $p-1$ is occupied by a car $c_m$ with preference at least $p$. In particular, all cars with preference at most $p-1$ coming before car $c_m$ park in some spots in $[1,p-2]$, so spot $p-1$ is reached even when only taking into account cars with preference at least $p$, and we get \emph{(i)}.
\end{proof}

For our last generalization, we take the parking problem with countably infinite cars and spots, and we allow cars to drive backwards as far as they need to find a free spot.

\begin{definition}
	Let $\alpha=(a_n)_{n\ge1}\in PP_{\infty}$, we say that a car $c_j$ follows the \emph{$\infty$-Naples parking rule} if it parks according to the generic Naples parking rule, except that there is no restriction on how many steps it can take backwards. In other words, car $c_j$ follows the $\infty$-Naples parking rule if it follows the $a_j$-Naples parking rule. \\
	Assuming all cars follow the $\infty$-Naples parking rule, if all spots are eventually occupied we call $\alpha$ an \emph{$\infty$-Naples parking function}, and we write $\alpha\in PF_{\infty,\infty}$.
\end{definition}

\begin{remark}
	When we defined $k$-Naples parking functions, we remarked that, clearly, any parking preference of length $n$ is also a $n$-Naples parking function. In other words, $PP_n=PF_{n,n}$ for all $n\ge1$. The same cannot be said for parking preference of infinite length. For example, let $\alpha=(a_n)_{n\ge1}$ be defined by $a_n=n+1$ for all $n\ge1$, i.e. $\alpha=(2,3,4,\dots)$. Trivially, each car $c_n$ parks in spot $a_n$, so spot $1$ is never filled by any car, and $\alpha\in PP_{\infty}\setminus PF_{\infty,\infty}$. Note that in this case, we get $\Ua=[2,\infty]$: in particular, $\Ua$ contains an infinite maximal interval. We prove in the next theorem that having no infinite maximal intervals in $\Ua$ is a necessary and sufficient condition for $\infty$-Naples parking functions.
\end{remark}

\begin{theorem}
	Let $\alpha\in PP_{\infty}$. Then $\alpha\in PF_{\infty,\infty}$ if and only if all maximal intervals of $\Ua$ are finite.
\end{theorem}

\begin{proof}
	Let $\alpha=(a_n)_{n\ge1}\in PP_{\infty}$. Let $[p,q]\subseteq\Ua$ be a finite maximal interval of $\Ua$, in particular $q<\infty$, and let $k=\abs{[p,q]}$. Since $\ua(q+1)\le0$, by definition $\sum_{i=1}^q\A_i\ge q$. Let $H=\{i\ge1\mid a_i\le q\}$, and consider $\alpha_{|_H}$. If $\abs{H}=h<\infty$, consider the parking problem with $h$ spots and $h$ cars, with preference $\alpha_{|_H}$. We get $U_{\alpha_{|_H}}=\Ua\cap[1,q]$. By Theorem \ref{sumuptheorem}, there exists a finite subset $J\subseteq H$, $q-p+2\le\abs{J}$, such that $\tau_{p-2}(\alpha_{|_J})$ is a complete $k$-Naples parking function. If instead $\abs{H}=\infty$, consider the problem with an infinite number of parking spots, and preference $\alpha_{|_H}$, and we get the same result using Theorem \ref{inftysumuptheorem}.\\
	Thus, there exists a finite subset $J\subseteq H\subseteq \N\setminus\{0\}$ such that it satisfies condition \emph{(ii)} in Theorem \ref{inftysumuptheorem}, using the $k$-Naples parking rule. Let $c_m$ be the last car in the sequence $\alpha_{|_J}$ (thus, $J\subseteq[1,m]$), and consider state $\mathcal{P}_m$ in the original problem. Note that for the first $m$ cars, following the $\infty$-Naples parking rule is the same as following the $m$-Naples parking rule. Moreover, since $\tau_{p-2}(\alpha_{|_J})$ is a $k$-Naples parking function, it is also a $m$-Naples parking function (recall that $k=\abs{[p,q]}<\abs{J}\le m$). In the same way as in the proof of Theorem \ref{theosuff}, this implies that spot $p-1$ is going to be occupied by a car among the first $m$ cars. Thus, spot $p-1$ is going to be occupied in the original problem as well.\\
	Theorem \ref{inftytheoprec} can easily be extended to the $\infty$-Naples parking rule as well, so if all maximal intervals are finite, we get that $\alpha\in PF_{\infty,\infty}$. Vice versa, if $\alpha\in PF_{\infty,\infty}$ we can extend Lemma \ref{noinfinitylemma} and obtain that all maximal intervals need to be finite.
\end{proof}

\chapter{Enumerative results}
\label{chapterenumeration}

\section{Enumeration of permutation invariant parking functions}

In this section, as well as in the following one, we cover the enumeration of permutation invariant $k$-Naples parking functions, through the characterization provided by Corollary \ref{perminv}. We approach this problem by expressing generic permutation invariant $k$-Naples parking functions in terms of other parking preferences with progressively more restrictions.

\begin{definition}
	We define $T_{n,k}$ as the number of parking preferences $\alpha\in\pp$ such that for each maximal interval $[p,q]\subseteq\Ua$, $q-p+1\le k$. That is, $T_{n,k}$ is the number of $k$-Naples parking functions of length n that are permutation invariant.\\
	We define $t_{n,k}$, $n\ge1$, as the number of parking preferences $\alpha\in\pp$ such that for each maximal interval $[p,q]\subseteq\Ua$, $q-p+1\le k$, and such that $2\in\Ua$. Note that $t_{1,k}=0$. Moreover, let $t_{0,k}=1$.
\end{definition}

Clearly, if $\alpha\in\pf$, then $\alpha\in\knap$, and is permutation invariant, since all rearrangements of a parking function are still parking functions. Suppose that $\alpha\in\knap$ is permutation invariant, such that $\alpha\notin\pf$. Then $\Ua$ is not empty; let $h=\min\Ua$, $H=\{i\in[n]\mid a_i\in[h-1,n]\}$ and $H^c=[1,n]\setminus H$. Note that $\abs{H}\ge2$.\\
By Proposition \ref{propUa}, $\ua(h-1)=0$. Consider $\alpha_{|_{H^c}}$: by Lemma \ref{translatelemma}, $\alpha_{|_{H^c}}\in PP_{h-2}$, and $U_{\alpha_{|_{H^c}}}=\Ua\cap[1,h-2]=\varnothing$, so $\alpha_{|_{H^c}}$ is a parking function of length $h-2$.\\
Moreover, again by Lemma \ref{translatelemma}, $\tau_{h-2}(\alpha_{|_H})\in PP_{n-h+2}$, and $U_{\tau_{h-2}(\alpha_{|_H})}=\tau_{h-2}(\Ua\cap[h-1,n])=\tau_{h-2}(\Ua)$. Since all maximal intervals in $\Ua$ are of size $k$ at most, the same holds for $U_{\tau_{h-2}(\alpha_{|_H})}$, thus $\tau_{h-2}(\alpha_{|_H})\in PF_{n-h+2,k}$. Furthermore, $u_{\tau_{h-2}(\alpha_{|_H})}(2)=\ua(h)\ge1$, thus $2\in U_{\tau_{h-2}(\alpha_{|_H})}$.\\

Hence, a permutation invariant $\alpha\in\knap\setminus\pf$ is uniquely determined by choosing:
\begin{itemize}
	\item an integer $j\in[0,n-2]$;
	\item a parking function $\alpha_1\in PF_j$;
	\item a permutation invariant $k$-Naples parking function $\alpha_2\in PF_{n-j,k}$ such that $2\in U_{\alpha_2}$;
	\item a set $J\subset[n]$ such that $\abs{J}=j$; that is, choosing $j$ cars to assign $\alpha_1$ to, and assigning $\alpha_2$ (suitably translated) to the remaining $n-j$ cars. 
\end{itemize}

Thus, we get
\begin{equation}
	\label{formulaT}
	T_{n,k}=\sum_{i=0}^n \binom{n}{i} (i+1)^{i-1} t_{n-i,k} .
\end{equation}
Note that the case $i=n-1$ is excluded since $t_{1,k}=0$, and the case $i=n$ ($t_{0,k}=1$ and $\binom{n}{0}=1$) corresponds to standard parking functions.

\begin{definition}
	We define $\Teq$, $k\ge1$, as the number of $\alpha\in\pp$ such that $\Ua=[2,k+1]$. Note that $\Theta_{1,k}^= =\Theta_{0,k}^= =0$ for all $k\ge1$.\\
	We define $\Tnk\coloneqq\sum_{i=1}^k \Theta_{n,i}^=$, i.e. $\Tnk$ is the number of parking preferences of length $n$ such that $\Ua=[2,j]$ for some $j=2,\dots,k+1$.
\end{definition}

Let $\alpha\in\knap$ be permutation invariant, and such that $2\in\Ua$. Then $2$ belongs to the maximal interval $[2,j]\subset\Ua$, such that $j\le k+1$. Suppose that $\Ua\ne[2,j]$, and let $h=\min(\Ua\setminus [2,j])$, $H=\{i\in[n]\mid a_i\in[h-1,n]\}$ and $H^c=[1,n]\setminus H$.\\
Then we again get that $\tau_{h-2}(\alpha_{|_H})$ is a permutation invariant $k$-Naples parking function, such that $2\in U_{\tau_{h-2}(\alpha_{|_H})}$. Moreover, $\alpha_{|_{H^c}}\in PP_{\abs{H^c}}$, and $U_{\alpha_{|_{H^c}}}=[2,j]$ (which also implies $\alpha_{|_{H^c}}\in PF_{\abs{H^c},k}$).\\
Thus, we get the recursive relation
\begin{equation}
	\label{formulat}
	t_{n,k}=\sum_{i=0}^n \binom{n}{i} \Theta_{i,k}^{\le} t_{n-i,k},
\end{equation}
where if $i=0$, $\Theta_{0,k}^{\le}=0$, so this is not an infinite recursion.

Note that using \eqref{formulat}, equation \eqref{formulaT} can be reformulated in the following way: let
\begin{equation}
	W=\biggl\{\beta=(b_1,b_2,\dots,b_m) \mid \sum_{i=1}^m b_i=n,\, b_i\ge2 \quad\! \textit{for all} \quad\! i\ne1\biggr\},
\end{equation}
then \eqref{formulaT} and \eqref{formulat} can be combined into

\begin{equation}
	\label{count Tnk}
	T_{n,k}=\sum_{\beta\in W} \binom{n}{b_0,b_1,\dots,b_m} (b_0+1)^{b_0-1} \prod_{i=1}^m \Theta_{b_i,k}^{\le}.
\end{equation}

This formula has the following interpretation: a permutation invariant $k$-Naples parking function $\alpha$ is uniquely determined by choosing:
\begin{itemize}
	\item the number $m$ of maximal intervals in $\Ua$;
	\item if $m\ne0$, the lower extreme of each maximal interval in $\Ua$: these are $b_0+2,b_0+b_1+2,\dots,\sum_{i=0}^{m-1}b_i +2$ (and $b_m=n-\sum_{i=1}^{m-1}$);
	\item if $b_0\ne0$, a parking function $\alpha_0\in PF_{b_0}$;
	\item if $m\ne0$, for each $1\le i\le m$, a permutation invariant $k$-Naples parking function $\alpha_i\in PF_{b_i,k}$ such that $U_{\alpha_i}=[2,j_i]$, $2\le j_i\le k+1$;
	\item sets $J_0,J_1,\dots,J_m\subseteq[n]$ such that $\abs{J_i}=b_i$ for all $i=0,\dots,m$, and $J_i\cap J_j=\varnothing$ for all $i\ne j$: thus, $\alpha_{|_{J_i}}=\tau_{-(\sum_{j=0}^{i-1}b_j)}(\alpha_i)$.\\
\end{itemize}

We now take a closer look at the objects $\Theta_{n,k}^=$. We start by showing the following symmetry.

\begin{lemma}
	\label{thetasymmetry}
	Let $n\ge2$ and $k\in[1,n-1]$. Then $\Theta_{n,k}^= =\Theta_{n,n-k}^=$.
\end{lemma}

\begin{proof}
	Define a function $\varphi$ such that for all $\alpha=(a_1,\dots,a_n)$, $\varphi(\alpha)=(n-a_1+2,\dots,n-a_n+2)$. Let $\alpha\in\pp$ be such that $\Ua=[2,k+1]$, and consider $\varphi(\alpha)$. Note that by Proposition \ref{propUa}, $\A_1=0$, so $\varphi(\alpha)\in\pp$. Observe that $\abs{\varphi(\alpha)}_i=\abs{\alpha}_{n-i+2}$, so for all $j\in[2,n]$
	\begin{multline*}
	u_{\varphi(\alpha)}(j)=\sum_{i=j}^n\abs{\varphi(\alpha)}_i-(n-j+1)=\sum_{i=j}^n\abs{\alpha}_{n-i+2}-(n-j+1)\\
	=\sum_{i'=1}^{n-j+2}\A_{i'}-(n-j+2)+1=-\ua(n-j+3)+1,
	\end{multline*}
	where \eqref{propua} was used in the last equality. Thus, $u_{\varphi(\alpha)}(j)\ge1$ if and only if $\ua(n-j+3)\le0$, meaning $j\in U_{\varphi(\alpha)}$ if and only if $(n-j+3)\notin\Ua=[2,k+1]$: since $n-j+3\ge k+2$ if and only if $j\le n-k+1$, we get $U_{\varphi(\alpha)}=[2,n-k+1]$. \\
	In the same way, if $\alpha\in\pp$ is such that $\Ua=[2,n-k+1]$, then $U_{\varphi(\alpha)}=[2,k+1]$, and since $\varphi$ is clearly an involution, we obtain $\Theta_{n,k}^= =\Theta_{n,n-k}^=$.\\
\end{proof}

\begin{definition}
	Let $n,m,h\ge1$ and $k\ge0$. We define $\tnk(m,h)$ as the number of $\alpha\in\pp$ that satisfy the following conditions: $\Ua=[2,k+1]$ (if $k=0$, $\Ua=\varnothing$), $\A_i=0$ for all $i<h$, and $\A_h=m$.
\end{definition}

Using the previous definition, we can express $\Theta^=_{n,k}$ ($n,k\ge1$) in terms of $\vartheta_{n,k}(m,h)$. Specifically:
\begin{equation}
	\label{recTheta}
\Theta^=_{n,k}=\sum_{m\ge1}\sum_{h\ge1}\vartheta_{n,k}(m,h),
\end{equation}
which immediately follows by noting that any parking preference $\alpha$ such that $\Ua=[2,k+1]$, satisfies $\A_i=0$ for all $i<h$, and $\A_h=m$ for some $m,h\ge1$.

\begin{lemma}
	\label{propertytheta}
	Let $n,m,h\ge1$ and $n> k\ge0$. Then:
	\begin{enumerate}[(i)]
		\item $\vartheta_{n,k}(m,h)=0$ for all $h>k+1$, in particular $\vartheta_{n,0}(m,h)=0$ for all $h\ne1$;
		\item if $k\ge1$, then $\tnk(m,1)=0$;
		\item if $k\ge1$, then $\tnk(m,k+1)=0$ for all $m<k+1$;
		\item if $k\ge1$ and $h\in[2,k]$, then $\tnk(m,h)=0$ for all $m\ge h$;
		\item if $k,m,h$ do not satisfy cases \emph{(i),(ii),(iii),(iv)}, then $\tnk(m,h)\ne0$.
	\end{enumerate}
\end{lemma}

\begin{proof}
	Let $\alpha\in\pp$, if $\alpha\in\pf$ then $\A_1\ne0$, so $\vartheta_{n,0}(m,h)=0$ for all $h>1$. From now on, suppose $k\ge1$, and so $\Ua=[2,k+1]$. By Proposition \ref{propUa} $\A_{k+1}\ne0$, so $\vartheta_{n,k}(m,h)=0$ for all $h>k+1$, thereby proving \emph{(i)}. By the same proposition, $\A_1=0$, and \emph{(ii)} follows.\\
	If $h=k+1$, then $\A_i=0$ for all $i\le k$. Thus
	\[
	\ua(k+2)=k+1-\sum_{i=1}^{k+1}\A_i=k+1-\A_{k+1}\le0,
	\]
	where the last inequality is due to $(k+2)\notin\Ua$. Hence, $\A_{k+1}\ge k+1$, so $\tnk(m,k+1)=0$ for all $m<k+1$, and we have \emph{(iii)}.\\
	If $h\in[2,k]$, then $(h+1)\in\Ua$, therefore
	\[
	\ua(h+1)=h-\sum_{i=1}^h\A_i=h-\A_h\ge1
	\]
	so $\A_h \le h-1$, and \emph{(iv)} follows.\\
	Finally, suppose that $k,m,h$ do not satisfy cases \emph{(i),(ii),(iii),(iv)}. If $k=0$, then $h=1$, and clearly $\alpha$ such that $a_i=1$ if $i\in[1,m]$, and $a_j=2$ otherwise, is one of the desired parking functions. If instead $k\ge1$, then it is easy to check that, if $h\in[2,k]$ then one of the desired preferences is $\alpha$ such that $a_i=h$ for $i\in[1,m]$, and $a_j=k+1$ otherwise; if $h=k+1$, then we can choose $\alpha$ such that $a_i=h=k+1$ for $i\in[1,m]$ and $a_j=k+2$ otherwise. Thus, we have \emph{(v)}.
\end{proof}

\begin{remark}
	By Lemma \ref{propertytheta}, we can rewrite \eqref{recTheta} by excluding all terms that equal zero, and obtain
	\begin{equation}
		\label{recTheta2}
	\Theta_{n,k}^= =\sum_{h=2}^k\sum_{m=1}^{h-1}\tnk(m,h) + \sum_{m=k+1}^n\tnk(m,k+1).
	\end{equation}
	Furthermore, formulas \eqref{count Tnk} and \eqref{recTheta2}, as well as the definition of $\Tnk$, can be put together to express $T_{n,k}$ just in terms of $\tnk(m,h)$, although the resulting formula is a bit too unwieldy to be of any use.
\end{remark}

The following results provide a recursive relation for $\tnk(m,h)$, thereby completing the enumeration of permutation invariant $k$-Naples parking functions.

\begin{proposition}
	Let $n\ge m\ge 1$. Then 
	\begin{equation}
		\label{count0}
		\vartheta_{n,0}(m,1)=\binom{n-1}{m-1}n^{n-m}.
	\end{equation}
\end{proposition}

\begin{proof}
	The value $\vartheta_{n,0}(m,1)$ counts the number of parking functions of length $n$ such that $m$ cars have preference $1$, which is given by Lemma \ref{1preference}.
\end{proof}

\begin{proposition}
	Let $n>k\ge1$, and $h\in[n]$. Then for $m=1$:
	\begin{equation}
		\label{ric1}
		\tnk(1,h)=n\sum_{h'\ge h}\sum_{m'}\vartheta_{n-1,k-1}(m',h');
	\end{equation}
	otherwise, if $m\ge2$:
	\begin{equation}
		\label{ric2}
		\tnk(m,h)=\frac{n}{m}\vartheta_{n-1,k-1}(m-1,h-1).
	\end{equation}
\end{proposition}

\begin{proof}
	Let $\alpha\in\knap$ be such that $\Ua=[2,k+1]$, and $h\in[n]$ be such that $\A_i=0$ for all $i<h$. Let $i_h$ be an index such that $a_{i_h}=h$, and let $J=\{j\in[n]\mid j\ne i_h\}$. Consider $\tau_1(\alpha_{|_J})$, note that $\abs{\tau_1(\alpha_{|_J})}_i=\A_{i+1}$ for all $i\ne h-1$, and $\abs{\tau_1(\alpha_{|_J})}_{h-1}=\A_h-1$. For all $j\in[h,n-1]$
	\[
	u_{\tau_1(\alpha_{|_J})}(j)=\sum_{i=j}^{n-1}\abs{\tau_1(\alpha_{|_J})}_i-(n-1-j+1)=\sum_{i=j+1}^n\A_i -(n-(j+1)+1)=\ua(j+1),
	\]
	By statement \emph{(i)} in Lemma \ref{propertytheta}, $h\le k+1$, thus $[h,k]\subseteq U_{\tau_1(\alpha_{|_J})}$, and $i\notin U_{\tau_1(\alpha_{|_J})}$ for all $i>k$. Furthermore, for all $j'\in[2,h-1]$
	\begin{multline*}
	u_{\tau_1(\alpha_{|_J})}(j')=\sum_{i=j'}^{n-1}\abs{\tau_1(\alpha_{|_J})}_i-(n-1-j'+1)\\
	=\sum_{i=j'+1}^n\A_i-1-(n-(j'+1)+1)=\ua(j'+1)-1.
	\end{multline*}
	Observe that for all $i\in[3,h]$, $\ua(i)\ge2$: if that was not the case, then since $\A_{i-1}=0$ we would have, using \eqref{propua3},
	\[
	\ua(i-1)=\ua(i)+\A_{i-1}-1\le0,
	\]
	which contradicts the fact that $(i-1)\in[2,h-1]\subset[2,k+1]=\Ua$. Thus $u_{\tau_1(\alpha_{|_J})}(j')=\ua(j'+1)-1\ge1$ for all $j'\in[2,h-1]$, so $[2,h-1]\subset U_{\tau_1(\alpha_{|_J})}$. Hence, $U_{\tau_1(\alpha_{|_J})}=[2,k]$.\\
	Conversely, let $\beta\in PP_{n-1}$ be a permutation invariant $(k-1)$-Naples parking function such that $U_{\beta}=[2,k]$, $\abs{\beta}_i=0$ for all $i\le h-2$, and $J=\{i\in[n]\mid i\ne j\}$ for some $j\in[n]$. In the same way, we can prove that $\alpha\in\pp$ constructed such that $\alpha_{|_J}=\tau_{-1}(\beta)$ and $a_j=h$, is a permutation invariant $k$-Naples parking function, such that $\A_i=0$ for all $i\le h-1$.\\ 
	Now, suppose that $\abs{\beta}_{h-1}=0$, then $j$ can be arbitrarily chosen in $n$ ways, and $\A_h=1$, so we never obtain any $\alpha$ more than once. Since all $\beta$ as above are counted by all the possible $\vartheta_{n-1,k-1}(m',h')$ such that $h'\ge h$, we get \eqref{ric1}.\\
	On the other hand, suppose that $\abs{\beta}_{h-1}=m-1$, $m\ge2$. Then $j$ can be chosen in $n$ ways, but since $\A_h=m\ge2$, each $\alpha$ is counted $m$ times. The number of such $\beta$ is $\vartheta_{n-1,k-1}(m-1,h-1)$, and \eqref{ric2} immediately follows.\\
\end{proof}

\begin{remark}
	Computing a few values of $\Theta^=_{n,k}$, $n>k\ge1$, gives the following:
	\[
	\begin{array}{cccccccc}
		\toprule
			&k=1	&k=2	&k=3	&k=4	&k=5	&k=6	&k=7	\\
		\midrule
		n=2	&1		&		&		&		&		&		&		\\
		n=3	&4		&4		&		&		&		&		&		\\
		n=4	&27		&21		&27		&		&		&		&		\\
		n=5	&256	&176	&176	&256	&		&		&		\\
		n=6	&3125	&1995	&1765	&1995	&3125	&		&		\\
		n=7	&46656	&28344	&23304	&23304	&28344	&46656	&		\\
		n=8	&823543	&482825	&378007	&351337	&378007	&482825	&823543	\\
		%n=9	&16777216 &9576160 &7238944 &6407680 &6407680 &7238944 &9576160 &16777216\\
		\bottomrule
	\end{array}.
	\]
	Note that, as we saw in Lemma \ref{thetasymmetry}, the above triangle is symmetric. Furthermore, the sequences $\Theta^=_{n,1}$ and $\Theta^=_{n,n-1}$ satisfy $\Theta^=_{n,1}=\Theta^=_{n,n-1}=(n-1)^{n-1}$: this will be proven in Lemma \ref{count1} and Corollary \ref{completetot}. As for the rest, they do not seem to appear anywhere on OEIS, the Online Encyclopedia of Integer Sequences \cite{oeis}. Values for $\Tnk$ also seem to be unknown sequences, except for $\Theta_{n,n-1}^{\le}$. For  $2\le n\le8$, we get $\Theta_{n,n-1}^{\le}=1, 8, 75, 864, 12005, 196608, 3720087$, which follows sequence \href{https://oeis.org/A071720}{A071720}, with an offset, with corresponding formula $(n-1)(n+1)^{n-2}$. It could be interesting to prove this formula for $\Theta_{n,n-1}^{\le}$, as it would imply that the number of parking preferences of length $n$ such that $\Ua=[2,k+1]$ for some $1\le k\le n-1$, is equal to the number of spanning trees in $K_{n+1}-e$, the complete graph on $n+1$ nodes minus an edge (that is, the objects counted by $A071720$).\\
\end{remark}

\begin{lemma}
	\label{count1}
	Let $n\ge2$, then $\vartheta_{n,1}(m,h)$ is nonzero if and only if $h=2$ and $m\in[2,n]$. Furthermore,
	\begin{equation}
		\label{count2}
		\Theta_{n,1}^= = (n-1)^{n-1}.
	\end{equation}
\end{lemma}

\begin{proof}
	Using Lemma \ref{propertytheta}, \emph{(i)} implies $\vartheta_{n,1}(m,h)=0$ if $h>2$, \emph{(ii)} implies $\vartheta_{n,1}(m,1)=0$, and \emph{(iii)} implies $\vartheta_{n,1}(m,2)=0$ if $m<2$. Finally, \emph{(v)} implies that in all other cases $\vartheta_{n,1}(m,h)$ is nonzero.\\
	Now, by \eqref{count0} and \eqref{ric2}, we have
	\begin{multline*}
		\Theta_{n,1}^= =\sum_{m,h}\vartheta_{n,1}(m,h)=\sum_{m=2}^n\vartheta_{n,1}(m,2)
		=\sum_{m=2}^n \frac{n}{m}\vartheta_{n-1,0}(m-1,1)\\
		=\sum_{m=2}^n \frac{n}{m}\binom{n-2}{m-2}(n-1)^{n-m}
		=\sum_{m'=0}^{n-2}\frac{n}{m'+2}\binom{n-2}{m'}(n-1)^{n-m'-2}\\
		=(n-1)^{n-1}\sum_{m'=0}^{n-2}\frac{n}{m'+2}\binom{n-2}{m'}(n-1)^{-(m'+1)}.
	\end{multline*}
	Let $n'=n-2$, and define
	\begin{equation}
	F(n')=\sum_{m'=0}^{n'}\binom{n'}{m'}\frac{n'+2}{m'+2}(n'+1)^{-(m'+1)}.
	\end{equation}
	To prove \eqref{count2}, it is sufficient to prove that $F(n')=1$ for all $n'\ge0$. For simplicity, for the rest of this proof we will use $n,m$ instead of $n',m'$.\\
	For $x\ge0$, define:
	\[
		A(x)=\int_0^x(y+1)^n dy;\qquad B(x)=\int_0^x A(y) dy. 
	\]
	Given Newton's formula $(y+1)^n=\sum_{m=0}^n\binom{n}{m}y^m$, we can write
	\begin{multline}
		\label{eqA1}
		A(x)=\int_0^x (y+1)^n dy=\int_0^x \Biggl(\sum_{m=0}^n\binom{n}{m}y^m \Biggr)dy\\
		=\sum_{m=0}^n\binom{n}{m}\int_0^x y^m dy=\sum_{m=0}^n\binom{n}{m}\frac{1}{m+1}x^{m+1}.
	\end{multline}
	On the other hand, simply solving the integral yields
	\begin{equation}
		\label{eqA2}
		A(x)=\int_0^x (y+1)^n dy= \frac{(x+1)^{n+1}}{n+1}-\frac{1}{n+1}.
	\end{equation}
	As for $B(x)$, using \eqref{eqA1} gives
	\begin{multline}
		\label{eqB1}
		B(x)=\int_0^x A(y)dy=\int_0^x\Biggl( \sum_{m=0}^n\binom{n}{m}\frac{1}{m+1}y^{m+1}\Biggr)dy\\
		=\sum_{m=0}^n\binom{n}{m}\frac{1}{m+1}\int_0^x y^{m+1}dy
		=\sum_{m=0}^n\binom{n}{m}\frac{1}{(m+1)(m+2)}x^{m+2},
	\end{multline}
	whereas \eqref{eqA2} yields
	\begin{equation}
		\label{eqB2}
		B(x)=\int_0^x\biggl(\frac{(y+1)^{n+1}}{n+1}-\frac{1}{n+1}\biggr)dy
		=\frac{(x+1)^{n+2}}{(n+1)(n+2)}-\frac{1}{(n+1)(n+2)}-\frac{x}{n+1}.
	\end{equation}
	Finally, evaluating the expression $A(x)-(n+1)B(x)$ in $x=(n+1)^{-1}$, and using \eqref{eqA1} and \eqref{eqB1}, we get
	\begin{align*}
	\bigl(A(x)-(n+1)B(x)\bigr)_{|_{x=(n+1)^{-1}}}
	=&\sum_{m=0}^n\binom{n}{m}\frac{1}{m+1}(n+1)^{-(m+1)}\\
	&-\sum_{m=0}^n\binom{n}{m}\frac{1}{(m+1)(m+2)}(n+1)^{-(m+1)}\\
	&=\sum_{m=0}^n\binom{n}{m}\frac{1}{m+2}(n+1)^{-(m+1)}=\frac{F(x)}{n+2}.
	\end{align*}
	Hence, by \eqref{eqA2} and \eqref{eqB2} we have
	\begin{multline*}
		F(x)=(n+2)\bigl(A(x)-(n+1)B(x)\bigr)_{|_{x=(n+1)^{-1}}}\\
		=(n+2)\Biggl(\frac{\bigl(\frac{n+2}{n+1}\bigr)^{n+1}}{n+1}-\frac{1}{n+1}- \frac{\bigl(\frac{n+2}{n+1}\bigr)^{n+2}}{n+2}+\frac{1}{n+2}+\frac{1}{n+1}\Biggr)\\
		=(n+2)\Biggl(\frac{(n+2)^{n+1}}{(n+1)^{n+2}}- \frac{(n+2)^{n+1}}{(n+1)^{n+2}}+\frac{1}{n+2}\Biggr)=1
	\end{multline*}
\end{proof}

\begin{corollary}
	\label{completetot}
	Let $n\ge2$. The number of parking preferences of length $n$ that are complete is $(n-1)^{n-1}$.
\end{corollary}

\begin{proof}
	The number of complete parking preferences of length $n$ is $\Theta_{n,n-1}^=$. By Lemmas \ref{thetasymmetry} and \ref{count1}, $\Theta_{n,n-1}^= =\Theta_{n,1}^= =(n-1)^{n-1}$
\end{proof}

We end this section by providing the first values of $T_{n,k}$, that is, the number of permutation invariant $k$-Naples parking functions of length $n$:
\[
\begin{array}{cccccccc}
	\toprule
	&k=1	&k=2	&k=3	&k=4	&k=5	&k=6	&k=7\\
	\midrule
	n=2	&\mathbf{4}	&	&	&	&	&	&\\
	n=3	&23	&\mathbf{27}	&	&	&	&	&\\
	n=4	&192	&229	&\mathbf{256}	&	&	&	&\\
	n=5	&2077	&2558	&2869	&\mathbf{3125}	&	&	&\\
	n=6	&27808	&35154	&40000	&43531	&\mathbf{46656}	&	&\\
	n=7	&444411	&572470	&662519	&726668	&776887	&\mathbf{823543}	&\\
	n=8	&8266240	&10815697	&12693504	&14055341	&15097600	&15953673	&\mathbf{16777216}\\
	\bottomrule
\end{array},
\]
where the numbers in bold are the number of parking preferences of length $n$, $\abs{\pp}=n^n$, and clearly for all $k\ge n-1$, $T_{n,k}=\abs{\pp}$.
Again, OEIS does not contain any of these sequences, except for the first diagonal, which is $n^n$, and the second one, which is $n^n-(n-1)^{n-1}$, i.e. the number of parking preferences minus the number of complete parking preferences of length $n$.

\section[Enumeration of complete $k$-Naples parking functions]{Enumeration of complete $k$-Naples parking\\ functions}

\begin{definition}
	Let $\unk(m)$ ($k\ge1$) be the number of $\alpha\in\knap$ such that $\ua(i)\ge1$ for all $i\in[m+1,n]$, and $\ua(j)\ge0$ for all $j\in[1,m]$. We also define $\unk^0(m)$ as the number of $\alpha\in\knap$ as above, with the additional condition that $\ua(m)=0$.
\end{definition}

\begin{remark}
	Note that since $\ua(1)=0$ for all $\alpha\in\pp$, the number of $k$-Naples parking functions that are complete is $\unk(1)=\unk^0(1)$. Furthermore, if $n\le k+1$, Corollary \ref{completetot} implies that $\unk(1)=(n-1)^{n-1}$.
\end{remark}

Consider a $k$-Naples parking function $\alpha=(a_1,\dots,a_n,a_{n+1})\in PF_{n+1,k}$ of length $n+1$ that is complete, and such that $n+1\ge k+1$. By Lemma \ref{lemmartl}, $a_{n+1}\le k+1$, moreover $a_{n+1}\ne1$ since $\A_1=0$. Let $\tilde{\alpha}\coloneqq(a_1-1,\dots,a_n-1)$, and $m=(a_{n+1}-1)\in[1,k]$. Note that $\abs{\tilde{\alpha}}_i=\A_{i+1}$ for all $i\in[1,n]\setminus\{m\}$, and $\abs{\tilde{\alpha}}_m=\A_{m+1}-1$.\\
For each $j>m$, we get
\begin{multline*}
u_{\tilde{\alpha}}(j)=\sum_{i=j}^n\abs{\tilde{\alpha}}_i-(n-j+1)=\sum_{i=j}^n\A_{i+1}-(n-j+1)\\
=\sum_{i=j+1}^{n+1}\A_i-((n+1)-(j+1)+1)=\ua(j+1)\ge1,
\end{multline*}
whereas for $j'\in[1,m]$
\[
u_{\tilde{\alpha}}(j')=\sum_{i=j'}^n\abs{\tilde{\alpha}}_i-(n-j'+1)
=\sum_{i=j'}^n\A_{i+1}-1-(n-j'+1)=\ua(j'+1)-1\ge0.
\]
Hence, for a fixed $m$, $\tilde{\alpha}$ is being counted by $\unk(m)$, and since $m$ can be any number in $[1,k]$ we have
\begin{equation}
	\label{ucount1}
	\upsilon_{n+1,k}^0(1)=\sum_{m=1}^k\unk(m).
\end{equation}
Now, consider $u_{\tilde{\alpha}}(m)$: if $u_{\tilde{\alpha}}(m)=0$, then $\tilde{\alpha}$ is being counted by $\unk^0(m)$; on the other hand, if $u_{\tilde{\alpha}}(m)\ge1$ (and thus necessarily $m\ge2$), then $\tilde{\alpha}$ gets counted by $\unk(m-1)$. Thus if $m\ge2$
\begin{equation*}
	\unk(m)=\unk^0(m)+\unk(m-1),
\end{equation*}
and by iterating that result, we get
\begin{equation}
	\label{ucount2}
	\unk(m)=\sum_{i=1}^m\unk^0(i)
\end{equation}
for all $m\in[1,k]$.

Suppose $m\ge2$, and that $\tilde{\alpha}$ as above is such that $u_{\tilde{\alpha}}(m)=0$, and thus is counted by $\unk^0(m)$. By Lemmas \ref{translatelemma} and \ref{restrlemma}, given $J=\{i\in[n]\mid a_i\ge m\}$, $\tau_{m-1}(\tilde{\alpha}_{|_J})$ is a $k$-Naples parking function, and it is complete since $u_{\tau_{m-1}(\tilde{\alpha}_{|_J})}(j)=u_{\tilde{\alpha}}(j+m-1)\ge1$ for all $j\in[2,n-m+1]$; therefore, $\tau_{m-1}(\tilde{\alpha}_{|_J})$ is counted by $\upsilon_{n-m+1,k}^0(1)$.\\
Furthermore, given $J^c=[1,n]\setminus J$, $\tilde{\alpha}_{|_{J^c}}$ is such that $u_{\tilde{\alpha}_{|_{J^c}}}(j)=u_{\tilde{\alpha}}(j)\ge0$ for all $j\in[1,m-1]$.\\
Denote $\tilde{\alpha}_{|_{J^c}}=(\hat{a}_1,\dots,\hat{a}_{m-1})$, and let $\beta=(m-\hat{a}_1,\dots,m-\hat{a}_{m-1})$. For all $i\in[1,m-1]$, $\abs{\beta}_i=\abs{\tilde{\alpha}_{|_{J^c}}}_{m-i}$, which implies
\begin{multline*}
u_{\beta}(j)=\sum_{i=j}^{m-1}\abs{\beta}_i-((m-1)-j+1)
=\sum_{i=j}^{m-1}\abs{\tilde{\alpha}_{|_{J^c}}}_{m-i}-((m-1)-j+1)\\
=\sum_{i=1}^{m-j}\abs{\tilde{\alpha}_{|_{J^c}}}_i -((m-1)-j+1)
=-u_{\tilde{\alpha}_{|_{J^c}}}(m-j+1)\le0
\end{multline*}
for all $j\in[2,m-1]$. Thus, $\beta$ is a parking function of length $m-1$; specifically, $\tilde{\alpha}_{|_{J^c}}$ is such that $u_{\tilde{\alpha}_{|_{J^c}}}(j)\ge0$ for all $j\in[2,m-1]$ if and only if $\beta\in PF_{m-1}$, moreover the transformation through which $\beta$ was obtained is clearly an involution. Thus, we have the following recursive relation
\begin{equation}
	\label{ucount3}
	\unk^0(m)=\binom{n}{m-1} m^{m-2} \upsilon_{n-m+1,k}^0(1).
\end{equation}

To conclude, we give the following result:

\begin{theorem}
	\label{completecount}
	Let $n\ge2$ and $k\ge1$. The number of complete $k$-Naples parking functions of length $n$ is 
	\begin{equation}
		\unk^0(1)=
		\begin{cases}
			(n-1)^{n-1}		& \text{if $n\le k+1$}\\
			\sum_{i=1}^k \binom{n-1}{i-1}i^{i-2}(k-i+1)\upsilon_{n-i,k}^0(1)	&\text{if $n\ge k+2$}
		\end{cases}
	\end{equation}
\end{theorem}

\begin{proof}
	We already observed that if $k\ge n+1$, meaning $\unk^0(1)$ counts all parking preferences of length $n$, then $\unk^0(1)=(n-1)^{n-1}$. If $n\ge k+2$, then \eqref{ucount1} and \eqref{ucount2} imply
	\[
	\unk^0(1)=\sum_{m=1}^k\sum_{i=1}^m\upsilon_{n-1,k}^0(i)
	=\sum_{i=1}^k(k-i+1)\upsilon_{n-1,k}^0(i),
	\]
	and we conclude by \eqref{ucount3}.
\end{proof}

\begin{remark}
	Computing $\unk^0(1)$ for $n\in[2,9]$ gives:
	\[
	\begin{array}{ccccccccc}
		\toprule
			&k=1	&k=2	&k=3	&k=4	&k=5	&k=6	&k=7	&k=8\\
		\midrule
		n=2	&\mathbf{1}		&		&		&		&		&		&		&\\
		n=3	&1		&\mathbf{4}		&		&		&		&		&		&\\
		n=4	&1		&11		&\mathbf{27}		&		&		&		&		&\\
		n=5	&1		&38		&131	&\mathbf{256}	&	&	&	&\\
		n=6	&1		&131	&783	&1829	&\mathbf{3125}	&	&	&\\
		n=7	&1		&490	&5136	&15634	&29849	&\mathbf{46656}	&	&\\
		n=8	&1		&1897	&34623	&148321	&332869	&561399	&\mathbf{823543}	&\\
		n=9	&1		&7714	&251817	&1505148 &4102015 &7735566 &11994247 &\mathbf{16777216}\\
		\bottomrule
	\end{array},
	\]
	where the numbers in bold are the number of complete parking sequences of length $n$, $(n-1)^{n-1}$.
	We remark that the diagonal $\upsilon_{n,n-2}^0(1)$, for $n\ge3$, is a known sequence (\href{https://oeis.org/A101334}{A101334} on OEIS, with an offset), with formula $\upsilon_{n,n-2}^0(1)=(n-1)^{n-1}-n^{n-2}$. We prove this formula in the following lemma. As for the rest, there do not seem to be any known sequences (aside from the obvious fact that $\unk^0(1)=(n-1)^{n-1}$ for $k\ge n-1$).
\end{remark}

\begin{lemma}
	Let $n\ge3$, then
	\begin{equation}
		\upsilon_{n,n-2}^0(1)=(n-1)^{n-1}-n^{n-2}.
	\end{equation}
\end{lemma}

\begin{proof}
	Let $n\ge3$. By Theorem \ref{completecount}, and by noting that $n-2\ge(n-i)-1$ for all $i\ge1$, we get
	\begin{multline*}
		\upsilon_{n,n-2}^0(1)=\sum_{i=1}^{n-2}\binom{n-1}{i-1}i^{i-2}(n-i-1)\upsilon_{n-i,n-2}^0(1)\\
		=\sum_{i=1}^{n-2}\binom{n-1}{i-1}i^{i-2}(n-i-1)(n- i -1)^{n-i-1}\\
		=\sum_{j=0}^{n-3}\binom{n-1}{j} (j+1)^{j-1} (n-j-2)^{n- j -1}.
	\end{multline*}
For ease of reference, we rewrite Abel's binomial formula \eqref{Abel}:
\begin{multline*}
(z+w+m)^m=\sum_{j=0}^m\binom{m}{j} w (w+m-j)^{m-j-1} (z+j)^j\\
=\sum_{j'=0}^m\binom{m}{j'} w (w+j')^{j'-1} (z+m-j')^{m-j'}.
\end{multline*}
Setting $m=n-1$, $z=-1$ and $w=1$ gives
\[
\sum_{j=0}^{n-1}\binom{n-1}{j}(j+1)^{j-1}(n-j-2)^{n-j-1}=(n-1)^{n-1}.
\]
Thus,
\begin{multline*}
	\upsilon_{n,n-2}^0(1)=\sum_{j=0}^{n-3}\binom{n-1}{j} (j+1)^{j-1} (n-j-2)^{n- j -1}\\
	=(n-1)^{n-1}-\binom{n-1}{n-2} ((n-2)+1)^{(n-2)-1} (n- (n -2)-2)^{n- (n-2) -1}\\
	-\binom{n-1}{n-1} ((n-1)+1)^{(n-1)-1} (n-(n-1)-2)^{n- (n-1) -1}\\
	=(n-1)^{n-1}-n^{n-2}.
\end{multline*}
\end{proof}

\chapter{Parking Strategies}
\label{chapterstrategies}
\section{Definition and properties of parking strategies}

In the previous chapters, we have only considered problems where cars all have the same $k$-Naples parking rule; in this chapter we will generalize this concept, allowing different cars to have different parking rules.

\begin{definition}
	Let $\alpha=(a_1,\dots,a_n)\in\pp$. To the parking preference $\alpha$, we associate a vector
	\begin{equation}
		\rho\coloneqq(r_1,\dots,r_n),
	\end{equation}
	where each $r_i$ is an integer chosen in $[0,n]$. We now consider the problem where for each $i\in[n]$, car $c_i$ follows the $r_i$-Naples parking rule, where by $0$-Naples we mean the standard parking rule. We define the \emph{parking relation} of length $n$ as the set $\pr\subseteq [1,n]^n\times [0,n]^n$ such that $(\alpha,\rho)\in\pr$ if and only if in the associated problem with $n$ cars, where each car $c_i$ has preference $a_i$ and follows the $r_i$-Naples parking rule, all cars are able to park. If $(\alpha,\rho)\in\pr$, we say that $\rho$ is a \emph{parking strategy} for $\alpha$. We denote by $\psi(\alpha,\rho)$ the generalization of the outcome map defined in \ref{outcomemap}.
\end{definition}

\begin{example}
	Let $\alpha=(4,3,3,4,1)$. Setting $\rho=(k,k,k,k,k)$ is the same as considering the $k$-Naples parking problem. In this case, $\alpha\notin PF_5$, but $\alpha\in PF_{5,k}$ for all $k\ge1$, thus $(\alpha,00000)\notin PR_5$, but $(\alpha,kkkkk)\in PR_5$. 
	Choosing vectors with non constant values can yield different results. For instance, $\psi(43341,11011)=(4,3,5,\infty,1)$, thus $(\alpha,11011)\notin PR_5$, while  $\psi(43341,00100)=(4,3,2,5,1)$, so $(\alpha,00100)\in PR_5$. Moreover, $(43341,00020)=(4,3,5,2,1)$, so $(\alpha,00020)\in PR_5$.
\end{example}

In the previous example, it is easy to check that, given $\alpha=(4,3,3,4,1)$, if $\rho$ is such that the third car is allowed to drive backwards at least one spot, or such that the fourth car can drive backwards at least two spots (i.e. $r_3\ge1$ and/or $r_4\ge2$), then $(\alpha,\rho)\in\pp$. The next Lemma generalizes this concept for an arbitrary parking preference $\alpha\in\pp$.

\begin{lemma}
	\label{orderlemma}
	Let $\alpha\in\pp$, and $\rho,\rho'\in[0,n]^n$  be such that $r_i\le r_i'$, for all $i\in[n]$. If $(\alpha,\rho)\in\pr$, then $(\alpha,\rho')\in\pr$.
\end{lemma}

\begin{proof}
	Let $(\alpha,\rho)\in\pr$ and $\rho'$ be such that $r_i\le r_i'$ for all $i\in[n]$. Suppose $(\alpha,\rho')\notin\pr$; for the sake of simplicity, let $\psi(c_i)$ be the outcome of car $c_i$ in $(\alpha,\rho)$, and $\psi'(c_i')$ be the outcome of car $c_i'$ in $(\alpha,\rho')$. Clearly $\psi(c_1)=a_1=\psi'(c_1')$. Suppose by induction that, for all $i\in[1,h-1]$, if $\psi'(c_i')\ge\psi(c_i)$ then $\psi'(c_i')\in\{\psi(c_1),\dots,\psi(c_i)\}$, and consider car $c_h'$. 
	The induction is trivial if $\psi'(c_h')=\psi(c_h)$ (since of course, $\psi'(c_h')\in\{\psi(c_1),\dots,\psi(c_h)\}$). Therefore, suppose that $\psi'(c_h')>\psi(c_h)$. If $\psi(c_h)<\psi'(c_h')\le a_h$, then car $c_h$ must have reached spot $\psi(c_h)$ by driving backwards from $a_h$, in particular it must have checked spot $\psi'(c_h')$, finding it occupied, which means that it was filled by a previous car, hence $\psi'(c_h')\in\{\psi(c_1),\dots,\psi(c_{h-1})\}$.\\
	If instead $a_h<\psi'(c_h')$, then car $c_h'$ must have checked spot $\psi(c_h)$, finding it full: this is trivial if $a_h\le\psi(c_h)<\psi'(c_h')$, since car $c_h'$ must find all spots in $[a_h,\psi'(c_h')-1]$ occupied before reaching $\psi'(c_h')$; if instead $\psi(c_h)<a_h$, then $\psi(c_h)\in[a_h-r_h,a_h-1]$, and since $r_h'\ge r_h$ car $c_h'$ has found all spots in $[a_h-r_h',a_h]\supseteq[a_h-r_h,a_h-1]$ already full. Thus, 
	\[
	\psi(c_h)\in\{\psi'(c_1)',\dots,\psi'(c_{h-1})\}.
	\]
	Let $j<h$ be such that $\psi'(c_j')=\psi(c_h)$: note that since car $c_h$ finds spot $\psi(c_h)$ not occupied, then $\psi'(c_j')\notin\{\psi(c_1),\dots,\psi(c_{h-1})\}$, and specifically, $\psi'(c_j')\notin\{\psi(c_1),\dots,\psi(c_{j})\}$. Hence by induction $\psi'(c_j')<\psi(c_j)$. Now, clearly if $\psi(c_j)=\psi'(c_h')$ the induction holds. Thus there are two cases: $\psi'(c_h')<\psi(c_j)$ or $\psi'(c_h')>\psi(c_j)$.\\
	Suppose $\psi'(c_h')<\psi(c_j)$: if $\psi'(c_h')<a_j$, then since $j<h$ and $\psi'(c_j')=\psi(c_h)<\psi'(c_h')<a_j$, then spot $\psi'(c_h')$ should be already occupied when car $c_j'$ parks, and thus also when car $c_h'$ tries to park, and we have a contradiction. Thus, $a_j\le\psi'(c_h')<\psi(c_j)$; consequently, car $c_j$ certainly finds spot $\psi'(c_h')$ occupied, so $\psi'(c_h')\in\{\psi(c_1),\dots,\psi(c_{j-1})\}$ and we conclude.\\
	Suppose instead that $\psi'(c_h')>\psi(c_j)$: since by induction $\psi'(c_h')>\psi(c_j)>\psi'(c_j')=\psi(c_h)$, by the same argument used on $\psi(c_h)$, spot $\psi(c_j)$ is found occupied by car $c_h'$, and we again find that there exists $k<h$ such that 
	\[
	\psi(c_h)=\psi'(c_j')<\psi(c_j)=\psi'(c_k')<\psi(c_k).
	\]
	In the same way as for car $c_j$, if $\psi'(c_h')\le\psi(c_k)$ then the induction holds, and if instead $\psi'(c_h')>\psi(c_j)$ then we can again reiterate this argument. Since there are at most $h-1$ cars with index lesser than $h$, at some point we must find
	\[
	\psi(c_h)=\psi'(c_{j_1}')<\psi(c_{j_1})=\psi'(c_{j_2}')<\dots<\psi(c_{j_m})
	\]
	such that $\psi'(c_h')\le\psi(c_{j_m})$, which therefore proves the induction.\\
	Hence, for all $h\in[n]$ either $\psi'(c_h')<\psi(h)$, or $\psi'(c_h')\in\{\psi(c_1),\dots,\psi(c_h)\}$. Furthermore, $(\alpha,\rho)\in\pr$, so $\psi(c_i)\in[n]$ for all $i\in[n]$: thus $\psi'(c_h')\in[n]$ for all $h\in[n]$, and so $(\alpha,\rho')\in\pr$.\\
\end{proof}

\begin{definition}
	Let $X$ be a set, a \emph{partial order} $\lesssim$ on $X$ is a binary relation such that: 
	\begin{itemize}
		\item for all $x\in X$, $x\lesssim x$ (\emph{reflexivity});
		\item for all $x,y\in X$, if $x\lesssim y$ and $y\lesssim x$, then $x=y$ (\emph{antisymmetry});
		\item for all $x,y,z\in X$, if $x\lesssim y$ and $y\lesssim z$, then $x\lesssim z$ (\emph{transitivity}).
	\end{itemize}
	Furthermore, if for all $x,y\in X$, either $x\lesssim y$ or $y\lesssim x$ holds (\emph{completeness}), then $\lesssim$ is called a \emph{total order}. If $\lesssim$ is a partial order on $X$, then ($X$,$\lesssim$) is called a \emph{partially ordered set}, or \emph{poset} for brevity.\\
	Let $(X,\lesssim)$ be a poset, a subset $A\subseteq X$ is called an \emph{up-set} if for all $x\in A$, if $y\in A$ is such that $x\lesssim y$, then $y\in A$.
\end{definition}

Consider the set $[0,n]^n$, and the order relation $\lesssim$ defined by declaring $\rho\lesssim\rho'$ if and only if $r_i\le r_i'$ for all $i\in[n]$. With this definition, $([0,n]^n,\lesssim)$ is a poset. Given $\alpha\in\pp$, Lemma \ref{orderlemma} implies that if $(\alpha,\rho)\in\pr$, then $(\alpha,\rho')\in\pr$ for all $\rho'\gtrsim\rho$. Thus, the set of all parking strategies for $\alpha$ is an up-set.

Consider the following subset of $[0,n]^n$:
\[
\R_n\coloneqq\{0\}\times[0,1]\times\dots\times[0,n-1].
\]
Clearly, $(\R_n,\lesssim)$ is also a poset, with maximum $(0,1,\dots,n-1)$, and minimum $(0,\dots,0)$. In particular, it is a lattice (as is $([0,n]^n,\lesssim)$), although this irrelevant for our purposes. Given $(\alpha,\rho)\in\pr$, define $\tilde{\rho}=(\tilde{r}_1,\dots,\tilde{r}_n)$ by:
\[
\tilde{r}_i=
\begin{cases*}
	0		&\text{if $\psi(c_i)\ge a_i$}\\
	a_i-\psi(c_i)	&\text{if  $\psi(c_i)<a_i$}
\end{cases*}
\quad \text{for all $i\in[n]$}.
\]
By construction, $\psi(\alpha,\rho)=\psi(\alpha,\tilde{\rho})$. Moreover, for each $i$, car $c_i$ has driven backwards at most $i-1$ spots, since only $i-1$ other cars already had their turn. Thus, $\tilde{\rho}\in\R_n$, $(\alpha,\tilde{\rho})\in\pr$, and clearly $\tilde{\rho}\lesssim\rho$.

As we have observed, $\alpha$ induces an up-set in $[0,n]^n$ containing all of its parking strategies. By what we have just noticed, every minimal element of such an up-set is in $\R_n$, so we can focus on that set, by remarking that $\alpha$ induces an up-set $\mathcal{A}\subseteq\R_n$, with the same minimal elements as the up-set induced in $[0,n]^n$. 

Clearly, $\mathcal{A}=\R_n$ if and only if $\alpha\in\pf$. The following lemma will show that in fact, all \emph{principal filters} of the lattice $\R_n$, i.e. all up-sets with a single minimal element, are represented by some parking preferences $\alpha\in\pp$.

\begin{theorem}
	Let $\rho\in\R_n$. Then there exists a parking preference $\alpha\in\pp$ such that: $(\alpha,\rho')\in\pr$ if and only if $\rho\le\rho'$. Specifically, one such parking preference is given by defining $\alpha$ as $a_i=n+1-i+r_i$, for all $i\in[n]$.
\end{theorem}

\begin{proof}
	Let $\rho\in\R_n$, and let $\alpha$ be defined as $a_i=n+1+r_i$ for all $i\in[n]$. We will prove by induction that $\psi(c_i)=n+1-i$ for all $i\in[n]$. \\
	Note that $r_1=0$, so $a_1=n$, and consequently $\psi(c_1)=n$. Let $j\in[n]$, and suppose that the induction holds for all $i\in[1,j-1]$. Thus, the previous cars have occupied spots $[n+2-j,n]$. Consider car $c_j$: if $r_j=0$, then $a_j=n+1-j$, and we simply get $\psi(c_j)=a_j=n+1-j$. On the other hand, if $r_j\ne0$, then $a_j\in[n+2-j,n]$: since spot $a_j$ is occupied, car $c_j$ will drive backwards for up to $r_j$ spots until it reaches the first available spot, which is $a_j-r_j=n+1-j$. Crucially, to be able to park, car $c_j$ \emph{needs} to drive backwards (for exactly $r_i$ spots), as all the following spots are already occupied, and since this was the case for all previous cars as well, any vector $\rho'\in\R_n$ such that $(\alpha,\rho')\in\pr$ must satisfy $r_j'\ge r_j$. Thus, we conclude.
\end{proof}

\begin{remark}
	Note that while we have just proved that any filter of $\R_n$ is represented by a parking preference, the same does not happen for a generic up-set $\mathcal{A}\in\R_n$. For example, let $n=3$, and consider $\mathcal{A}=\{011,002,012\}$, which is the up-set in $\R_n$ generated by the elements $(0,1,1)$ and $(0,0,2)$. \\
	Suppose that there existed $\alpha\in PP_3$ such that $(\alpha,\rho)\in PR_3$ if and only if $\rho\in\mathcal{A}$.
	Suppose $a_1\ne a_2$, then certainly $\psi(c_1)=a_1$ and $\psi(c_2)=a_2$ regardless of the choice of $\rho$, thus certainly if $(\alpha,011)\in PR_3$, then also $(\alpha,001)\in PR_3$, which is not the case since $(0,0,1)\notin\mathcal{A}$.
	Thus, $a_1=a_2$, however:
	\begin{itemize}
		\item if $a_1=a_2=1$, then $\psi(c_1)=1$ and $\psi(c_2)=2$ independently of the choice of $\rho$, thus we again get $(\alpha,001)\in PR_3$ and a contradiction;
		\item if $a_1=a_2=2$, then considering $(\alpha,011)$, $\psi(c_1)=2$ and $\psi(c_2)=1$, thus no matter what the values for $a_3$ and $r_3$ are, we obtain $\psi(c_3)=3$ and $(\alpha,010)\in PR_3$, which is not the case;
		\item if $a_1=a_2=3$, then considering $(\alpha,002)$, car $c_2$ does not park, which is absurd since $(0,0,2)\in\mathcal{A}$.
	\end{itemize}
	Hence, not all up-sets in $\R_n$ are induced by some parking preference.
\end{remark}

\section{Strategies minimizing the number of steps}

Considering the result found in Lemma \ref{orderlemma}, given a parking preference $\alpha\in\pp$, it can be worthwhile to study cases where $\rho$ is minimized in some sense.
Thus, we are going to analyse the following problem: given $\alpha\in\pp$, can we find $\rho$ such that $(\alpha,\rho)\in\pr$, and $\rho$ minimizes $\abs{\rho}\coloneqq\sum_{i=1}^n r_i$? Clearly, it is optimal to choose $\rho$ such that if car $c_i$ parks $j$ spaces backwards from $a_i$, then $r_i=j$, so from another point of view our problem becomes: "find a parking strategy $\rho$ for $\alpha$, such that $\rho$ minimizes the number of backward steps taken by the cars".

Observe that, in a way, we are considering the reverse of the usual parking problem: instead of evaluating whether a given vector $\rho$ is a parking strategy for $\alpha\in\pp$, we are \emph{choosing} parking strategies that satisfy certain properties.

\begin{remark}
Note that $\abs{\rho}=\sum_{i=1}^n r_i$ is none other than the rank of $\rho$ in the lattice $\R_n$, so, given $\alpha\in\pp$, we are trying to find the minimal elements in $\mathcal{A}$ with the smallest rank (note, however, that these are generally not all of the minimal elements of $\mathcal{A}$).
\end{remark}

\begin{remark}
	The problem considered above is \emph{not} the same as the problem: "given $\alpha\in\pp$, assuming all cars follow the $n$-Naples parking rule, how many backward steps do the cars take?". For example, consider $\alpha=44323$: $\psi_n(\alpha)=43215$, meaning that cars $c_2$, $c_3$ and $c_4$ each take a backward step to park, and indeed $(44323,01110)\in PR_5$. However, the optimal way to choose $\rho$ is actually $\rho=00002$, where $(44323,00002)\in PR_5$ and cars only need to take a total of two steps backwards.\\
\end{remark}

A direct consequence of Proposition \ref{leastprop} is the following lower bound:

\begin{lemma}
	\label{greaterlemma}
	Let $\alpha\in\pp$, and $\rho\in[0,n]^n$. If $(\alpha,\rho)\in\pr$, then 
	\begin{equation}
		\abs{\rho} \ge \sum_{j\in\Ua} \ua(j).
	\end{equation}
\end{lemma}

\begin{proof}
	Let $j\in\Ua$; by Proposition \ref{leastprop} at least $\ua(j)$ cars with preference in $[j,n]$ park in spots in $[1,j-1]$. In particular, at least $\ua(j)$ cars have to take the backward step from spot $j$ to spot $(j-1)$. Thus, the total number of backward steps taken is at least $\sum_{j\in\Ua}\ua(j)$.
\end{proof}

We will now prove that this bound can be reached, and in fact we can construct a parking strategy such that $\abs{\rho}=\sum_{j\in\Ua}\ua(j)$.

\begin{theorem}
	\label{constructrho}
	Let $\alpha\in\pp$. For each maximal interval $[p,q]\subseteq\Ua$, and $j\in[p,q]$, define set $T_j$ in the following way:
	\begin{itemize}
		\item if $j=q$, then $T_q$ is the set of the last $\ua(q)$ indices $i\in[n]$ such that $a_i=q$;
		\item if $j\in[p,q-1]$ , then $T_j$ is the set of the last $\ua(j)$ indices $i\in[n]$ such that either $a_i=j$, or $i\in T_{j+1}$.
	\end{itemize}
	Define
	\begin{equation}
		\rho=\sum_{j\in\Ua}\sum_{i\in T_j} e_i,
	\end{equation}
	where $e_i\in[0,n]^n$ is the vector with all zeroes except for a $1$ in the $i$-th place.\\
	Then $(\alpha,\rho)\in\pr$, and $\abs{\rho}=\sum_{j\in\Ua}\ua(j)$.
\end{theorem}

\begin{proof}
	Let $[p,q]\subseteq\Ua$ be a maximal interval, then $(q+1)\notin\Ua$, so $\ua(q+1)\le0$ and consequently $\A_q=\ua(q)-\ua(q+1)+1\ge\ua(q)+1$, so set $T_q$ is well defined and contains $\ua(q)$ elements. By induction, given $j\in[p,q]$, suppose that for all $i\in[j+1,q]$ set $T_i$ is well defined and contains $\ua(i)$ elements; then the number of indices such that either $a_i=j$ or $i\in T_{j+1}$ is $\A_j+\ua(j+1)=\ua(j)+1$, so set $T_j$ is well defined and contains $\ua(j)$ elements.
	Thus, for all $j\in\Ua$ set $T_j$ is well defined, and $\abs{T_j}=\ua(j)$, therefore since each $e_i$ adds $1$ to $\abs{\rho}$, $\abs{\rho}=\sum_{j\in\Ua}\ua(j)$.\\
	Let $[p,q]\subseteq\Ua$, and consider the cars with preference $q$: we remarked that $\A_q\ge\ua(q)+1$. After the first of such cars attempts to park, spot $q$ is certainly occupied, and there are at least $\ua(q)$ such cars remaining. Thus, for all $i\in T_q$, car $c_i$ has preference $q$ and is provided by a backward step by $e_i$ ($q\in\Ua$, $i\in T_q$), so it will reach spot $q-1$. Suppose by induction that, given $j\in[p,q-1]$, exactly $\ua(j+1)$ cars reach spot $j$ by using backward steps provided by sets $T_{j+1},\dots,T_q$. Then, there is a total of $\A_j+\ua(j+1)$ cars that either have preference $j$, or reach spot $j$ using backward steps provided by sets $T_{j+1},\dots,T_q$: the first among these cars fills spot $j$ (if required), whereas all the remaining $\A_j+\ua(j+1)-1=\ua(j)$ cars are granted a backward step by virtue of set $T_j$, and thus reach spot $j-1$.\\
	Hence, the induction holds for all $j\in[p,q]$. In particular, since $\ua(p)=1$, a car reaches spot $p-1$, and fills it if it is empty. Observe that the proof provided for Theorem \ref{theoprec} can easily be extended to the case with varied parking rules, so to prove $(\alpha,\rho)\in\pr$ it is sufficient to check that for every maximal interval $[p,q]\subseteq\Ua$, spot $(p-1)$ gets occupied. Thus, we conclude.
\end{proof}

\begin{example}
	\label{examplemin1}
	To better illustrate the construction of $\rho$ according to Theorem \ref{constructrho}, we provide an example. Let $\alpha=(4,9,5,9,5,8,7,9,4,6)\in PP_{10}$: $\Ua=[2,9]$.
	\[
	\begin{array}{ccccc}
		\toprule
		&\ua(j)	&\{i\in[n]\mid a_i=j\}	&\{i\in[n]\mid a_i=j\}\cup T_{j+1}	&T_j\\
		\midrule
		j=9	&1		&\{2,4,8\}				&\{2,4,8\}							&\{8\}\\
		j=8	&1		&\{6\}					&\{6,8\}							&\{8\}\\
		j=7	&1		&\{7\}					&\{7,8\}							&\{8\}\\
		j=6	&1		&\{10\}					&\{8,10\}							&\{10\}\\
		j=5	&2		&\{3,5\}				&\{3,5,10\}							&\{5,10\}\\
		j=4	&3		&\{1,9\}				&\{1,5,9,10\}						&\{5,9,10\}\\
		j=3	&2		&\varnothing			&\{5,9,10\}							&\{9,10\}\\
		j=2	&1		&\varnothing			&\{9,10\}							&\{10\}\\
		\bottomrule
	\end{array}
	\]
	Thus, $\rho=2e_5+3e_8+2e_9+5e_{10}=(0,0,0,0,2,0,0,3,2,5)$ and we can easily check that
	\[
	\psi(4959587946,0000200325)=(4,9,5,10,3,8,7,6,2,1).
	\]
	Hence, $(4959587946,0000200325)\in PR_{10}$.
\end{example}

\begin{remark}
	A natural follow-up question to Theorem \ref{constructrho} is whether there is uniqueness over the construction of $\rho$ such that it minimizes $\abs{\rho}$. The answer is, generally, negative. However, as we will now see, any parking strategy $\rho$ such that it minimizes $\abs{\rho}$ still follows a similar construction to that found in Theorem \ref{constructrho}.\\
	Let $(\alpha,\rho)\in\pr$ be such that $\abs{\rho}=\sum_{j\in\Ua}\ua(j)$, and let $j\in\Ua$. By Proposition \ref{leastprop} at least $\ua(j)$ cars with preference at least $j$ must park in spots in $[1,j-1]$: each of those cars has to make the backward step from $j$ to $(j-1)$, so for $\rho$ to minimize $\abs{\rho}$ the number of cars with preference in $[j,n]$ that park in $[1,j-1]$ must be \emph{exactly} $\ua(j)$. In particular, the number of cars that make the backward step from $j$ to $(j-1)$ must be $\ua(j)$.\\
	Let $[p,q]\subseteq\Ua$ be a maximal interval, note that $\ua(p-1)=0$, so $\sum_{i=p-1}^n\A_i=n-p+2$. If there existed a car $c_i$ with preference $a_i<p-1$ such that it parks in $[p-1,n]$, then one of the cars with preference in $[p-1,n]$ needs to park in $[1,p-2]$, so it must make the backward step from $(p-1)$ to $(p-2)$, which is distinct from the ones counted in $\sum_{j\in\Ua}\ua(j)$, so $\rho$ would not minimize $\abs{\rho}$. Thus, all spots in $[p-1,q]$ are filled by cars in $[p-1,n]$.\\
	Consider spot $(p-1)$: since $\A_{p-1}=0$, it is filled by a car that drives backwards from $p$. By what we just remarked, the number of cars that reach spot $p$ from $(p+1)$ is $\ua(p+1)$, so the total number of cars that ever reach spot $p$ is $\ua(p+1)+\A_p$. Among those, one must fill spot $p$, whereas \emph{all} the remaining $\ua(p+1)+\A_p-1=\ua(p)$ cars will have to take the backward step from $p$ to $(p-1)$. Thus, a backward step (and a corresponding added value to $\rho$) is associated to each of the last $\ua(p)$ cars to reach spot $p$. Observe that as a consequence, cars with preference $p$ park in $[p-1,p]$.\\
	In an analogous way, for all $j\in[p,q-1]$ exactly $\ua(j)+1$ cars reach spot $j$, either because they have preference $j$ or because they come backwards from $(j+1)$, and a backward step is associated to all but the first of those cars. For $j=q$ the situation is similar, however no cars from $(q+1)$ will reach $q$, and since $\ua(q+1)\le0$ we have
	\[
	\A_q=\ua(q)-\ua(q+1)+1\ge \ua(q)+1,
	\]
	and this is not necessarily an equality. Thus, the $\ua(q)$ cars that take the backward step from $q$ to $(q-1)$ can be chosen between any of the $\A_q$ cars that have preference $q$ except from the first, and this choice is not necessarily unique. Hence, $\rho$ is obtained in the same way as in Theorem \ref{constructrho}, with the only difference being that for each maximal interval $[p,q]\in\Ua$, the indices in set $T_q$ can be freely chosen between any of the cars with preference $q$ except for the first, instead of necessarily choosing the last $\ua(j)$ indices. Moreover, the proof provided in \ref{constructrho} can clearly also be applied to this case to prove that any vector constructed in such a way is a parking strategy.
	Thus, we have proved the following:
\end{remark}

\begin{theorem}
	\label{choiceminimize}
	Let $\alpha\in\pp$, and $\rho\in[0,n]^n$. Then $\rho$ is a parking strategy such that  $\abs{\rho}=\sum_{j\in\Ua}\ua(j)$ if and only if for all maximal intervals $[p,q]\subseteq\Ua$ and $j\in[p,q]$, there exist some sets $T_j$ such that
	\begin{equation}
		\rho=\sum_{j\in\Ua}\sum_{i\in T_j} e_i,
	\end{equation}
	and defined as follows:
	\begin{itemize}
		\item if $j=q$, then $T_q$ is  a set of any $\ua(q)$ distinct indices chosen among all cars with preference $q$, except for the first such car;
		\item if $j\in[p,q-1]$ , then $T_j$ is the set of the last $\ua(j)$ indices $i\in[n]$ such that either $a_i=j$, or $i\in T_{j+1}$.
	\end{itemize}
\end{theorem}

\begin{corollary}
	\label{minimalrhocount}
	Let $\alpha\in\pp$, and $\Ua=\bigcup_{i=1}^m [p_i,q_i]$ be the decomposition of $\Ua$ in maximal intervals. Then the number of distinct parking strategies such that $\abs{\rho}=\sum_{j\in\Ua}\ua(j)$ is 
	\begin{equation}
		\prod_{i=1}^m \binom{\A_{q_i}-1}{\ua(q_i)}.
	\end{equation}
\end{corollary}

\begin{proof}
	Two different ways to choose sets $T_{q_i}$, for $i=1,\dots,m$, generate different minimizing vectors, thus the number of ways to choose $\rho$ such that $(\alpha,\rho)\in\pr$ and $\abs{\rho}=\sum_{j\in\Ua}\ua(j)$ is equal to the number of ways to choose such sets, each of which consists in choosing $\ua(q_i)$ cars among $\A_{q_i}-1$.
\end{proof}

\begin{example}
	\label{examplemin2}
	Going back to the previous example, let $\alpha=(4,9,5,9,5,8,7,9,4,6)$, where $\Ua=[2,9]$. Here $\A_9=3$ is strictly greater than $\ua(9)+1=2$, so there is an additional way to construct $\rho$ to minimize $\abs{\rho}$, and that is to choose $T_9={4}$ instead of $T_9={8}$. In that case, we get
	\[
	\begin{array}{ccccc}
		\toprule
		&\ua(j)	&\{i\in[n]\mid a_i=j\}	&\{i\in[n]\mid a_i=j\}\cup T_{j+1}	&T_j\\
		\midrule
		j=9	&1		&\{2,4,8\}				&\{2,4,8\}							&\{4\}\\
		j=8	&1		&\{6\}					&\{4,6\}							&\{6\}\\
		j=7	&1		&\{7\}					&\{6,7\}							&\{7\}\\
		j=6	&1		&\{10\}					&\{7,10\}							&\{10\}\\
		j=5	&2		&\{3,5\}				&\{3,5,10\}							&\{5,10\}\\
		j=4	&3		&\{1,9\}				&\{1,5,9,10\}						&\{5,9,10\}\\
		j=3	&2		&\varnothing			&\{5,9,10\}							&\{9,10\}\\
		j=2	&1		&\varnothing			&\{9,10\}							&\{10\}\\
		\bottomrule
	\end{array}
	\]
	and thus $\rho=e_4+2e_5+e_6+e_7+2e_9+5e_{10}=(0,0,0,1,2,1,1,0,2,5)$, and $(4959587946, \\0001211025)\in PR_{10}$.\\
	Since $[2,9]$ is the only maximal interval of $\Ua$, and the only way to choose $T_9$ are $T_9=\{4\}$ and $T_9=\{8\}$, $\rho=0001211025$ and $\rho=0000200325$ are the only parking strategies such that $(4959587946,\rho)\in PR_{10}$ and $\abs{\rho}=\sum_{j\in\Ua}\ua(j)=12$.\\
\end{example}

In the previous section, we remarked that for $n=3$, there exists an up-set in $\R_3$ that is not induced by any parking preference. That result can be generalized for all $n\ge3$.

\begin{lemma}
	Let $n\ge3$, and consider the up-set $\mathcal{A}\subset\R_n$ generated by all elements of rank $\frac{n(n-1)}{2}-1$. Then $\mathcal{A}$ is not induced by any parking preference $\alpha\in\pp$. 
\end{lemma}

\begin{proof}
	Suppose that there exists $\alpha\in\pp$ such that it induces $\mathcal{A}$ defined as above. Then any element $\rho\in\R_n$ of rank $\abs{\rho}=\frac{n(n-1)}{2}-1$, is such that it minimizes $\abs{\rho}$, thus 
	\begin{equation}
		\label{gorb}
		\frac{n(n-1)}{2}-1=\abs{\rho}=\sum_{j\in\Ua}\ua(j).
	\end{equation}
	Suppose $\ua(n)=n-1$, then by definition $\A_n=\ua(n)+1=n$, and thus $\alpha=(n,n,\dots,n)$. As a consequence,
	\[
	\ua(j)=\sum_{i=j}^n\A_i-(n-j+1)=\A_n-(n-j+1)=j-1
	\]
	for all $j\in[2,n]$. Thus,
	\[
	\sum_{j\in\Ua}\ua(j)=\sum_{j=2}^n (j-1)=\sum_{i=1}^{n-1} i =\frac{n(n-1)}{2},
	\]
	which contradicts \eqref{gorb}.\\
	Thus, $\ua(n)\le n-2$. By \eqref{propua4}, $\ua(j)\le j-1$ for all $j\in[2,n-1]$, so given all of these restrictions, the only way to get $\sum_{j\in\Ua}\ua(j)=\frac{n(n-1)}{2}-1$ is if $\ua(j)=j-1$ for all $j\in[2,n-1]$, and $\ua(n)=n-2$. Specifically, $\Ua=[2,n]$ is a single maximal interval, and $\A_n=\ua(n)+1=n-1$. By Corollary \ref{minimalrhocount}, the number of distinct vectors $\rho\in\R_n$ such that $(\alpha,\rho)\in\pr$ and $\abs{\rho}=\frac{n(n-1)}{2}-1$ is $\binom{\A_n-1}{\ua(n)}=\binom{n-2}{n-2}=1$. However, the number of elements of rank $\frac{n(n-1)}{2}-1$ in $\R_n$ is $n-1\ge3-1=2$, so we have a contradiction and we conclude.\\
\end{proof}

As a final observation, consider $\alpha\in\pp$ and a parking strategy $\rho\in\R_n$ for $\alpha$ such that $\abs{\rho}=\sum_{j\in\Ua}\ua(j)$ (that is, $\rho$ minimizes the backward steps taken by the cars). For all $r_i\ne0$, $\psi(c_i)=a_i-r_i$, and all other cars only take forward steps. For the purposes of counting steps, we assume that each car starts from its preference, and disregard how it got there. Thus, if we define $\tilde{\alpha}$ such that $\tilde{a}_i=a_i-r_i$, $\tilde{\alpha}$ is a parking function, and the number of forward steps taken by the cars in $\tilde{\alpha}$ is the same as that for $(\alpha,\rho)$.\\
Now, considering $\tilde{\alpha}$, note that each forward step taken by the cars is such that it increases a value in $\tilde{\alpha}$, until every car reaches the spot it will occupy: therefore, accounting for all forward steps, $\tilde{\alpha}$ will be increased up to a permutation of $\{1,2,\dots,n\}$ (i.e. the outcome $\psi(\tilde{\alpha})$). Hence, the number of forward steps taken is 
\begin{align*}
\sum_{i=1}^n i -\biggl(\sum_{i=1}^n \tilde{a}_i\biggr)&=\frac{n(n+1)}{2}-\biggl(\sum_{i=1}^n a_i -\sum_{j\in\Ua}r_j  \biggr)\\
&=\frac{n(n+1)}{2}-\sum_{j=1}^n j\A_j+\sum_{j\in\Ua}\ua(j)\\
&=\frac{n(n+1)}{2}-\sum_{j=1}^n j\A_j+\sum_{j=1}^n\ua(j)-\sum_{j\notin\Ua}\ua(j)\\
=\frac{n(n+1)}{2} &-\sum_{j=1}^n j\A_j+\sum_{j=1}^n\biggl(\sum_{k=j}^n\A_k-(n-j+1)\biggr)-\sum_{j\notin\Ua}\ua(j)\\
=\frac{n(n+1)}{2} &-\sum_{j=1}^n j\A_j+\sum_{j=1}^n j\A_j -\sum_{j=1}^n(n-j+1)-\sum_{j\notin\Ua}\ua(j)\\
=\frac{n(n+1)}{2} &-\sum_{j'=1}^n j'-\sum_{j\notin\Ua}\ua(j)=-\sum_{j\notin\Ua}\ua(j).
\end{align*}
Note that $\ua(j)\le$ for all $j\notin\Ua$, so the result is indeed non-negative.

Now, define $\beta\in\pp$ by $b_i=n+1-a_i$ for all $i\in[n]$, and define $\rho'=\{r_1'\dots,r_n'\}$ by $r_i'=0$ if $r_i\ne0$, and $r_i'=\psi(c_i)-a_i$ otherwise. Denote by $\psi'$ the outcome map for $(\beta,\rho')$.\\
Observe that the problem for $(\beta,\rho')$ is the reflection of the problem for $(\alpha,\rho)$, in the sense that for all $i\in[n]$, $\psi'(c_i')=n+1-\psi(c_i)$. Furthermore, if car $c_i$ takes $\psi(c_i)-a_i=r_i'$ forward steps to park, car $c_i'$ takes $r_i'$ backward steps; conversely, if $c_i$ takes $r_i$ backward steps, then $c_i'$ takes $r_i$ forward steps. 
We can prove this by induction. Clearly, $\psi'(c_1')=b_1=n+1-a_1=n+1-\psi(c_1)$. Moreover, let $h\ge2$, and assume that $\psi'(c_i')=n+1-\psi(c_i)$ for all $i\in[1,h-1]$.
\begin{itemize}
	\item If $r_h=0$, then $r_h'=\psi(c_h)-a_h$. The first $h-1$ cars in $(\alpha,\rho)$ occupy spots $[a_h,\psi(c_h)-1]$, and leave spot $\psi(c_h)$ free, by induction the first $h-1$ cars for $(\beta,\rho')$ occupy spots $[n+1-(\psi(c_h)-1),n+1-a_h]=[b_h-r_h'+1,b_h]$, but do not occupy spot $n+1-\psi(c_h)=b_h-r_h'$. Thus, car $c_h'$ performs $r_i'$ backward steps, and parks in spot $b_h-r_h'=n+1-\psi(c_h)$.
	\item If $r_h>0$, then $r_h'=0$ and car $c_h$ parks by driving forwards. The first $h-1$ cars in $(\alpha,\rho)$ occupy spots $[a_h-r_h+1,a_h]$, and leave spot $a_h-r_h=\psi(c_h)$ free, by induction the first cars in $(\beta,\rho')$ occupy spots $[n+1-a_h,n+1-(a_h-r_h+1)]=[b_h,b_h+r_h-1]$, but do not occupy spot $b_h+r_h$. Thus, car $c_h'$ performs $r_i$ forward steps, and parks in spot $b_h+r_h=n+1-\psi(c_h)$.
\end{itemize}
Now, note that $\abs{\beta}_i=\A_{n+1-i}$ for all $i\in[n]$. For all $j\in[2,n]$, we get
\begin{multline*}
	u_{\beta}(j)=\sum_{i=j}^n\abs{\beta}_i -(n-j+1)
	=\sum_{i=j}^n\A_{n+1-i}-(n-j+1)\\
	=\sum_{i=1}^{n-j+1}\A_i-(n-j+1)=-\ua(n-j+2).
\end{multline*}
Thus, $u_{\beta}(j)\ge0$ if and only if $n-j+2\notin\Ua$. Note that since $(\beta,\rho')$ is the reflection of $(\alpha,\rho)$, the total number of backward steps taken in $(\beta,\rho')$ is equal to the total number of forward steps taken in $(\alpha,\rho)$, i.e. $-\sum_{j\notin\Ua}\ua(j)$. Furthermore,
\begin{equation*}
-\sum_{j\notin\Ua}\ua(j)=\sum_{j:u_{\beta}(n-j+2)\ge0}u_{\beta}(n-j+2)
=\sum_{j:n-j+2\in U_{\beta}}u_{\beta}(n-j+2)=\sum_{j'\in U_{\beta}}u_{\beta}(j'),
\end{equation*}
which means that $\rho'$ minimizes the number of backward steps taken in the problem for $\beta$. By symmetry, $\rho$ \emph{minimizes the number of forward steps} taken in the original problem for $\alpha$, and the minimal number of forward steps is $-\sum_{j\notin\Ua}\ua(j)=\sum_{j\notin\Ua}\abs{\ua(j)}$.

\begin{theorem}
	Let $\alpha\in\pp$, and consider a vector $\rho\in[0,n]^n$ such that $\abs{\rho}=\sum_{j\in\Ua}\ua(j)$ (i.e., $\rho$ is constructed as shown is Theorem \ref{choiceminimize}). Then $\rho$ minimizes the total number of steps taken by the cars, which amounts to $\sum_{j=1}^n \abs{\ua(j)}$.
\end{theorem}

\begin{remark}
	The previous construction is significant, as it highlights the dual nature of the problem (in cases where each car takes either only backward steps, or only forward steps). However, the fact that vectors minimize forward steps if and only if they minimize backward steps, can also be shown in a more straightforward way. Let $(\alpha,\rho)\in PR_n$, and denote  $f_{\alpha}(\rho)$ as the number of forward steps taken by the cars. Just like we have previously observed, defining $\tilde{\alpha}$ by $\tilde{a}_i=a_i-r_i$ for all $i\in[n]$ yields a parking function, such that $f_{\tilde{\alpha}}(\underline{0})=f_{\alpha}(\rho)$. Thus,
	\[
	f_{\alpha}(\rho)=f_{\tilde{\alpha}}(\underline{0})=\sum_{i=1}^n i -\sum_{i=1}^n \tilde{a}_i
	=\sum_{i=1}^n i -\sum_{i=1}^n a_i +\sum_{i=1}^n r_i.
	\]
	Observe that the values $\sum_{i=1}^n=\frac{n(n+1)}{2}$ and $\sum_{i=1}^n a_i$ are fixed by the choice of $\alpha$, thus minimizing $f_{\alpha}(\rho)$ corresponds to minimizing $\abs{\rho}=\sum_{i=1}^n r_i$. As we have seen above, in that case we get $f_{\alpha}(\rho)=-\sum_{j\notin\Ua}\ua(j)$.
\end{remark}

\begin{remark}
	Let $\alpha\in\pp$, and $\beta\in\pp$ defined by $b_i=n+1-a_i$ for all $i\in[n]$. By Corollary \ref{minimalrhocount}, the number of vectors $\rho\in[0,n]^n$ that minimize the total number of steps taken in $(\alpha,\rho)$ is $\prod_{i=1}^m \binom{\A_{q_i}-1}{\ua(q_i)}$, where $\Ua=\bigcup_{i=1}^m [p_i,q_i]$ is the decomposition of $\Ua$ in maximal intervals. This also implies that the number of vectors $\rho'\in[0,n]^n$ that minimize the total number of steps taken in $(\beta,\rho')$ is the same.\\
	We can check that directly. Recall that $\abs{\beta}_i=\A_{n+1-i}$ for all $i\in[n]$, and $u_{\beta}(j)=-\ua(n-j+2)$ for all $j\in[2,n]$. Consider a maximal interval $[p,q]\subset U_{\beta}$ (thus, if $q<n$, $u_{\beta}(q+1)\le0$). There are three cases:
	\begin{itemize}
		\item if $q=n$, then by definition $u_{\beta}(n)=\abs{\beta}_n-1$, thus there is only one way to choose $u_{\beta}(n)$ cars with preference $n$, excluding the first;
		\item if $q<n$, and $u_{\beta}(q+1)=0$, then by \eqref{propua3}, $u_{\beta}(q)=u_{\beta}(q+1)+\abs{\beta}_q-1=\abs{\beta}_q-1$, and again there is only one way to choose $u_{\beta}(n)$ cars with preference $n$, excluding the first;
		\item if $q<n$, and $u_{\beta}(q+1)<0$, then $\ua(n-q+2)=-u_{\beta}(q)<0$, and $\ua(n-q+1)=-u_{\beta}(q+1)>0$, thus $n-q+1$ is the upper extreme of a maximal interval in $\Ua$. Moreover, $\A_{n+1-q}=\abs{\beta}_{q}$, and $u_{\beta}(q)=u_{\beta}(q+1)+\abs{\beta}_q-1=\abs{\beta}_q-1$, thus
		\begin{equation*}
		\binom{\abs{\beta}_q-1}{u_{\beta}(q)}=\binom{\abs{\beta}_q-1}{\abs{\beta}_q-1-(-u_{\beta}(q+1)}
		=\binom{\abs{\beta}_q-1}{-u_{\beta}(q+1)}=\binom{\A_{n+1-q}-1}{\ua(n+1-q)}.
		\end{equation*}
	\end{itemize}
	Thus, if there is more than one way to choose $u_{\beta}(q)$ among $\abs{\beta}_q$ cars with preference $q$ in $\beta$, excluding the first, then there is a maximal interval in $\Ua$, with upper extreme $n-q+1$, such that $\binom{\abs{\beta}_q-1}{u_{\beta}(q)}=\binom{\A_{n+1-q}-1}{\ua(n+1-q)}$. Reversing the roles of $\alpha$ and $\beta$, for each factor $\binom{\A_{q_i}-1}{\ua(q_i)}$ strictly greater than one in $\prod_{i=1}^m \binom{\A_{q_i}-1}{\ua(q_i)}$, there is also a corresponding factor $\binom{\abs{\beta}_{n+1-q_i}-1}{u_{\beta}(n+1-q_i)}$ for $\beta$. As result, the number of ways to choose a minimizing $\rho'$ for $\beta$ is equal to the number of ways to choose a minimizing $\rho$ for $\alpha$.
\end{remark}

\begin{remark}
	Note that the values $\ua(j)$, for $j\in[n]$, are independent of the order of a parking preference $\alpha\in\pp$. As a result, the minimal number of steps that need to be taken by the cars ($\sum_{j=1}^n \abs{\ua(j)}$) remains unchanged if the parking preference is rearranged. Furthermore, by Corollary \ref{minimalrhocount}, the number of vectors that minimize the number of steps taken  remains unchanged as well.
\end{remark}

\section{Strategies minimizing the number of cars driving backwards}
Consider the following problem: let $\alpha\in\pp$, what is the \emph{minimum number of cars} that must be allowed to drive backwards, so that all cars are able to park, and how can we choose them? In other words, how can we find a parking strategy $\rho\in\{0,n\}^n$ for $\alpha$, such that it minimizes $\abs{\rho}$?

Note that this is different from the problem analysed in the previous section. For example, let $\alpha=(3,3,2,4)$: the vector that minimizes the number of steps is $(0,1,1,0)$, but the one that minimizes the number of cars is $(0,0,0,3)$.

The intuitive way to minimize the number of cars driving backwards is to allow a car $c_i$ to drive backwards if and only if all spots in $[a_i,n]$ are occupied. As we will show, this kind of "greedy" strategy is indeed the correct way to make the choice (or at least, it is one of the correct ways).

\begin{theorem}
	\label{minimizecars}
	Let $\alpha\in\pp$. For each $j\in\Ua$, define $\tilde{T}_j$ as the set of the last $\ua(j)$ cars with preference at least $j$, and let $\tilde{T}=\bigcup_{j\in\Ua}\tilde{T}_j$. Then if $\rho\in\{0,n\}^n$ is such that $(\alpha,\rho)\in\pr$, then $\abs{\rho}\ge n\abs{\tilde{T}}$.\\ 
	Moreover, define $\rho'\in\{0,n\}^n$ by $r_i'=n$ if $i\in\tilde{T}$, and $r_i'=0$ otherwise. Then $(\alpha,\rho')\in\pr$.
\end{theorem}

\begin{proof}
	Let $\tilde{T}_j$, $j\in\Ua$, and $\tilde{T}$ be defined as above. Let $\rho\in\{0,n\}^n$, and consider $(\alpha,\rho)$. For all $i\in[n]$, define 
	\[
	\lambda_i=\min\{j\in\Ua\mid \tilde{T}_j\cap[1,i]\ne\varnothing\},
	\]
	where if $\{j\in\Ua\mid \tilde{T}_j\cap[1,i]\ne\varnothing\}=\varnothing$, we set $\lambda_i=\infty$. Moreover, note that for all $i\in[2,n]$, $\lambda_i\le\lambda_{i-1}$.
	\begin{comment}
	As a preliminary observation, we prove the following: let $i\in[n]$, if the first $i-1$ cars fill spots $[a_i,n]$, then $i\in\tilde{T}$. Let $h$ be the highest unoccupied spot by the start of car $c_i$'s turn to park, then clearly $h< a_i$. Spots $[h+1,n]$ have been occupied by $n-h$ cars with preference in $[h+1,n]$, since spot $h$ has remained unoccupied. Car $c_i$ also has preference in $[h+1,n]$, therefore
	\[
	\ua(h+1)=\sum_{k=h+1}^n\A_k-(n-(h+1)+1)\ge n-h+1-(n-h)=1.
	\]
	Thus, $h+1\in\Ua$. Moreover, since $a_i\ge h+1$, and there have already been at least $n-h$ cars with preference in $[h+1,n]$, car $c_i$ is one of the last $\sum_{k=h+1}^n\A_k-(n-h)=\ua(h+1)$ cars with preference at least $h+1$, so we get $i\in\tilde{T}_{h+1}\subseteq\tilde{T}$.
	\end{comment}
	
	We will prove by induction that, for all $i\in[n]$, if the first $i$ cars are able to park, then
	\begin{equation}
		\label{mincarinduction}
	\abs{\{j\le i\mid a_j\ge \lambda_i,\, r_j=n \}}\ge\abs{\tilde{T}\cap [1,i]},
	\end{equation}
	and furthermore, if $\lambda_i\le n$ and $\abs{\{j\le i\mid a_j\ge \lambda_i,\, r_j=n \}}=\abs{\tilde{T}\cap [1,i]}$, then $[\lambda_i,n]$ is filled by the first $i$ cars.\\

	If $i=1$, then car $c_i$ parks in spot $a_i$. Moreover, let $j\in\Ua$ such that $j\le i$, then 
	\[
	\abs{\{k\in[n]\mid a_k\in[j,n]\}}=\sum_{k=j}^n\A_k=\ua(j)+(n-j+1),
	\]
	where $n-j+1\ge1$, thus $c_i$ is not one of the first $\ua(j)$ cars with preference at least $j$, so $i\notin\tilde{T}_j$. Therefore, $\tilde{T}\cap[1]=\varnothing$, and \eqref{mincarinduction} is trivially satisfied.

	Let $i\in[2,n]$, and suppose that the first $i$ cars park successfully, and that \eqref{mincarinduction} holds for all values in $[1,i-1]$. 
	\begin{itemize}
	\item
	If $i\notin\tilde{T}$, then since $\lambda_i=\lambda_{i-1}$, and we get
	\begin{multline*}
	\abs{\{j\le i\mid a_j\ge \lambda_i,\, r_j=n \}}\\
	\ge\abs{\{j\le i-1\mid a_j\ge \lambda_{i-1},\, r_j=n \}}\ge\abs{\tilde{T}\cap [1,i-1]}=\abs{\tilde{T}\cap [1,i]}.
	\end{multline*}
	Moreover, if $\lambda_i\le n$ and $\abs{\{j\le i\mid a_j\ge \lambda_i,\, r_j=n \}}=\abs{\tilde{T}\cap [1,i]}$, then also $\abs{\{j\le i-1\mid a_j\ge \lambda_{i-1},\, r_j=n \}}=\abs{\tilde{T}\cap [1,i-1]}$, thus by induction all spots in $[\lambda_i,n]=[\lambda_{i-1},n]$ are filled.
	\item
	Now assume that $i\in\tilde{T}$. In particular, we show that in this case, $i\in\tilde{T}_{\lambda_i}$. Clearly $i\in\tilde{T}_j$ for some $j\in\Ua$, thus, by definition of $\lambda_i$ and $\tilde{T}_j$, $a_i\ge j\ge\lambda_i$. Moreover, if $i\notin\tilde{T}_{\lambda_i}$, then there exists $h\in[1,i-1]$ such that $h\in\tilde{T}_{\lambda_i}$. Thus, $h$ is one of the last $\ua(\lambda_i)$ cars with preference at least $\lambda_i$, but car $c_i$ ($i>h$) with preference $a_i\ge\lambda_i$ is not, which is absurd.\\
	There are now two cases: either the first $i-1$ cars have completely filled spots $[\lambda_i,n]$, or they have not. 
	\begin{itemize}
	\item
	In the former case, for car $c_i$ to be able to park, we must needs have $r_i=n$, thus
	\begin{multline*}
	\abs{\{j\le i\mid a_j\ge \lambda_i,\, r_j=n \}}=1+\abs{\{j\le i-1\mid a_j\ge \lambda_i,\, r_j=n \}}\\
	\ge1+\abs{\{j\le i-1\mid a_j\ge \lambda_{i-1},\, r_j=n \}}\ge1+\abs{\tilde{T}\cap [1,i-1]}
	=\abs{\tilde{T}\cap [1,i]}.
	\end{multline*}
	Furthermore, the second condition for the induction holds trivially for $i$.
	\item
	On the other hand, suppose that there exits a spot $h\in[\lambda_i,n]$ that is not filled by one of the first $i-1$ cars (and thus the condition $r_i=n$ may not be required). We consider two cases: either $\abs{\tilde{T}\cap[1,i-1]}=0$ or $\abs{\tilde{T}\cap[1,i-1]}\ge1$.
	\begin{itemize}
	\item
	If $\abs{\tilde{T}\cap[1,i-1]}=0$, suppose that $r_j=0$ for all $j\in[1,i-1]$. Since car $c_i$ is the first of the last $\ua(\lambda_i)$ cars with preference at least $\lambda_i$, there have already been exactly $\sum_{j=\lambda_i}^n\A_j-\ua(\lambda_i)=n-\lambda_i+1$ cars with preference in $[\lambda_i,n]$. Since all of those cars have been able to park, $\abs{[\lambda_i,n]}=n-\lambda_i+1$, and none among those cars park backwards, then all spots in $[\lambda_i,n]$ have been filled by the first $i-1$ cars, and we have a contradiction. Thus, $r_j=n$ for some $j\in[1,i-1]$, and consequently
	\[
	\abs{\{j\le i\mid a_j\ge \lambda_i,\, r_j=n \}}\ge1=\abs{\tilde{T}\cap [1,i]}.
	\]
%	Moreover, if $\abs{\{j\le i\mid a_j\ge \lambda_i,\, r_j=n \}}=1$, then $r_j=n$ for only one $j\in[1,i-1]$, and $r_i=0$. Thus, the first $i-1$ cars leave exactly one free spot in $[\lambda_i,n]$, and car $c_i$ will fill it (car $c_i$ is able to park forwards, and $\lambda_i\le a_i$), so all spots in $[\lambda_i,n]$ are filled by the first $i$ cars.
	\item
	Suppose instead that $\abs{\tilde{T}\cap[1,i-1]}\ge1$.
	\begin{itemize}
	\item 
	If $\lambda_{i-1}=\lambda_i$, by induction, since not all spots in $[\lambda_i,n]$ are filled, necessarily
	\[
	\abs{\{j\le i-1\mid a_j\ge \lambda_{i-1},\, r_j=n \}}>\abs{\tilde{T}\cap [1,i-1]},
	\]
	and thus,
	\begin{multline*}
	\abs{\{j\le i\mid a_j\ge \lambda_i,\, r_j=n \}}\\
	\ge\abs{\{j\le i-1\mid a_j\ge \lambda_{i-1},\, r_j=n \}}\\
	\ge1+\abs{\tilde{T}\cap [1,i-1]}=\abs{\tilde{T}\cap [1,i]}.
	\end{multline*}
	\item
	Assume $\lambda_{i-1}>\lambda_i$ (note that $\lambda_{i-1}<\infty$ since $\abs{\tilde{T}\cap[1,i-1]}\ge1$). If $[\lambda_{i-1},n]$ is not filled by the first $i-1$ cars, then by induction  $\abs{\{j\le i-1\mid a_j\ge \lambda_{i-1},\, r_j=n \}}>\abs{\tilde{T}\cap [1,i-1]}$, and \eqref{mincarinduction} follows just like in the previous case. Otherwise, $[\lambda_{i-1},n]$ is filled by the first $i-1$ cars, therefore $h\in[\lambda_i,\lambda_{i-1}-1]$. Since car $c_i$ is the first of the last $\ua(\lambda_i)$ cars with preference in $[\lambda_i,n]$, there have already been exactly $n-\lambda_i+1$ cars with preference in $[\lambda_i,n]$ and since there is an empty spot in $[\lambda_i,\lambda_{i-1}-1]$, and all cars were able to park, a car \emph{with preference in $[\lambda_i,\lambda_{i-1}-1]$} (since not all spots in $[\lambda_i,\lambda_{i-1}-1]$ have been occupied) has parked by driving backwards. Thus,
	\begin{multline*}
	\abs{\{j\le i\mid a_j\ge \lambda_i,\, r_j=n \}}\\
	\ge 1+\abs{\{j\le i\mid a_j\ge \lambda_{i-1},\, r_j=n \}}\\
	\ge1+\abs{\{j\le i-1\mid a_j\ge \lambda_{i-1},\, r_j=n \}}\\
	\ge1+\abs{\tilde{T}\cap [1,i-1]}=\abs{\tilde{T}\cap [1,i]}.
	\end{multline*}
	\end{itemize}
	\end{itemize}
	\end{itemize}

	Now, suppose that $\abs{\{j\le i\mid a_j\ge \lambda_i,\, r_j=n \}}=\abs{\tilde{T}\cap [1,i]}$. Furthermore, suppose that not all spots in $[\lambda_i,n]$ are filled by the first $i$ cars, and let $h$ be the smallest unoccupied spot in $[\lambda_i,n]$. Since $i\in\tilde{T}_{\lambda_i}$, car $c_i$ is one of the last $\ua(\lambda_i)$ cars with preference at least $\lambda_i$, at least $\sum_{j=i}^n\A_{\lambda_i}-\ua(\lambda_i)=n-\lambda_i+1$ cars among the first $i-1$ cars had preference in $[\lambda_i,n]$, and since $a_i\in[\lambda_i,n]$, the total number of cars with preference in $[\lambda_i,n]$ among the first $i$ cars is at least $n-\lambda_i+2$. There are only $n-\lambda_i+1$ spots in $[\lambda_i,n]$, thus there has been at least one car with preference in $[\lambda_i,n$ that has parked in $[1,\lambda_i-1]$ by driving backwards; furthermore, since $h\in[\lambda_i,n]$ is an unoccupied spot, the preference of that car must have been in $[\lambda_i,h-1]$.\\
	Let $m\in[1,i]$ be such that $c_m$ is the last car (among the first $i$ cars) with preference in $[\lambda_i,h-1]$ that parks in $[1,\lambda_i-1]$. Clearly, $r_m=n$. Define $\tilde{\rho}$ by $\tilde{r}_m=0$, and $\tilde{r}_i=r_i$ for all $i\ne m$, and consider $(\alpha,\tilde{\rho})$. Trivially, all of the first $m-1$ cars park in the same spots they occupied in the original problem. Focus on all cars with indices in $[m,i]$, and with preference smaller than $h$. In $(\alpha,\rho)$:
	\begin{itemize}
		\item the first of those cars (i.e. car $c_m$) is such that it has preference $h>a_m\ge\lambda_i$, and parks in $[1,\lambda_i-1]$: thus, all spots in $[\lambda_i,a_m]$ have already been occupied by the first $m-1$ cars;
		\item among the others, there are $k_1$ cars with preference in $[a_m+1,h]$: by definition of $h$ and $m$, and by what we have just remarked, these cars all park in $[a_m+1,h-1]$;
		\item there are $k_2$ cars with preference in $[1,a_m]$ that reach spot $a_m$ by driving forwards, and thus park in $[a_m+1,h-1]$.
	\end{itemize}
	Thus, $k_1+k_2=\abs{[a_m+1,h-1]}$. As a result, focusing on the same cars in $(\alpha,\tilde{\rho})$:
	\begin{itemize}
		\item the first car (i.e. $c_m$) parks in $[a_m+1,h]$;
		\item among the others, there are $k_1$ cars with preference in $[a_m+1,h]$;
		\item there are at most $k_2$ cars with preference in $[1,a_m]$ that reach spot $a_m$ by driving forwards.
	\end{itemize}
	Thus, there are at most $1+k_1+k_2=1+\abs{[a_m+1,h-1]}\abs{[a_m,h]}$ cars parking in $[a_m,h]$. In particular, in $(\alpha,\tilde{\rho})$, among the first $i$ cars, all cars with preference in $[1,h]$ are able to park in $[1,h]$. Consequently, among the first $i$ cars,  all those with preference in $[h+1,n]$, occupy the same spots they did in the original problem.\\
	Thus, the first $i$ cars in $(\alpha,\tilde{\rho})$ are able to park. However, by definition of $\tilde{\rho}$:
	\[
	\abs{\{j\le i\mid a_j\ge \lambda_i,\, \tilde{r}_j=n \}}<\abs{\{j\le i\mid a_j\ge \lambda_i,\, r_j=n \}}=\abs{\tilde{T}\cap [1,i]},
	\]
	which contradicts \eqref{mincarinduction}. Thus, we have proved the second part of the induction.
	\end{itemize}
	As we have seen, \eqref{mincarinduction} holds for all $i\in[n]$. In particular, for $i=n$ we get
	\begin{equation*}
	\abs{\{j\le n\mid r_j=n \}}\ge\abs{\{j\le n\mid a_j\ge \lambda_i,\, r_j=n \}}
	=\abs{\tilde{T}\cap [1,n]}=\abs{\tilde{T}},
	\end{equation*}
	and consequently, $\abs{\rho}=n\abs{\{j\le n\mid r_j=n \}}\ge n\abs{\tilde{T}}$.\\
	
	Finally, define $\rho'\in\{0,n\}^n$ by $r_j'=n$ if $j\in\tilde{T}$ and $r_j'=0$ otherwise, and let $i\in[n]$. If the first $i-1$ cars have not occupied all spots in $[a_i,n]$, then clearly car $c_i$ is able to park. Otherwise, suppose that $[a_i,n]$ is filled by the first $i-1$ cars. Let $h$ be the largest unoccupied spot by the start of car $c_i$'s turn to park, then clearly $h< a_i$. Spots $[h+1,n]$ have been occupied by $n-h$ cars with preference in $[h+1,n]$, since spot $h$ has remained unoccupied. Car $c_i$ also has preference in $[h+1,n]$, therefore
	\[
	\ua(h+1)=\sum_{k=h+1}^n\A_k-(n-(h+1)+1)\ge n-h+1-(n-h)=1.
	\]
	Thus, $h+1\in\Ua$. Moreover, since $a_i\ge h+1$, and there have already been at least $n-h$ cars with preference in $[h+1,n]$, car $c_i$ is one of the last $\sum_{k=h+1}^n\A_k-(n-h)=\ua(h+1)$ cars with preference at least $h+1$, so we get $i\in\tilde{T}_{h+1}\subseteq\tilde{T}$. Therefore, $r_i'=n$ and car $c_i$ is able to park. Thus, $(\alpha,\rho')\in\pr$. 
\end{proof}

\begin{example}
	Let $\alpha=(4,9,5,9,5,8,7,9,4,6)$, where $\Ua=[2,9]$. We get:
	\[
	\begin{array}{ccccc}
		\toprule
		&\ua(j)	&\{i\in[n]\mid a_i\ge j\}	&\tilde{T}_j\\
		\midrule
		j=9	&1		&\{2,4,8\}				&\{8\}\\
		j=8	&1		&\{2,4,6,8\}			&\{8\}\\
		j=7	&1		&\{2,4,6,7,8\}			&\{8\}\\
		j=6	&1		&\{2,4,6,7,8,10\}		&\{10\}\\
		j=5	&2		&\{2,3,4,5,6,7,8,10\}	&\{8,10\}\\
		j=4	&3		&[1,10]					&\{8,9,10\}\\
		j=3	&2		&[1,10]					&\{9,10\}\\
		j=2	&1		&[1,10]					&\{10\}\\
		\bottomrule
	\end{array}.
	\]
	Thus, $\tilde{T}=\bigcup_{j\in\Ua} \tilde{T}_j=\{8,9,10\}$, so at least three cars need to drive backwards for all cars to be able to park. Let $\rho=(0,0,0,0,0,0,0,n,n,n)$, that is, $\rho$ is the specific vector found in Theorem \ref{minimizecars}. Then
	\[
	\psi(4959587946,0000000nnn)=(4,9,5,10,6,8,7,3,2,1).
	\]
	As we can observe, for all $i\in[1,7]=[1,10]\setminus\tilde{T}$, car $c_i$ is able to park driving forwards. On the other hand, car $c_8$ needs to drive backwards since all spots in $[9,10]=[a_8,10]$ are filled, and the same goes for cars $c_9$ and $c_{10}$. Furthermore, the total number of backward steps actually taken is $\sum_{j\in\tilde{T}}(a_i-\psi(c_i))=6+2+5=13$, which is larger than $\sum_{j\in\Ua}\ua(j)=12$, and that is why the vector $(0,0,0,0,0,0,0,6,2,5)$ could not be found using the method shown in the previous section (minimizing $\abs{\rho}$), like in Examples \ref{examplemin1} and \ref{examplemin2}.
\end{example}

\begin{remark}
	The parking strategy constructed in Theorem \ref{minimizecars} is not necessarily the only one such that the least number of cars drive backwards. For example, let $\alpha=(3,3,3,2,3,3,3)\in PP_7$: it is easy to check that in this case, $\tilde{T}=\{7\}$, and indeed $(3332333,0000007)\in PR_7$. However, in this case we also get $(3332333,0000700)\in PR_7$ and $(3332333,00000070)\in PR_7$.\\
	Furthermore, note that rearranging the parking preference changes the number of vectors that minimize the number of cars driving backwards: 
	\begin{itemize}
		\item considering $\alpha_1=(3,3,3,3,3,2,3)$, $\rho_1=(0,0,0,0,0,0,7)\in\{0,7\}^7$ is the only parking strategy for $\alpha_1$ such that $\abs{\rho_1}=7$;
		\item considering $\alpha_2=(2,3,3,3,3,3,3)$, $\rho_{2,i}=(0,0,r_{3,i},r_{4,i},r_{5,i},r_{6,i},r_{7,i})\in\{0,7\}^7$, for $i\in[3,7]$, defined by $r_{j,i}=7$ if and only if $j=i$, are five different parking strategies for $\alpha_2$ such that $\abs{\rho_{2,i}}=7$.
	\end{itemize} 
	Moreover, rearranging the parking preference can also change the size of $\tilde{T}$: for example, $\alpha_3=(3,3,3,3,3,3,2)$ is such that $\tilde{T}=\{6,7\}$, and thus at least $2$ cars need to drive backwards for all cars to be able to park.
\end{remark}

\section{The case $\rho\in[0,1]^n$}

\label{section1nap}

It is interesting to study what happens if we put some restrictions on $\rho$: specifically, we are going to analyse the case where the rules in $\rho$ can only be either the standard parking rule, or the $1$-Naples parking rule (i.e. $\rho\in[0,1]^n$). \\
We begin by noting that for $k=1$, complete parking preferences have an especially simple form.

\begin{lemma}
	\label{cardlemma}
	Let $\alpha\in\pp$ be such that $\ua(j)\le1$ for all $j\in[n]$. Let $[p,q]\subseteq\Ua$ be a maximal interval, then $\A_q\ge2$, and $\A_j=1$ for all $p\le j<q$. In particular, if $\alpha$ is complete then $\A_n=2$, and $\A_j=1$ for all $j\in[2,n-1]$.
\end{lemma}

\begin{proof}
	By Proposition \ref{propUa}, $\A_q\ge2$. For all $j\in[p,q]$, $j\in\Ua$ so $\ua(j)\ge1$, thus using the hypothesis $\ua(j)=1$. Hence for all $j\in[p,q-1]$
	\[
	\A_j=\ua(j)-\ua(j+1)+1=1.
	\]
	Finally, if $\alpha$ is complete then $n\in\Ua$, thus $\ua(n)=1$ and by definition $\A_n=\ua(n)+1=2$.
\end{proof}

\begin{lemma}
	\label{char1napUcompl}
	Let $\alpha\in\nap$ be complete. Then $a_1=n$, and $a_j=n+2-j$ for all $j\in[2,n]$. In other words, $\alpha=(n,n,n-1,n-2,\dots,3,2)$.
\end{lemma}

\begin{proof}
	Note that since $\Ua=[2,n]$, by Theorem \ref{necessk} $\ua(j)=1$ for all $j\in[2,n]$. Since $\alpha\in\nap$, by Theorem \ref{charUcomp}, for all $h\in[2,n]$ there exists $0\le\lambda_h\le 1-\ua(h)=0$ such that
	\[
	\abs{\{i<j_h\mid a_i=h-1\}}=0,
	\]
	where $j_h$ is the last index such that $a_{j_h}=h$. Since for all $h\in[3,n]$, $\A_{h-1}\ne0$ by the previous Lemma, this implies that the sequence $a_1,\dots,a_n$ is in decreasing order. We conclude by remarking that again by the previous Lemma, $\A_n=2$ and $\A_j=1$ for all $j\in[2,n-1]$, so $\alpha=(n,n,n-1,\dots,3,2)$.
\end{proof}

\begin{comment}
\begin{remark}
	Note that if $\alpha=(n,n,n-1,\dots,3,2)$, then it is easy to see that $\alpha\in\nap$ and is complete.
\end{remark}
\end{comment}

\begin{theorem}
	\label{char1nap}
	Let $\alpha\in\pp$. Then $\alpha\in\nap$ if and only if for all maximal intervals $[p,q]\subseteq\Ua$ there exist $m$ indices $j_1<\dots<j_m$, with $m\ge q-p+2$, such that $a_{j_1}=a_{j_2}$, and for all $i=2,\dots,m$, $a_{j_i}=m-i+p$.
\end{theorem}

\begin{proof}
	By Theorems \ref{theosuff} and \ref{theoexist}, $\alpha\in\nap$ if and only if for all maximal intervals $[p,q]\subseteq\Ua$ there exists $J\subseteq[n]$ such that $\tau_{p-2}(\alpha_{|_J})$ is complete, and $\tau_{p-2}(\alpha_{|_J})\in PF_{\abs{J},1}$. Let $j_1<\dots<j_m$ be such that $J=\{j_1,\dots,j_m\}$, then by Lemma \ref{char1napUcompl}
	\[
	\tau_{p-2}(\alpha_{|_J})=(m,m,m-1,\dots,3,2),
	\]
	and consequently $a_{j_i}=m+2-i+p-2=m-i+p$ for all $i=2,\dots,m$, and $a_{j_1}=a_{j_2}$.
\end{proof}

This result shows that in the $1$-Naples case we have a much simpler characterization than for the general $k$-Naples parking functions. Moreover, Theorem \ref{char1nap} can be slightly modified to address the problem for $(\alpha,\rho)$ with $\rho\in[0,1]^n$. The proof we provide for the following theorem also doubles as an alternative, more direct way to prove Theorem \ref{char1nap}, without needing to rely on most of the more general results found in the previous chapters.

\begin{theorem}
	\label{char1napchoice}
	Let $\alpha\in\pp$ and $\rho\in[0,1]^n$. Then $(\alpha,\rho)\in\pr$ if and only if for all maximal intervals $[p,q]\subseteq\Ua$ there exist $m$ indices $j_1<\dots<j_m$, with $m\ge q-p+2$, such that $a_{j_1}=a_{j_2}$, and for all $i=2,\dots,m$, $a_{j_i}=m-i+p$ and $r_{j_i}=1$.
\end{theorem}

\begin{proof}
	Suppose that for all maximal intervals $[p,q]\subseteq\Ua$, there exist $m$ indices $j_1<\dots<j_m$, as above. Let $[p,q]\subseteq\Ua$, and $j_1,\dots,j_m$ be such indices. After car $c_{j_1}$ attempts to park, spot $a_{j_1}$ is certainly occupied; thus, on car $c_{j_2}$'s turn, it finds spot $a_{j_2}=a_{j_1}$ occupied, and since $r_{j_2}=1$, it drives backwards to spot $a_{j_2}-1=a_{j_3}$, filling it if necessary.\\
	By induction, suppose that for all $i\in[2,h-1]$, after car $c_{j_i}$'s turn spot $a_{j_{i+1}}$ is certainly occupied. Then car $c_{j_h}$ finds spot $a_{j_h}$ occupied, and since $r_{j_h}=1$, it fills spot $a_{j_h}-1$ if necessary (and $a_{j_h}-1=a_{j_{h+1}}$ if $h<m$). Hence, the induction holds for all $h\in[2,m]$, and in particular spot $a_{j_m}-1=m-m+p-1=p-1$ is certainly filled, and this suffices to prove $(\alpha,\rho)\in\pr$.\\
	Suppose $(\alpha,\rho)\in\pr$, and let $[p,q]\subseteq\Ua$ be a maximal interval. Note that $\A_{p-1}=0$, and $\ua(p-1)=0$. Since $\rho\in[0,1]^n$, no car can drive backwards for more than one spot, so all cars with preference in $[p,n]$ can never reach any spot prior to $(p-1)$; conversely, note that $\sum_{i=1}^{p-2}\A_i=p-2-\ua(p-1)=p-2$, thus since $(\alpha,\rho)\in\pr$ all cars with preference in $[1,p-2]$ have to park in $[1,p-2]$ to be able to fill it, so $[p-1,n]$ is all occupied by cars with preference in $[p,n]$.\\
	Note that by Lemma \ref{greaterlemma}, $(\alpha,\rho)\in\pr$ implies $(\alpha,(1,1,\dots,1))\in\pr$, and thus $\alpha\in\nap$. In particular by Theorem \ref{necessk}, $\ua(j)=1$ for all $j\in[p,q]$ and thus by Lemma \ref{cardlemma}, $\A_q\ge2$, and $\A_j=1$ for all $j\in[p,q-1]$.\\
	Consider spot $(p-1)$: it must be filled by a car having preference in $[p,n]$, and since the maximum number of backward steps a single car can make is $1$, it must be filled by a car with preference $p$. Thus, there exists a car $c_{h_1}$ such that $a_{h_1}=p$, $r_{h_1}=1$ and $\psi(c_{h_1})=p-1$.\\
	Suppose that for all $j<h_1$, $a_j\ne p$. Since car $c_{h_1}$ must find spot $p$ already occupied to be able to drive backwards; furthermore, spot $p$ cannot have been occupied by a car with preference $p$ (because $a_j\ne p$ for all $j<h_1$), nor by a car with preference at most $(p-1)$, otherwise spot $(p-1)$ would have been already occupied by the time car $c_{h_1}$ parks. Thus, spot $p$ must have been occupied by a car $c_{h_2}$ ($h_2<h_1$) such that $a_{h_2}=p+1$, $r_{h_2}=1$ and $\psi(c_{h_2})=p$. We can therefore iterate this argument on $h_1$.\\
	Suppose that there exist $h_{\lambda}<h_{\lambda-1}<\dots<h_1$ such that for each $i\in[1,\lambda]$, car $c_{h_i}$ is such that $a_{h_i}=p+i-1$, $r_{h_i}=1$ and $\psi(c_{h_i})=p+i-2$. Suppose $j<h_{\lambda}$, $a_{j}\ne p+\lambda-1$: spot $p+\lambda-1$ must already be occupied before car $c_{h_\lambda}$ parks, moreover it cannot have been occupied by a car with preference $p+\lambda-1$, nor by a car with preference at most $p+\lambda-2$ (otherwise spot $p+\lambda-2$ would also already be occupied). Thus, there exists a car $c_{h_{\lambda+1}}$ ($h_{\lambda+1}<h_{\lambda}$) such that $a_{h_{\lambda+1}}=p+\lambda$, $r_{h_{\lambda+1}}=1$ and $\psi(c_{h_{\lambda+1}})=p+\lambda-1$.\\
	The sequence $h_1>h_2>\dots$ cannot be infinitely long, since there are only $n$ cars, so at some point we will get $h_1>h_2>\dots>h_m$ such that $a_{h_m}=p+m-2=a_{h_{m-1}}$, and for all $i\in[1,m-1]$, $a_{h_i}=p+i-1$, $r_{h_i}=1$ and $\psi(c_{h_i})=p+i-2$. Substituting $h_i=j_{m+1-i}$ for each $i\in[1,m]$, we get a sequence $j_1<\dots<j_m$ such that $a_{j_1}=a_{j_2}$, and for all $i\in[2,m]$, $a_{j_i}=p+m-i$ and $r_{j_i}=1$. We conclude by remarking that since $\A_j=1$ for all $j\in[p,q-1]$, then $a_{j_1}=a_{j_2}=p+m-2\ge q$, and thus $m\ge q-p+2$.
\end{proof}

\begin{example}
	Let $\alpha=(5,10,11,1,5,11,10,4,3,9,10,8,3)\in PP_{13}$, where $\Ua=\{3\}\cup[8,10]$. It is a $1$-Naples parking function, and the outcome map gives $\psi_1(\alpha)=(5,10,11,1,\\4,12,9,3,2,8,13,7,6)$.
	Consider $\rho=(0,0,0,0,0,1,1,0,0,1,0,1,1)$. For the maximal interval $[3,3]\subseteq\Ua$, we have the subsequence $\alpha_{|_{\{9,13\}}}=(3,3)$, where $r_{13}=1$. For the other maximal interval $[8,10]\subseteq\Ua$, $\alpha_{|_{\{3,6,7,10,12\}}}=(11,11,10,9,8)$ is such that $r_6=r_7=r_{10}=r_{12}=1$. Moreover, note that there also exists $\alpha_{|_{\{2,7,10,12\}}}=(10,10,9,8)$, such that $r_7=r_{10}=r_{12}=1$. Thus, by Theorem \ref{char1napchoice}, $(\alpha,\rho)\in PR_{13}$. Checking manually, we get $\psi(\alpha,\rho)=(5,10,11,1,6,12,9,4,3,8,13,7,2)$.
\end{example}

\begin{remark}
	Let $\alpha\in\pp$, such that $\ua(j)\le1$ for all $j\in[n]$, and let $\Ua=\bigcup_{i=1}^m[p_i,q_i]$ the decomposition of $\Ua$ in maximal intervals. Let 
	\[
	A_i=\{j\in[n]\mid a_j\in[p_i-1,p_{i+1}-2] \}
	\]
	for all $i=1,\dots,m$, and $A_m=\{j\in[n]\mid a_j\in[p_m-1,n] \}$. Then $\alpha\in\nap$ if and only if $\alpha_{|_{A_i}}$ is a $1$-Naples parking function for all $i=1,\dots,m$, when properly translated. This follows from the fact that, given a maximal interval $[p,q]\subseteq\Ua$, no car with preference at least $p-1$ can park in any spot in $[1,p-2]$, since $\A_{p-1}=0$ and cars can drive backwards at most one spot. In the same way, $(\alpha,\rho)\in\pr$ if and only if $(\tau_{p-2}(\alpha_{|_{A_i}}),\rho_{|_{A_i}})\in PR_{\abs{J_i}}$ for all $i=1,\dots,m$.\\
	If we want to find a parking strategy $\rho\in[0,1]^n$ that minimizes $\abs{\rho}$, we can therefore focus on each $\alpha_{|_{A_i}}$. In this case, taking the shortest sequence satisfying Theorem \ref{char1napchoice} for $[p_i,q_i]\subseteq\Ua$, and then setting $r_j=1$ only where it is necessary for that sequence, clearly yields a parking strategy that minimizes $\abs{\rho}$.\\
\end{remark}

Note that the results found in this section cannot be easily generalized for $\rho\in[0,k]^n$, $k\ge2$, for multiple reasons. Firstly, unlike what we just observed for $k=1$, we cannot simply focus on $\alpha_{|_{A_i}}$. For example, $\alpha=(6,6,6,4,3,2,4)\in PP_7$ is a $2$-Naples parking function, where $\Ua=[2,4]\cup \{6\}$. However, $\alpha_{|_{\{4,5,6,7\}}}=(4,3,2,4)$ is not a $2$-Naples parking function, even though $\{4,5,6,7\}=\{j\in[7]\mid a_j\in[1,6-2]\}$.\\
Moreover, while we can generalize Theorem \ref{theosuff} for this case (that is, given $\alpha\in\pp$ and $\rho\in[0,k]^n$, the existence of certain complete parking subsequences, such that their restriction is such that all cars are able to park, implies $(\alpha,\rho)\in\pr$), the same cannot be done for Theorem \ref{theoexist}. For example, let $\alpha=(2,2,3)$ and $\rho=(0,0,2)$, we have $\Ua=\{2\}$. Although we get $(\alpha,\rho)\in PR_3$, there is no $J\subseteq[3]$ such that $\alpha_{|_J}$ is complete, and $(\alpha_{|_J},\rho_{|_J})\in PR_{\abs{J}}$. As a consequence, the results found in Section \ref{sectionknaples} for $k$-Naples parking functions, cannot be generalized for $\rho\in[0,k]^n$. Incidentally, this is the reason why, as we already remarked, the proof in Theorem \ref{char1napchoice} does not really rely on results found in previous chapters.

\addcontentsline{toc}{chapter}{Bibliography}
\printbibliography

\end{document}